\documentclass{article}
\usepackage[english]{babel}
\usepackage{todonotes}
\usepackage{listings}
\usepackage{color} %red, green, blue, yellow, cyan, magenta, black, white
\definecolor{mygreen}{RGB}{28,172,0} % color values Red, Green, Blue
\definecolor{mylilas}{RGB}{170,55,241}
\usepackage{graphics} %new
\usepackage{graphicx}
\usepackage{subcaption}
\usepackage{amsmath}
\usepackage{amssymb}
\usepackage{mathtools}
\usepackage{overpic}
\usepackage{amsthm}
\usepackage{changepage}
\usepackage{hyperref}
\usepackage{enumerate}%
\usetikzlibrary{decorations.pathreplacing,calligraphy}

\usepackage[utf8]{inputenc}

\usepackage{booktabs,array}

\usepackage[
  font = small,
  labelfont = bf,
  tableposition = top
]{caption}

\usepackage{pgfplots}
\pgfplotsset{grid style={dotted,gray}}
\pgfplotsset{compat=newest}
\pgfplotsset{plot coordinates/math parser=false}
\newlength\figureheight
\newlength\figurewidth

\usepackage{stfloats} % for positioning of figure* on the same page
\usepackage{array, multirow}

\usepackage[pass]{geometry}
\usepackage[outdir=./]{epstopdf}
\usepackage{url}
\usepackage{tabstackengine}
\setstacktabbedgap{1.5ex}%               sets gap between columns
\setstackgap{L}{1.2\normalbaselineskip}% sets baselineskip of rows
\theoremstyle{definition}
\usepackage{thmtools}

\declaretheorem[style=definition]{example}

\newtheorem{remark}{Remark}

\definecolor{KTHblue}{RGB}{25,105,188}
\definecolor{KTHlblue}{RGB}{22,159,219}
\definecolor{KTHyellow}{RGB}{251,186,0}
\definecolor{KTHred}{RGB}{176,9,48}
\definecolor{KTHlred}{RGB}{231,51,57}
\definecolor{KTHgreen}{RGB}{98,146,46}
\definecolor{KTHlgreen}{RGB}{175,202,11}
\definecolor{KTHpink}{RGB}{219,81,151}

\definecolor{mydraw}{RGB}{255,0,0}

\addto{\captionsenglish}{}
\renewcommand{\vec}[1]{\boldsymbol{#1}}
\DeclarePairedDelimiter{\ceil}{\lceil}{\rceil}
\newcommand{\vecmc}[1]{\boldsymbol{\mathcal{#1}}}

\addto{\captionsenglish}{}

\renewcommand\thmcontinues[1]{Continued}

\newcommand{\R}{\mathbb{R}}

\newgeometry{left=2.5cm,right =2.5cm,top= 2.5cm,bottom = 2.5cm}

\providecommand{\keywords}[1]{\textbf{Key words: } #1}

\usepackage[english]{nomencl}
\makenomenclature
\newcommand{\thickhline}{%
	\noalign {\ifnum 0=`}\fi \hrule height 1pt
	\futurelet \reserved@a \@xhline
}

\newcommand*{\revTwo}[1]{{#1}}%Changes from reviewer 2
\newcommand*{\revThree}[1]{{#1}}%Changes from reviewer 3
\newcommand*{\revFour}[1]{{#1}} % Changes from reviewer 4
\newcommand*{\ownchange}[1]{{#1}} %our own changes while 

\title{Accurate close interactions of Stokes spheres using lubrication-adapted image systems}
\author{Anna Broms$^{*,1)}$, Alex H. Barnett$^{2)}$ and Anna-Karin Tornberg$^{1)}$\\\\ $^{1)}$ 
	KTH Royal Institute of Technology, Department of Mathematics\\
	Stockholm, Sweden\\\\	
	$^{2)}$ Center for Computational Mathematics, Flatiron Institute,
	New York, United States\\\\
	$^*$e-mail: annabrom@kth.se}% https://www.kth.se/profile/annabrom}
\date{\today}

\begin{document}

\renewcommand\thmcontinues[1]{Continued}

	\maketitle
	\keywords{Method of fundamental solutions, method of images, elliptic PDE, Stokes flow, collocation}
	\begin{abstract}
  %Closely interacting particles moving relative to each other are affected by hard-to-resolve lubrication forces.

   Stokes flows with near-touching rigid particles induce near-singular lubrication forces under relative motion, making their accurate numerical treatment challenging. With the aim of controlling the accuracy with a computationally cheap method, we present a new technique that combines the method of fundamental solutions (MFS) with the method of images.  For rigid spheres, we propose to represent the flow using Stokeslet proxy sources on interior spheres, augmented by lines of image sources adapted to each near-contact to resolve lubrication. Source strengths are found by a least-squares solve at contact-adapted boundary collocation nodes.
   %enforcing boundary conditions through boundary collocation. 
 %Stokes flows with near-touching rigid particles  induce near-singular lubrication forces that makes an accurate numerical solution challenging. <- ab idea 	
  % EVEN OLDER VERSION: From a boundary integral perspective, two major challenges occur for particles in near-contact: lubrication forces between close to touching particles moving relative to each other are hard to resolve and near-singularities in integrals must be accurately evaluated. These two effects combined require a high resolution grid on each particle in near contact, combined with a special quadrature method, both adding to the computational expense.  
 % move to 1st sentence and expand?
 % Anna's old text:
 %In the novel singularity-free technique, sources on inner proxy-surfaces are complemented as needed by a small set of extra sources located along the lines connecting the centers of every pair of close-to-touching spheres, to resolve lubrication.
 % ab Instead: We propose representing the flow using interior spheres of stokeslets, augmented by lines of sources adapted to each near-contact, and using boundary collocation as in the MFS. 
 We include extensive numerical tests, and validate against reference solutions from a well-resolved boundary integral formulation.
% using a special quadrature method.  
With less than 60  additional % = 3*n_im  
  image sources per particle per contact, % via single-body preconditioning, 
 we show controlled uniform accuracy to three relative digits in surface velocities,
 and up to five digits in particle forces and torques,
 %for all physically relevant separations between the particles, i.e.~down to gaps of a 
 for all separations
 %between the particles
 down to a thousandth of the radius.
 In the special case of flows around fixed particles, the proxy sphere alone gives controlled accuracy. A one-body preconditioning strategy allows acceleration with the fast multipole method, hence close to linear scaling in the number of particles. This is demonstrated by solving problems of up to 2000 spheres on a workstation using only 700
 %unknown
 proxy sources per particle. \\	
		
	\end{abstract}

	\section{Introduction}\label{sec:intro}
	The Stokes equations constitute a linearization of the Navier-Stokes equations, valid in the limit of small Reynolds numbers, where viscous forces are dominant over inertial forces. Consequently, the fluid adjusts instantly to any applied force without damping, as the mass of the fluid does not affect the fluid motion. Such fluid systems with immersed particles have a wide range of applications, from the design of a lab-on-a-chip for diagnostics or drug delivery, to water purification and sedimentation processes. In biology, the flow of red blood cells in our veins, the motion of swimming bacteria \cite{Stone1996, Cortez2005,Delmotte2015} and the mechanics of the cytoskeleton in our cells, where a complex network of microtubules, actin filaments and molecular motors provide structural support and enable movement within cells \cite{Maxian2022, Yan2022}, can be understood based on a Stokesian model \cite{Pral2022, Lu2019}. In materials science, the flow of fibers in e.g.~the paper making industry can be characterized \cite{Niskanen}, and novel materials developed \cite{Hakansson2014,Calabrese2023}. 
	
	We focus on the Stokes resistance boundary value problem (BVP) for a suspension of rigid spheres immersed in a given background flow $\vec u_\infty$. Let $\Omega \coloneqq \bigcup_{k=1}^P \Omega^{(k)}$ denote a collection of $P$ balls, and $\partial \Omega$ denote the union of all of their sphere boundaries $\partial \Omega^{(k)}$.
 %The translational and angular velocities are given for each of these spheres. 
 We assume that $\vec u_\infty$ is a homogeneous Stokes solution throughout $\R^3$, and that in the presence of the spheres, the flow field becomes
 $\vec u_\infty + \vec u$. The unknown disturbance flow field, $\vec u$, therefore solves the BVP
%presence of the particles modifies the flow velocity $\vec u(\vec x)\in\mathbb R^3$ at each location $\vec x$ and the BVP for the disturbance field caused by the particles in the fluid domain exterior to the spherical bodies is
	\begin{equation}\label{stokeseq}
		\begin{aligned}
			-\mu\Delta \vec u + \nabla p &= \vec 0, && \text{ in }\mathbb R^3\backslash\overline{\Omega},\\
			\nabla\cdot \vec u &= 0,&&\text{ in }\mathbb R^3\backslash\overline{\Omega},\\
			\vec u &= \vec u_{\text{bc}},  &&\text{ on }\partial\Omega, \\
   \vec u(\vec x) &\to \vec 0, &&\|\vec x\|\to\infty.
		\end{aligned}
	\end{equation} 
	Here, $p$ is the disturbance pressure field, and $\mu$ the constant viscosity of the fluid. %\question{Introduce equation earlier?}
 %With $\vec u^{\text{bc}}$ decaying at infinity and
We assume no-slip boundary conditions, with given translational and angular velocities for each sphere. For sphere $k$, this means that the velocity boundary data is prescribed as
	\begin{equation}\label{bc}
		\vec u_{\text{bc}}^{(k)}(\vec x) = \vec v^{(k)}+\boldsymbol\omega^{(k)}\times(\vec x-\vec c^{(k)})-\vec u_{\infty}(\vec x), \quad %\forall
  \vec x\in \partial\Omega^{(k)},
	\end{equation}
	with $\vec v^{(k)}$ its given translational velocity, $\boldsymbol \omega^{(k)}$ its given angular velocity, and $\vec c^{(k)}$ its center. Aside from solving for the above flow field, it is also of interest to extract the net force and torque on each particle.
 
 As particles move relative to each other in Stokes flow, their interactions become increasingly difficult to resolve accurately with decreasing particle separations, due to growing lubrication forces \cite{Kim1991}. Small gaps and strong lubrication come with high resolution requirements and a high computational cost. This is intractable for simulations of systems with many particles or of systems with small Brownian particles, where many realizations are needed to draw statistical conclusions. With simulation methods based on a volumetric grid, such as finite element methods, the narrow gap between particles becomes very costly to mesh. For a boundary integral method (BIE), on the other hand, it is enough to discretize only the particle surfaces, but two other challenges have to carefully addressed: As one solves for a so called layer-density\footnote{The simplest interpretation of the layer density is that of a force density over the surface, which holds if a single layer formulation is used in the boundary integral formulation.} over the particle surfaces, near-singularities in integrals over the surfaces must be accurately evaluated \cite{Ying06,AfKlinteberg2016,Bagge2021,Wu2021b,Zhu2022}. %\note{The later refs not specifically Stokes... There are other ones that might be better to mention and/or the original QBX paper?}
 Additionally, strong lubrication forces manifest themselves as a sharply peaked layer density in regions close to touching with another particle, which, to be resolved, requires a very fine surface grid,  resulting in large linear systems with many unknowns to solve.  There is hence a need for more efficient techniques. 

% ab: I notice this Intro is directed at JFM fluids readers rather than numerical PDE (applied math).

 Several approximate methods exist to capture lubrication \cite{Sangani1994,Sangani1996}. This is often done by explicitly adding a contribution to particle pair-interactions to compensate for what is not captured with the ``regular'' computational method or grid \cite{Mammoli2006,Lefebvre-Lepot2015,Orsi2023}. The most famous of these approximate correction strategies is the widely used Stokesian dynamics \cite{Brady1988,Brady2001,Swan2011,Fiore2019}. The surfaces of the particles are then not explicitly represented and the interaction between particles is governed by a multipole expansion.  The corrected hydrodynamic resistances between particles are assumed to be pair-wise additive, which is an approximation that holds only for small particle gaps 
 %and hence, the resulting computed interactions are deteriorating in accuracy compared to a reference based on the method of %reflections 
 \cite{Wilson2013}. Including more terms in the multipole expansion can provide higher accuracy also for moderately separated particles, but implies an increased algebraic complexity. Stokesian dynamics can be accelerated with fast summation techniques both with and without Brownian motion, which enable simulations of large particle systems over long time-spans to a manageable cost \cite{Fiore2019}. An alternative is the rigid multiblob method \cite{Usabiaga2016,Broms2022}, where a particle is described as a collection of spherical blobs or beads constrained to move as a rigid body and the interaction between blobs is regularized by the Rotne-Prager-Yamakawa tensor. In \cite{Broms2022}, we improve on the accuracy of the multiblob method in the far-field and also reduce errors for close interactions. Neither with Stokesian dynamics nor the rigid multiblob method is the accuracy however possible to \emph{control}. Using a BIE method equipped with special quadrature, controllable accuracy can be obtained for moderately separated particles \cite{AfKlinteberg2016,Bagge2021}. One such technique based on spherical harmonics enables fast large-scale simulations for spheres specifically \cite{Corona2018,Yan2020}. However, since they are not spatially adaptive, such methods demand an excessive number of surface unknowns at small separation distances, and thus become inefficient.
 %However, the peaked layer density hinders controllable accuracy for very dense suspensions if the particle resolution is not intractably fine, and either a spatially adaptive handling of the mesh or an adaptively set degree of expansion for the particle would be required for these methods to be applicable for small separation distances.    % alex simplified.

In the description above, it has been assumed that closely interacting particles also move relative to each other; however, only under relative motion do lubrication forces diverge with vanishing separation distance between particles \cite[p.~175]{Kim1991}. For bodies at fixed locations subject to a background flow (e.g., a porosity simulation), the singularities are much weaker, so that a non-adaptive BIE method is adequate \cite{AfKlinteberg2016}.
The difference in difficulty to resolve problems with or without relative motion
has been exploited numerically. In \cite{Lefebvre-Lepot2015}, velocities for pairs of particles are split into relative and rigid body motion. Only relative motions then require special treatment as a consequence of a peaked force density over the particle surfaces, while joint rigid body motions are easy to resolve \cite{Lefebvre-Lepot2015}. In \cite{Orsi2023}, lubrication-corrections for particle pairs, similar to those in Stokesian dynamics, are constructed to be frame-invariant, meaning that only relative particle motions are affected by the correction.  

%\subsection{Image enhancement}\label{reflections} 
 %%%%%%%%%%%%%%%%%%%%%%%%%%%%%%%%%%%%%%%%%%%%%%%%%%
%\subsection{Sphere reflections}
%\begin{remark}[Image enhancement]
When formulating an integral equation, a free space Green's function is typically used, and boundary conditions are enforced by solving for the layer density in the solution ansatz. If one can instead build a Green's function that satisfies the boundary conditions, the solution is easily constructed. Such domain specific Green's functions can sometimes be found using an {\em image} or {\em reflection} technique, the simplest case being that of a homogeneous Dirichlet condition on a flat plane for the Laplace equation. %\cite{jackson}.
%Here, the boundary condition will hold if, for each original source, an image source with opposite sign is placed at the same distance to the plane, but on the opposite side. 
Such a Laplace image system with known strengths and locations can also be built for two interacting conducting discs or spheres, but here the number of reflections, or image sources, is infinite
(this problem was solved 170 years ago by Lord Kelvin \cite{kelvin2spheres}; also see
\cite[\S172--\S173]{maxwell1} and \cite{Cheng1998,Cheng2000}).
With more objects, there is an exponential explosion of multiple reflections---even for three disks or spheres the image points form a fractal \cite{Mumford}---rendering a numerical image sum impossible. %fractal dust on a circle that cuts the spheres orthogonally in addition to image lines... \cite{Klein}.

Cheng \& Greengard \cite{Cheng1998} considered closely spaced conducting discs (the Laplace equation) in two dimensions, followed by an extension to three dimensions for conducting spheres by Cheng \cite{Cheng2000}, and dielectric spheres by Gan et al.~\cite{Gan2016}. Here, interactions of the objects get increasingly difficult to resolve with decreasing separation distances: a standard numerical solution in a multipole (Rayleigh expansion) basis centered on each object then demands
a large and wasteful number of unknowns per body.
% when material properties (conductivities) of the objects differ significantly from that of the background. 
The authors introduce a hybrid method which adds to each multipole basis function its (truncated) pairwise image reflection series
in each sufficiently close neighbor,
%for each pair of objects
thus keeping the number of unknowns per body small regardless of the closeness of separation.
%resolution requirements when solving the integral equation. This allows for a solution to the problem to high accuracy also at very small separation distances.  
%The formulae are already complicated for Laplace spheres, and since analytic image formulae become significantly more complicated for the Stokes case, a generalization of this hybrid method seems out of reach. 
Yet, image methods are significantly more complicated for the Stokes case.
For a flat wall with a no-slip condition, analytic formulae were first derived for the Stokeslet by Blake \cite{Blake1971}, and later in a different form by Gimbutas et al.~\cite{Gimbutas2015}. 
The methodology for deriving the image system for a Stokeslet in the exterior of one sphere, that enforces a no-slip condition on the surface of the sphere, was introduced in \cite{Fuentes1988,Fuentes1989}, and expressions are found in \cite{Maul1996} and \cite[\S10.2]{Kim1991}.  They involve multiple types of singularities: the Stokeslet, rotlet, stresslet, and potential dipole,
including {\em continuous distributions} along the line from image point to center. Although various approximations exist for two distant spheres \cite[Ch.~10-11]{Kim1991}, we do not believe that an exact formula for multiple reflections in two nearby spheres would be practical.
Thus a Stokes generalization of the above hybrid method seems out of reach.
 However, we will sacrifice analytic formulae in favor of a more numerical approach that still draws inspiration from Lord Kelvin, Cheng \& Greengard,
 and the Stokes image for a single sphere.

%% At end of section: 
The approach of this paper is, broadly speaking, a method of fundamental solutions (MFS). The solution flow $\vec u$ is written as a linear combination of fundamental solutions to the Stokes PDE, located in the non-physical region (the sphere interiors), whose coefficients are found by enforcing the boundary data \revThree{in the least-squares sense at a similar number of} surface collocation points.
This idea has a long history, and sometimes also goes under the names of the charge simulation method, the method of auxiliary sources, or proxy-sources \cite{Alves2004,Fairweather2005,Barnett2007,Alves2009,Liu2016,Karageorghis2019,Antunes2022}.  %Chen2014
A choice of locations and collocation points leading to small errors depends on the domain and data, is heuristic, and is not always possible to determine.
However, there are rigorous convergence results in certain smooth geometries, such as the disk with source locations on a concentric circle \cite{Ka89}, and when the solution may be analytically continued as a PDE solution up to the source curve \cite{Stein2022}.
The MFS has the disadvantage of \ownchange{involving} ill-conditioned rectangular linear systems, when compared against the well-conditioned square ones arising from discretizing 2nd-kind Fredholm
%boundary integral equations.
BIEs. In special geometries such as close-to-touching spheres, however, we will see that the advantage of much smaller numbers of unknowns at a given accuracy makes the MFS rather favorable compared to BIEs. %\ownchange{Since the preprint of this paper was submitted, another paper on Stokes MFS appeared \cite{Jordan2025}, who approximate particle-particle interactions via Gauss-Siedel iterations instead of solving the ill-conditioned linear system. 
%}.

Two sources of numerical errors affect a dynamic simulation of particles in a non-Brownian Stokesian suspension: the error in the Stokes BVP, and the time-step discretization error in evolving the particle coordinates. To properly capture the strong lubrication forces in a dynamic simulation of close-to-touching particles, the time-step size often has to be extremely small, which adds substantially to the computational cost. Larger time steps reduce the computational effort but can lead to particle overlaps.
In addition, at very small separation distances, the surface roughness of physical particles becomes non-negligible, and the Stokes equations are no longer a good model for their interactions \cite{Ball1995,Morris2020}.
To avoid investing a huge computational effort in this regime, a contact-avoiding strategy is needed, such that a minimum particle separation is enforced \cite{Lu2019,Bystricky2019,Yan2020,Broms2024}. In the present work, focus is on  solving the Stokes BVP,
%at any given instance of time
and time-stepping is left for future work. However, with any future time-stepping strategy in mind, its resolution requirements determine what we here consider a \emph{physically relevant} particle separation for which we seek accurate solutions.

\begin{figure}[h!]
		\centering
		\begin{subfigure}[t]{0.49\textwidth}
			\centering
   %\hspace*{-4ex}
			\includegraphics[trim={0.5cm 0cm 0cm 0.6cm},clip,width=1\textwidth]{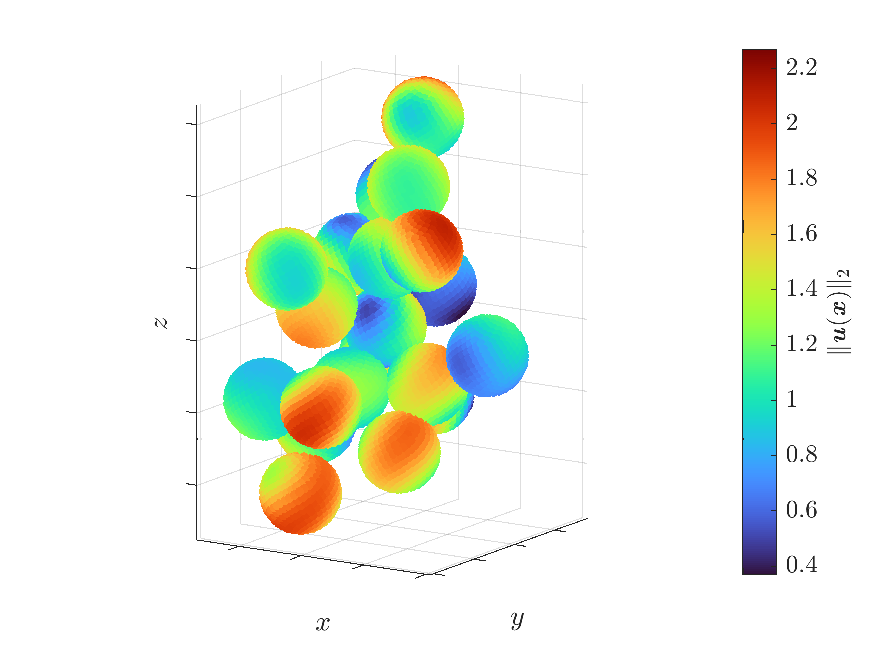}
            \caption{Surface velocities. }		
			\label{surf_vel}
		\end{subfigure}~~
  		\begin{subfigure}[t]{0.49\textwidth}
			\centering
			\includegraphics[trim={0.5cm 0cm 0cm 0.55cm},clip,width=1.03\textwidth]{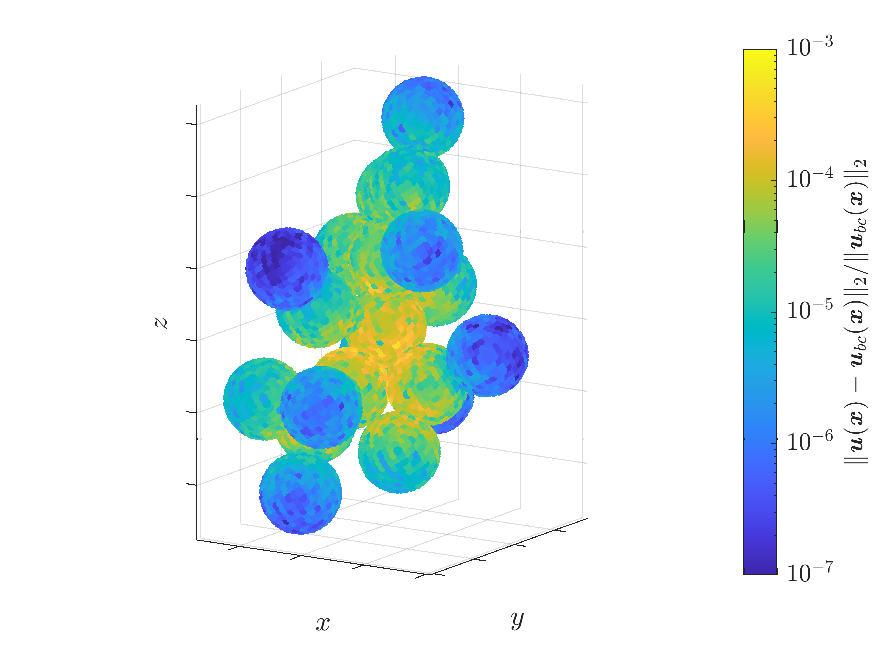}
  % \tikz \fill [orange] (0.1,0.1) rectangle (0.2,0.2);
			\caption{Relative residual at the particle surfaces. }
			\label{cluster_resim}
		\end{subfigure}
  \caption{By using an adaptive number of additional source points per closely interacting particle pair, lubrication forces between 20 closely interacting spheres in a cluster (many separated by $10^{-3}$ radii) are resolved with our new technique. Each sphere travels with a random translational and angular velocity. The target accuracy of $10^{-3}$ in the relative residual is reached everywhere on the particle surfaces.}
  \label{image_example}
  \end{figure}

%%%%%%%%%%%%%%%%%%%%%%%%%%%%%%%%%%%
\subsection{Contribution and outline}
 We present an
 %strategy
 efficient numerical scheme
% for rigid particles in Stokes flow
for the Stokes resistance BVP compatible with linear-scaling
fast solution algorithms.
 % in order to control the errors in solving the Stokes BVP.
 We target 3-digit (0.1\%) uniform error level in the surface velocities, chosen to be small enough that other sources of error discussed above would dominate in practical settings. Moreover, we will show that controllable accuracy can be obtained at a competitive cost for %\emph{any physically relevant} particle separation, here meaning 
 any separation down to $10^{-3}$ radii,
far smaller than has been shown in other multi-sphere Stokes solvers \cite{AfKlinteberg2016,Corona2017,Fiore2019,Yan2020,Sprinkle2020,HU2001}. %the last pretty old...
 %which is assumed the smallest distance allowed by the model %\footnote{Smaller particle separations than those considered to be physically relevant here are not allowed by any contact-avoiding strategy in a dynamic simulation.}. 
 Our cost is higher than that of more approximate techniques such as Stokesian dynamics or the rigid multiblob method, but lower than that of a BIE method accurate to the same level.
 Our technique is a variant of the MFS discussed above,
 hence avoids volume discretization, and allows evaluation of
 potentials right up to the boundary without specialized techniques.
 %Our technique is a variant of the method of fundamental solutions (MFS): the idea is to express the flow field as a sum of Stokes fundamental solutions located inside each body, then enforce boundary conditions by collocation on a set of surface nodes. Like BIE, this is a boundary-based method that avoids volumetric discretization.
%However, since the resulting MFS overdetermined linear system is very ill-conditioned,
To tackle ill-conditioning,
we apply so-called {\em one-body preconditioning}
\cite{Liu2016,Stein2022}
to create a well-conditioned square system amenable to FMM-accelerated iterative solution in quasi-linear time. This allows large numbers of spheres to be handled on a single workstation.
The MFS source locations, types, and collocation nodes must be carefully chosen to capture the near singular force densities involved. One of our main contributions is a recipe for this---inspired by the analytic image distributions discussed above---that adapts to each pairwise particle separation.
We propose Stokeslets located on interior concentric ``proxy'' spheres, complemented by
point sources along segments of lines connecting the centers of close-to-touching pairs. We call the latter \emph{image points}, as their locations discretize an image line source.
See Figure \ref{colloc_points}.
The resulting number of unknowns needed to achieve the target accuracy is much less than for a BIE-method, as we show in our experiments. The result is a higher efficiency scheme.
  
A small-scale example of the type of problem we solve is visualized in Figure \ref{image_example} with  a cluster of 20 closely interacting unit spheres with relative motion. It is a challenging problem: each sphere is $\delta = 10^{-3}$ from at least one neighbor, with velocities set to randomly sampled rigid body motion, as visualized in panel \ref{surf_vel}. With the number of source points per particle set adaptively based on separation distances, panel \ref{cluster_resim} shows that the target accuracy of $10^{-3}$ is reached everywhere on the particle surfaces.

\ownchange{
\begin{remark}[Recent Stokes MFS work]
    Since the preprint of this paper appeared on arXiv, %it influenced 
    another MFS Stokes paper appeared by Jordan \& Lockerby \cite{Jordan2025}. They present results for the resistance and mobility problems, including nonspherical particles. However, in contrast to our work, only problems with weak lubrication effects (large separations) are treated. %they do not consider space-filling suspensions of particles and as , the scheme is formally non-convergent. 
    They apply Gauss--Seidel iterations to account for particle-particle effects, giving a $\mathcal O(P^2)$ cost, contrasting the linear scaling obtained via the FMM in our work. Results in \cite{Jordan2025} are shown for particles along a line, with interactions greatly simplified beyond a certain distance.
    In contrast, we present a method with %without such %simplifications and
    full interactions and %convergent method and 
    test it with space-filling configurations.  
 %whose long-range quadratic cost would be huge.
 \end{remark}
}

%As there are several parallels to draw to a BIE method, we will in the following point out some key differences and similarities. 
 
 % \subsection{Paper overview}
The outline of the paper is as follows. The Stokes MFS is
%mathematically
introduced in Section \ref{MFS}, including the basic discretization, the one-body preconditioning scheme, and the enhanced MFS image sources needed for problems with strong lubrication forces.
The parameter selection for the basic MFS discretization is discussed in Section \ref{param_sec}, and the resulting performance for large-scale problems with small lubrication effects is discussed in Section \ref{large_ex}. 
%This category includes dilute particle suspensions, where each particle is sufficiently far away from its neighbors, and fixed-in-space closely interacting particles not moving relative to each other, only affected by a background flow.  
For the more challenging problems involving strong lubrication forces, fundamental solutions at the additional source points have to be chosen such that their combined approximation power is sufficient for resolving the solution to the set of PDEs. The types of fundamental solutions needed at the image points are investigated in Section \ref{lub}, together with distributions of additional collocation points to capture the near-singular lubrication flow. The investigations lead to our adaptive scheme for the number of image points, depending on the particle-particle separations only, which is then applied to problems with clusters of closely interacting spheres affected by lubrication in Section \ref{Num_res}. 
To demonstrate that a target accuracy of $10^{-3}$ can be reached for all physically relevant particle separations, errors in the flow fields, forces and torques computed with MFS are assessed by comparing to reference solutions computed with a BIE-method equipped with special quadrature, by self-convergence tests, and by measuring the residual at the particle surfaces at points different from the collocation points. 
%To demonstrate that a target accuracy of $10^{-3}$ can be reached for all physically relevant particle separations, the flow fields, forces and torques computed with MFS are compared to reference solutions computed with a BIE method equipped with special quadrature. 
%As reference solutions are very costly to compute for the smallest particle gaps, self-convergence tests are also made. In addition, and as our most common test, we check the relative residual at the particle surfaces. %, that is 
Our findings are summarized in Section \ref{conclusions}.

 %%%%%%%%%%%%%%%%%%%%%%%%%%%
	\section{The method of fundamental solutions}\label{MFS}
In this section, we first introduce the basic discretization that will be used in the paper, followed by a description of the preconditioning strategy that we employ. In subsection \ref{image_enhance}, we discuss how to enhance this basic discretization to allow for accurate solutions also in the case of strong lubrication forces, i.e.~for particle interactions with relative motion at small separation distances.

 \subsection{The basic discretization for one particle}\label{base_disc}
%We next introduce the method of fundamental solutions in its most basic setting, suitable for problems where lubrication effects are small to moderate. 

Possibly the simplest solution to the PDE and decay condition at infinity in \eqref{stokeseq} is given by
	\begin{equation}\label{stokes_flow}
		\vec u(\vec x) = \mathbb S(\vec x-\vec y)\vec f,
	\end{equation}
 valid in $\mathbb R^3\backslash \{\vec y\}$.
	The tensor $\mathbb S$ is the Stokes fundamental solution known as the Stokeslet, and $\vec f$ is a point force located at $\vec y$. In free-space, without confinements, defining $\vec r\coloneqq\vec x-\vec y$, $\vec I \in \mathbb R^{3\times 3}$ the identity matrix, and $\left(\vec r\vec r^T\right)_{ij} = r_ir_j\in \mathbb R^{3\times 3}$ the dyadic product, the Stokeslet is
	\begin{equation}\label{stokeslet}\index{Stokeslet}
		\mathbb S(\vec r) = \frac{1}{8\pi\mu \|\vec r\|}\left(\vec I+\frac{\vec r \vec r^T}{ \|\vec r\|^2}\right).
	\end{equation}
 There is of course a pressure solution $p$ accompanying $\vec u$ in \eqref{stokes_flow}, but we will not need it. 
 To represent the boundary condition  $\vec u_{\text{bc}}$ on a single particle, however, a linear combination of fundamental solutions at $N$ source points $\vec y_j$ in the particle interior is needed, so that
	\begin{equation}\label{sum}
		\vec u(\vec x) = \sum_{j=1}^N\mathbb S(\vec x-\vec y_j)\vec \lambda_j.
	\end{equation}
   %To write the solution field as a sum of the form in \eqref{sum} is the basic idea behind the MFS.
   %With this representation \eqref{sum},
   Without loss of generality we fix a sphere radius $R=1$.
   One imposes boundary conditions via collocation at a set of points $\vec x_i$, $i=1,\dots,M$,
   covering the surface, where $M>N$.
   This gives the overdetermined linear system
\begin{equation}\label{lsq}
\sum_{j=1}^N \mathbb S(\vec x_i-\vec y_j)\vec \lambda_j
\; = \;
\vec u_{\text{bc}}(\vec x_i),   \qquad i=1,\dots, M,
\end{equation}
   which is solved in the least-squares sense to find the set of
   strength vectors $\vec \lambda_{j}$.
   %are then determined as a solution to the least-squares problem,
   %\begin{equation}\label{lsq}
	%\min\limits_{\vec\lambda}\sum_{i=1}^M\left(\sum_{j=1}^N\mathbb G(\vec x^i-\vec y^j)\vec \lambda^j -\vec u_{\text{bc}}(\vec x^i)\right)^2.
%	\end{equation}
 In the basic version of our algorithm, the $N$ source points $\vec y_j$ cover a concentric interior proxy-sphere of radius $R_p<1$.
 % suggest don't mention n_p until absolutely have to :)
 As source and collocation points live on separated surfaces, the singularities in the Stokeslet at $\vec x = \vec y_j$, and thus the singular quadratures needed for BIE methods, are avoided. The tall rectangular matrix
 $\vec G$ of size $3M\times 3N$,
 whose $3\times 3$ block elements are $\vec G_{ij}\coloneqq\mathbb S(\vec x_i-\vec y_j)$,
 we refer to as the target-from-source matrix.
 %We have observed its ill-conditioning to be worsened by a nonuniform point distribution: clustered source points lead to nearly linearly dependent columns, while clustering in collocation points often results in less information being obtained from the boundary data sampled on the particle surface.   % alex simplified into the start of next paragraph...
 
 %the distance between the proxy-surface and the surface where boundary conditions are imposed 
	%the %ill-conditioning
 %dynamic range
 %of $\vec G$ is also more severe the more clustered the source and collocation points are on their respective surfaces. A uniform discretization of points is therefore preferred.
 % We thus seek as uniform as possible a distribution on
 % surface and proxy spheres, and have found that
 % {\em spherical design} nodes perform as least as well as any other
 % distribution.
 % \footnote{Spherical design nodes are available in double-precision accuracy for $N$ up to around 16000 at
	% 	\url{https://web.maths.unsw.edu.au/~rsw/Sphere/EffSphDes/}.} \cite{sphdesign}. 
  % (We have compared this to other sphere point distributions such as Fibonacci grids, but not found any distribution that exceeds spherical designs in efficiency.) 

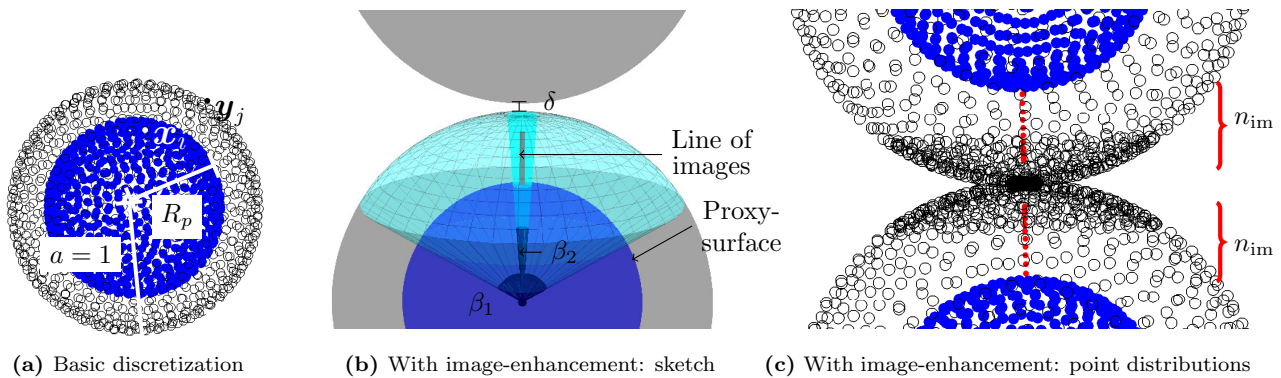
\begin{figure}[h!]
  \begin{subfigure}[t]{0.2\textwidth}
% \vspace*{-20.5ex}
 
		\centering
  %\vspace*{-5.5ex}
  \begin{tikzpicture}
    \hspace*{-23ex}  
  \includegraphics[trim = {4.4cm 2cm 3.5cm 2cm},clip,width=1.1\textwidth]{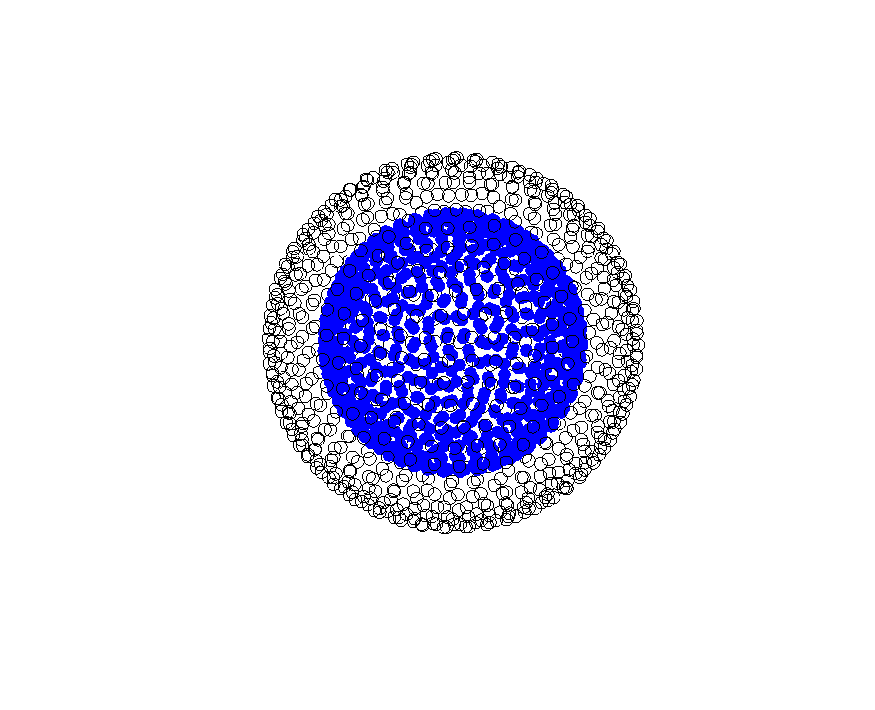}

  \coordinate (start) at (-2.0,2.32);
			\coordinate (target) at (-0.85,2.8);
			%
			%% Draw the arrow
			\draw[|-|,white, line width = 1.5pt] (start) to (target);
   \coordinate (target2) at (-1.7,-0.8);
   \draw[|-|,white, line width = 1.5pt] (start) to (target2);
     \coordinate (start2) at (-1.7,2.1);
			\node[right,fill=white] at (start2) {{$R_p$}};
\tikz \fill [white] (-1.5,0.5) rectangle (-1.6,0.6);
\coordinate (source) at (-0.85,3.55);
\node at (source) {\huge{$\cdot$}\large{$\vec y_j$}};
\coordinate (target) at (-1.65,3.20);
\node[white] at (target) {\huge{$\cdot$}\Large{$\vec x_i$}};
 \coordinate (start3) at (-2.2,1.6);
			\node[left,fill=white] at (start3) {{$a = 1$}};

    \end{tikzpicture}

		% \includegraphics[trim = {4cm 1cm 3cm 2cm},clip,width=1.15\textwidth]{figures/sphere.eps}
  %\vspace*{5.9ex}
		\caption{Basic discretization}
		\label{basic_disc}
\end{subfigure}~~~~
%\hspace*{1ex}
       \begin{subfigure}[t]{0.39\textwidth}
   %    \vspace*{-27ex}
		\begin{tikzpicture}		
			% Draw a circle using TikZ
			%\draw (0,0) circle (1);
			%
			%
			%% Label the start and target points
			%\node[below] at (start) {Start};
			%\node[above] at (target) {Target}
			%
			% Insert an external graphic (e.g., a PNG file)
			\node at (3,0) {\includegraphics[trim={21cm 8cm 19cm 10cm},clip,width=0.95\textwidth]{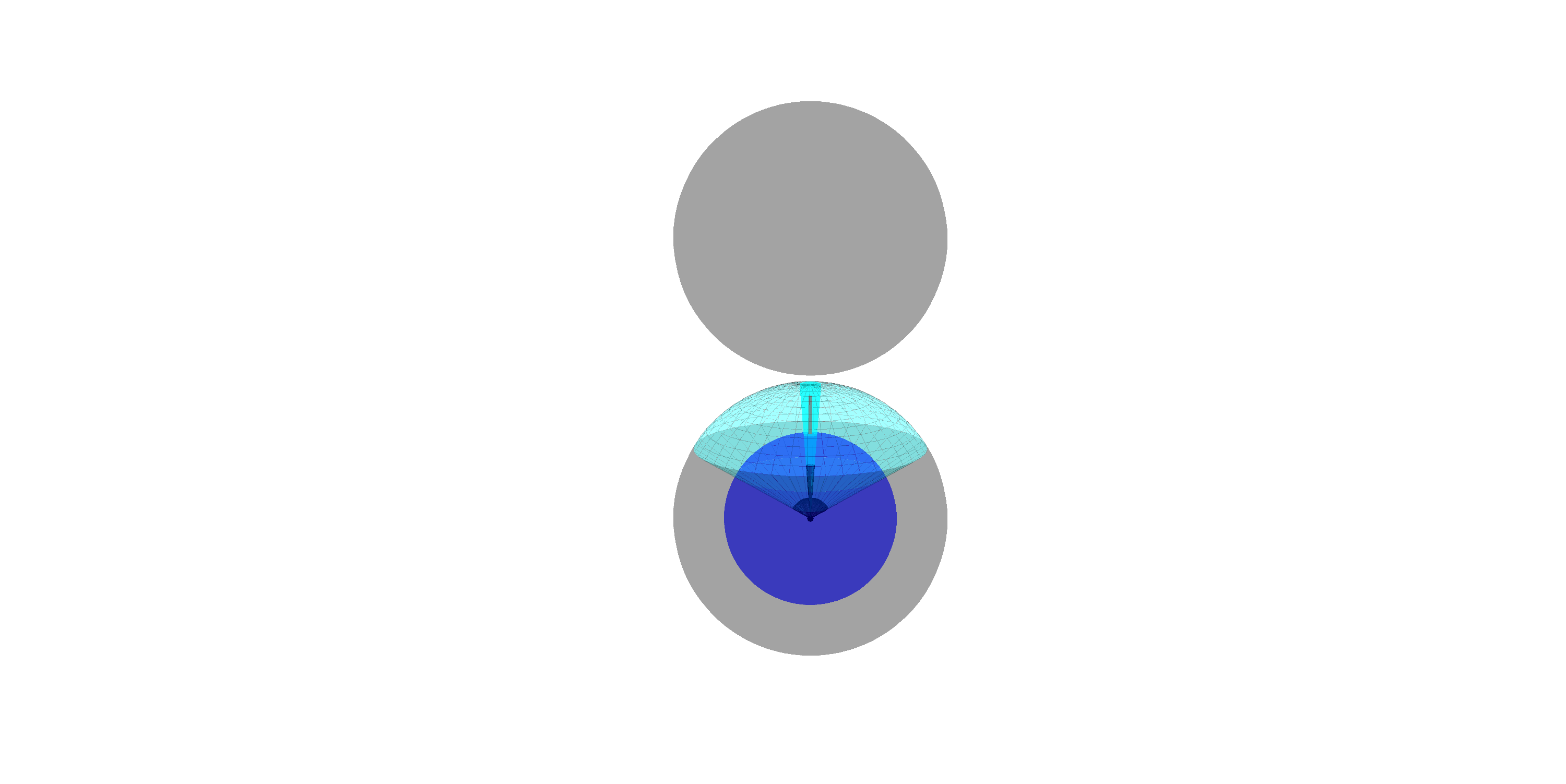}};
			\coordinate (start) at (4.8,0.2);
			\coordinate (target) at (2.9,0.2);
			%
			%% Draw the arrow
			\draw[->] (start) to (target);
			\node[right] at (start) {\shortstack{Line of\\images}};%Line of\\ images};
   		\coordinate (start) at (3.2,-1.1);
			\coordinate (target) at (2.9,-1.1);
   \node[right] at (start) {$\beta_2$};
			\draw[->] (start) to (target);
   \coordinate (start) at (2.1,-1.8);
     \node[right] at (start) {$\beta_1$};
     \coordinate (start) at (2.9,0.77);
			\coordinate (target) at (2.9,0.91);
   \draw[|-|] (start) to (target);
    \coordinate (start) at (3.1,0.92);
   \node[right] at (start) {$\delta$};
   			\coordinate (start) at (5.2,-0.75);
			\coordinate (target) at (4.4,-1.2);
			%
			%% Draw the arrow
			\draw[->] (start) to (target);
			\node[right] at (start) {\shortstack{Proxy- \\surface}};
		\end{tikzpicture}
  \caption{With image-enhancement: sketch}
  \end{subfigure}
  \hspace*{-2ex}
  \begin{subfigure}[t]{0.39\textwidth}
\hspace*{1ex}
\begin{tikzpicture}
			%Draw a circle using TikZ
			%\draw (0,0) circle (1);

			% Label the start and target points
			%\node[below] at (start) {Start};
			%\node[above] at (target) {Target}
			
			%Insert an external graphic (e.g., a PNG file)
			\node at (40,4) {\includegraphics[trim={3.2cm 3cm 2.5cm 2cm},clip,width=0.98\textwidth]{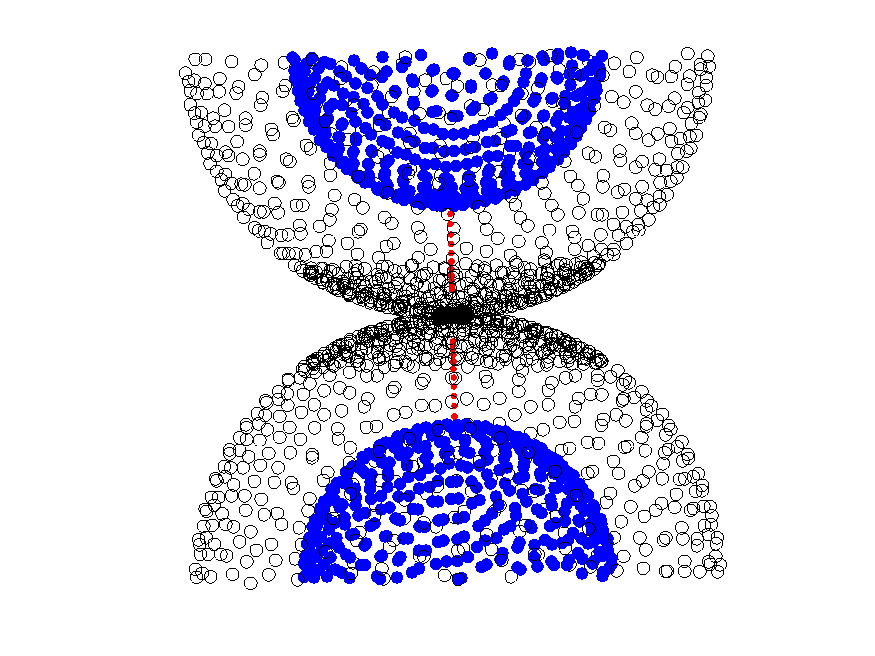}};
   \coordinate (start) at (43.0,4.65);
			%\coordinate (target) at (40.4,-1.1);
			%
			%% Draw the arrow
			%\draw[->] (start) to (target);
			%\node[right] at (start) {Proxy-surface};
     \draw [pen colour=mydraw,decorate,
      decoration = {brace,mirror}, ultra thick,mydraw] (42.5,4) --  (42.5,5.15);
       \draw [pen colour=mydraw,decorate,
      decoration = {brace,mirror}, ultra thick,mydraw] (42.5,2.5) --  (42.5,3.55);
      %     \draw [pen colour=mydraw,decorate,
     % decoration = {brace,mirror}, ultra thick,mydraw] (40.8,4) --  (40.8,5.2);
     %     \draw [pen colour=black,decorate,
     % decoration = {calligraphic brace,mirror}, ultra thick] (40.7,4) --  (40.7,5.15);
      \node [fill=white,inner sep=1pt] at (start) {$n_{\text{im}}$};
       \coordinate (start2) at (43.0,2.95);
        \node [fill=white,inner sep=1pt] at (start2) {$n_{\text{im}}$};
      %\coordinate (line) at (2.95,0.92);
    %  \node[right] at (line) {Line of images};
      
    %\draw [red,decorate,
   % decoration = {calligraphic brace,mirror},ultra thick,red] (50,1) --  (50,2);
  \end{tikzpicture}
 \caption{With image-enhancement: point distributions}
 \label{colloc_images}
 \end{subfigure}
 % \caption{Sketch}
 % \end{minipage}
	%	\caption}		
\caption{Panel (a) displays the basic discretization, with the original source points $\vec y_j$, $j = 1,\dots,N$, on the radius $R_p$ proxy-surface in blue, and collocation points  $\vec x_i$, $i = 1,\dots,M$, in black on the unit sphere surface. Both source and collocation points are obtained from spherical design nodes. Panels (b) and (c) display the image-enhanced grid of source and collocation points for a pair of close-to-touching spheres. Additional image points are distributed along a line segment in the interior of each sphere, while additional collocation points are distributed  \revFour{``above'' each line of images. The collocation points are located on two spherical caps determined by angles $\beta_1$ and $\beta_2$ surrounding the point of closest approach to another particle, as in panel (b).} If a particle is closely interacting with more than one neighbor, one region of extra source and collocation points of this type is added for every near-contact. In panel (c), collocation points are indicated in black and image points in red,  generated via \eqref{line_points}-\eqref{tj}. }
\label{colloc_points}
\end{figure}

  In order to minimize the boundary condition residual and the ill-conditioning,
  we seek as uniform as possible a point distribution on the physical sphere and proxy sphere, respectively. The many available quasi-uniform point distributions include Fibonacci grids \cite{fibogrid06,Marques2021}, HEALpix grids \cite{healpix},
 quadratures with a (possibly different) equispaced grid on each line of latitude,
 %\cite[Sec.~5]{Stein2022},
and Gaussian quadrature schemes for spherical harmonics
such as Lebedev points (based on the octahedron)
or Ahrens--Beylkin points (based on the icosahedron);
see \cite{ahrens09} and references within.
However, we recommend so-called ``spherical design points'' \cite{sphdesign}
 for surface- and proxy-points,
  since with {\em constant} weights they
integrate spherical harmonics up to a certain degree $p$ exactly,
even though their efficiency (they need $N\approx p^2/2$ points)
is not quite as high as Gaussian ($N\approx p^2/3$).
The constancy of weights means that they are quasi-uniform, unlike variable-weight quadrature schemes which are more clustered; we find that this improves stability.\footnote{Note that certain of the above grids could enable
exploiting azimuthal symmetry, but for the basic proxy discretization the SVD
costs are too small for it to matter.}  Spherical designs only exist for certain $N$ values.
  We use points computed\footnote{These are available in double-precision accuracy for $N$ up to around 16000 at
\url{https://web.maths.unsw.edu.au/~rsw/Sphere/EffSphDes/}.}
by Womersley \cite{womersley18}.
In our tests, proxy-surfaces of radii $R_p$ are each discretized with $N$ spherical design nodes, while collocation surfaces are discretized with $M\approx1.2N$ spherical design nodes, as in Figure \ref{basic_disc}.

% missing paragraph indentation due to color macro...
\revThree{The ratio $M/N\approx 1.2$ was found by a convergence study, and is typical for MFS applications \cite{Barnett2007,Stein2022}. We do not know of rigorous analysis of this ratio, apart from the case $M/N=1$ for the 2D Laplace equation \cite{Ka89}. Intuitively, a ratio greater than $1$ is beneficial because the potential generated by somewhat nearby proxy points contains frequency components higher than merely the Nyquist frequency associated with their spacing.}

\begin{remark}[Conditioning]     % rrrrrrrrrrrrrrrrrrrrrrrrrrrrrrrrrrrrrrrrrrrrrrrrrrrrrrrrr
\label{r:cond}
The condition number of the target-from-source matrix for a single sphere
is expected to grow exponentially like $(1/R_p)^{c\sqrt{N}}$ for some constant $c$ close to unity.
To see this heuristically, note that
\cite[(3.9)--(3.11)]{Corona2018} show that the continuous single-layer operator from a source sphere of radius $R_p<1$ to a concentric sphere of radius 1 has each eigenfunction a linear combination of vector spherical harmonics within a given degree $n$ and order $m$, with eigenvalue $R_p^n$ (up to algebraic prefactors).
Furthermore, $N$ spherical design points accurately discretize all harmonics of degree $n\le p\approx c\sqrt{N}$, for some constant $c$ close to 1. Thus the minimum eigenvalue of the matrix $\vec G$ is expected to be $R_p^{c\sqrt{N}}$, and the maximum eigenvalue $\mathcal{O}(1)$, justifying the claim.
 Numerically, the decay of singular values matches this model, for $c\approx 1$, as visualized in Figure \ref{sing_without} (Appendix \ref{decay}). %The condition number growth is explained well by $(R/R_p)^p$, for $p\approx\sqrt{N}$, the amplitude ratio for the degree-$p$ multipole $Y_p^m(\theta,\phi)r^{-(1+p)}$ between the concentric spheres.
With a typical choice $R_p = 0.7$, choosing $N$ in the range $500$ to $2000$ gives $\kappa(\vec G)$ in the range $10^5$ to $10^9$.
% \begin{remark}\label{nullspace}
 %     The smallest singular value without image points corresponds to the null space of the single layer potential for a sphere, where $\int \mathbb S(\vec x-\vec y)\vec n(\vec y)\,\mathrm dS_{\vec y} = \sum_{i = 1}^N \mathbb S(\vec x-\vec y_i)\vec n(\vec y_i)=\vec 0$, with $\vec n$ the normal at the sphere surface \cite{Hsiao2008}. 
 % \end{remark}
\end{remark}

The effects of the parameters $N$ (proxy-points per sphere) and $R_p$ (proxy-radius) are investigated in Section \ref{param_sec}. %for problems where the sources at the proxy-surface alone are enough to resolve the hydrodynamic interactions.
We also give numerical examples that show that the ill-conditioning of the MFS system matrix does not pose a problem, as long as the norm of the resulting solution vector $\vec\lambda$ is bounded,
and that the latter is insured by setting $R_p$ large enough to enclose any singularities
(i.e. images) in the continuation of the solution into the interior of the sphere.
This is standard for the MFS in other settings \cite{Barnett2007,Stein2022}.

\subsection{One-body preconditioned iterative solution for many particles}\label{solving}
 Next, we describe a technique for solving the large overdetermined and ill-conditioned system for the unknown MFS coefficients $\vec \lambda$ on $P$ particles. Let $\vec x^{(k)}_i$, $i = 1,\dots,M$, denote the collocation points on particle $k$, and $\vec y^{(k)}_j$, $j = 1,\dots, N$ denote the sources in its interior. Denoting by
 $$\vec X = \begin{bmatrix} {\vec x^{(1)}}^T & {\vec x^{(2)}}^T & \dots &{\vec x^{(P)}}^T\end{bmatrix}^T$$%\rbrace_{k=1}^P$, 
 the stacked vector of collocation points for all particles,  %.In compact mathematical notation
we seek to solve the overdetermined linear system
	\begin{equation}\label{lq}
 \vec G \vec \lambda = \vec b, \qquad\text{with }\quad \vec b = \vec u_{\text{bc}}(\vec X),
		% \min\limits_{\boldsymbol\lambda\in\mathbb R^{3NP}}\|\vec G\boldsymbol\lambda-\vec u(\vec b)\|^2_2,
	\end{equation}
 in the least-squares sense. %We here  assuming  every particle has the same number of collocation points.
 % The problem can be solved backward-stably by using e.g.~an SVD:  With $\vec G = \vec W\boldsymbol \Sigma\vec V$, the pseudo-inverse of $\vec G$ is given by 
	% \begin{equation}
	% 	\vec G^{\dagger} = \vec V\vec\Sigma^{-T}\vec W^T,
	% \end{equation}
	% where the matrix $\vec\Sigma^{-T}$ is computed as a pseudo-inverse in the sense that singular values below a threshold $\epsilon_{\text{trunc}}$ are truncated, so that the elements of $\vec\Sigma^{-T}$ are limited in magnitude. 
 % A solution to the least squares problem in \eqref{lq} is constructed by 
	% \begin{equation}
	% 	\boldsymbol \lambda = \vec V\vec\Sigma^{-T}\left(\vec U^T\vec u(\vec b)\right),
	% \end{equation} 
	% where the order of multiplication is important for backward-stability \cite{Lai2015,Malhotra2015,Stein2022,Parolin2022}. 
Here, the entire target-from-source matrix $\vec G\in\mathbb R^{3MP\times 3NP}$ has blocks $\vec G^{(kk')} \in \mathbb R^{3M\times 3N}$, each representing the effect on the surface velocities $\vec u$ of body $k$ due to the Stokeslet sources within body $k'$, giving the block form
	\begin{equation}\label{Gmat}
	\vec G = \begin{bmatrix} \vec G^{(11)} & \vec G^{(12)} & \dots & \vec G^{(1P)} \\
		\vec G^{(21)} & \vec G^{(22)} & \hdots & \hdots \\
		\vdots & \vdots & \ddots & \vdots \\
		\vec G^{(P1)} & \hdots & \hdots & \vec G^{(PP)}
	\end{bmatrix}.
\end{equation}
A dense direct least squares solution would be unfeasibly slow for large $P$, with a cost $\mathcal O\left(MP\cdot(NP)^2\right)$.  % note 3's not meaningful in big-O notation...
However, the ill-conditioning inherited from that of the MFS for each particle prevents an iterative solution in the form \eqref{lq}.
To address this we use {\em one-body right preconditioning}. The diagonal block $\vec G^{(kk)}$ takes the form of the matrix in \eqref{lsq}, with sources at $\vec y^{(k)}_j$, $j = 1,\dots,N$,  in body $k$ evaluated at the collocation points $\vec x^{(k)}_i$, $i = 1,\dots, M$, on its surface. For identical particles considered here, denote this one-body block by $\vec B$. The block-diagonal preconditioning matrix contains the blocks $\vec B^{+}$,
and hence represents the pseudo-inverses of the target-from-source matrices for each body in isolation.
This turns the ill-conditioned rectangular system with $3N$ unknowns per body into a well-conditioned square system with $3M$ unknowns per body,
which is solved via GMRES.
This idea was used by Stein \& Barnett in \cite{Stein2022}, and first considered in an MFS Helmholtz problem by Liu \& Barnett in \cite{Liu2016}. It has not to our knowledge been proposed in the 3D
Stokes MFS setting before.

One-body preconditioning has subtleties, so we now present it in
detail. For each particle $k$, we work with the unknown coefficient vector $\vec \mu^{(k)} := \vec B\vec \lambda^{(k)}\in\mathbb R^{3M}$, representing surface collocation velocities, instead of its source strength vector $\vec\lambda^{(k)}\in\mathbb R^{3N}$.
Substituting $\vec\lambda^{(k)} = \vec B^+\vec\mu^{(k)}$ directly into \eqref{lq} would result in diagonal blocks $\vec B\vec B^+\in\mathbb R^{3M\times3M}$ that are rank deficient because $M>N$, rather than the identity matrix $\vec I \in \mathbb R^{3M\times3M}$ that expresses the true surface velocity self-interaction.
Thus we replace each such diagonal block by $\vec I$, and recover a well-conditioned formulation. The preconditioned linear system to solve, with square system matrix $\hat{\vec G}$ of size $\mathbb R^{3MP}\times\mathbb R^{3MP}$, is then
	\begin{equation}\label{resistance}
 \hat{\vec G}{\vec \mu} := 
		\begin{bmatrix}
			\vec I & \vec G^{(12)} \vec B^{+}& \dots & \vec G^{(1P)}\vec B^{+} \\
			\vec G^{(21)}\vec B^{+} & \vec I & \hdots & \hdots \\
			\vdots & \vdots & \ddots & \vdots \\
			\vec G^{(P1)}\vec B^{+} & \hdots & \hdots & \vec I
		\end{bmatrix} 
	\begin{bmatrix} \vec \mu^{(1)} \\ \vec \mu^{(2)} \\ \vdots \\ \vec \mu^{(P)}
		\end{bmatrix} = \revTwo{\begin{bmatrix}
		\vec b^{(1)} \\ \vec b^{(2)} \\ \vdots \\ \vec b^{(P)}
	\end{bmatrix}.}	
	\end{equation}
 Taking the SVD, $\vec B = \vec U\boldsymbol \Sigma\vec V$, and the pseudo-inverse $\vec B^{+}$
 is
	\begin{equation}
		\vec B^{+} = \vec V\vec\Sigma^{+}\vec U^T,
	\end{equation}
 where the matrix $\vec\Sigma^{+}$ is computed as a pseudo-inverse in the sense that singular values below a threshold $\epsilon_{\text{trunc}}$ are truncated, so that the elements of $\vec\Sigma^{+}$ are limited in magnitude.
 However, the matrix $\vec B^{+}$ cannot be formed and then applied in a backward stable fashion \cite{Lai2015}; two factors must be stored and applied sequentially to prevent catastrophic cancellations. The matrix in \eqref{resistance} could therefore not be built explicitly, even if speed and storage limitations were ignored.
	Instead, to perform the matrix-vector multiply with $\hat{\vec G}$ in the GMRES iteration to solve \eqref{resistance}, we propose the following three steps:
	\begin{enumerate}
		\item Apply the self-interaction pseudo-inverse $\vec B^{+}$ backward-stably (using a pair of dense matrix-vector multiplies as in \cite{Lai2015,Malhotra2015,Stein2022,Parolin2022}) to the current coefficient vector for each particle,
\begin{equation}\label{transform}
	\hat{\vec\lambda}^{(k)} = \vec V\vec\Sigma^{+}\left({\vec U}^T\vec \mu^{(k)}\right),\quad k = 1,\dots,P.
\end{equation}
		\item Compute, either using direct summation or the fast multipole method (FMM), the all-to-all interaction
		\begin{equation}
			\vec u_{\text{All}}= \vec G\hat{\vec \lambda},
   \label{uAll}
		\end{equation}
		between sources and collocation points for all particles, recalling \eqref{Gmat}.
		\item %To get the matrix-vector multiply $\hat{\vec G}{\vec\mu}$, subtract from each block of $\vec u_{\text{All}}$ the self-interaction 		\begin{equation}
 	%\hat{\vec u}^{(i)} = \vec B\hat{\vec \lambda}^{(i)}, \quad i = 1,\dots,P
 %\end{equation}and add $\tilde{\vec u}^{(i)} = \vec\mu^{(i)}$, corresponding to the identity blocks in \eqref{resistance}.
 %With more compact mathematical notation, 
 Subtract off the incorrect self-interactions and replace them by the effect of identity blocks, to give
  \begin{equation}
      \hat{\vec G}{\vec \mu} = \vec u_{\text{All}}-\vec B_{\text{blkdiag}}\hat{\vec \lambda}+\vec\mu,
  \end{equation}
  where $\vec B_{\text{blkdiag}}$ is the block-diagonal matrix with diagonal blocks $\vec B$.

		%\vspace*{-5ex}
	\end{enumerate}
    \revTwo{If direct summation is used, of course steps 2 and 3 could be combined, but in order to tackle the largest problems we will focus on the use of the FMM where the steps must be separate.}
	Upon GMRES convergence, the vector $\vec\mu$ is again via \eqref{transform} transformed back to the source strength vector $\vec\lambda$. The latter can then be used for solution velocity or pressure field evaluation throughout the fluid.
 %Note specifically that step 2 above is particularly well suited for a fast summation method; for the results presented in Section \ref{large_ex} we will use the fast multipole method (FMM). 

%%%%%%%%%%%%%%%%%%%%%%%%%%%%%%%%%%%%%
\subsection{Image enhancement for problems with strong lubrication forces}\label{image_enhance}
    % Section intro
	For challenging problems with strong lubrication, Stokeslets on interior proxy-spheres cannot provide enough resolution unless $R_p$ is made very close to $1$ and the number of proxy-sources $N$ is made impractically large. 
 It is much more efficient to keep $N$ and $R_p$ small while adding 
 other types of fundamental solutions at locations informed by
 image methods discussed in the introduction. Hence, only with the basic discretization does $N$ denote the total number of vector valued sources.
 
	Due to the linearity of the set of PDEs in \eqref{stokeseq}, not only the Stokeslet in \eqref{stokeslet}, but also differentiated Green's functions are fundamental solutions that satisfy the set of equations everywhere except at the source point. The most basic ones are the rotlet, $\mathbb R$, stresslet, $\mathbb T$ and potential dipole, $\mathbb D$, defined in the following. We start by the gradient of the Stokeslet, the force dipole tensor \cite{Pozrikidis2017,Graham}. The symmetric part of this tensor can be written in terms of the rank-3 tensor known as the stresslet, $\mathbb T$, with
    \revTwo{%r2
 \cite{Graham}%\cite{Pozrikidis2017}
	\begin{equation}
		\mathbb T_{ijk}(\vec r) = -\frac{1}{8\pi\mu}\frac{3r_ir_jr_k}{\|\vec r\|^5}
  % was just a scaling of -6
	\end{equation}
    }  %r2
	and a flow field represented by the stresslet can be expressed as
	\begin{equation}
		 u_i(\vec x) = \mathbb T_{ijk}(\vec x-\vec y)C_{jk}.	
	\end{equation}
	Here, Einstein's summation convention is used and $\vec C$ is a constant strength matrix for the stresslet located at $\vec y$  \revTwo{\cite{Graham}}. %\cite{Pozrikidis2017}.
	The anti-symmetric part of the force dipole tensor is the rotlet, $\mathbb R$, with 
 \revTwo{
	\begin{equation}\index{rotlet}
		\mathbb R(\vec r) = \frac{1}{8\pi\mu}\frac{1}{\|\vec r\|^3}\begin{bmatrix}
			0 & r_3 & -r_2 \\ -r_3 & 0 & r_1 \\ r_2 & -r_1 & 0
		\end{bmatrix}.
	\end{equation}
	The associated flow field is given by
	\begin{equation}
		\vec u(\vec x) = \mathbb R(\vec x-\vec y)\vec t, 
	\end{equation}}
	where $\vec t$ is a point torque at $\vec y$. The fourth type of singularity,
	%By instead differentiating a point source of fluid mass, 
 the potential dipole, $\mathbb D$, can be interpreted as an opposite source and sink of fluid mass separated by an infinitesimal distance along a given unit vector \cite{Graham}. 
	A potential dipole with strength $\vec d$ at $\vec y$ gives rise to a fluid velocity field of the form 
	\revTwo{\begin{equation}
		\vec u(\vec x) = \mathbb D(\vec x-\vec y)\vec d,
	\end{equation}
	where \begin{equation}\label{doublet}
		\mathbb D(\vec r) =\frac{1}{4\pi}\left( -\frac{\vec I}{\|\vec r\|^3}+3\frac{\vec r \vec r^T}{ \|\vec r\|^5}\right).
	\end{equation}}
	%Higher order tensors can also be derived, defined in terms of higher order derivatives of the Stokeslet and point source singularities, but here we will not continue beyond $\mathbb S$, $\mathbb T$, $\mathbb R$ and $\mathbb D$.
	With $\mathbb S$, $\mathbb T$, $\mathbb R$ and $\mathbb D$ in place,
 %and following linearity,
 a generalization of the representation in \eqref{sum} is to write the solution to \eqref{stokeseq} due to a single particle  as 
	\begin{equation}\label{lin_comb}		\begin{aligned}
			u_i(\vec x) = \sum_{\alpha=1}^{N_{\text{S}}}\mathbb S_{ij}(\vec x-\vec y_{\alpha}^{\text{S}}) f_{\alpha,j} + \sum_{\alpha=1}^{N_\text{T}} \mathbb T_{ijk}(\vec x-\vec y_{\alpha}^{\text{T}})C_{\alpha,{jk}}+ \\+ \sum_{\alpha=1}^{N_\text{R}}\mathbb R_{ij}(\vec x-\vec y_{\alpha}^{\text{R}}) t_{\alpha,j} + \sum_{\alpha=1}^{N_\text{D}}\mathbb D_{ij}(\vec x-\vec y_{\alpha}^{\text{D}}) d_{\alpha,j},
		\end{aligned}
	\end{equation} 
where the singularity locations, $\lbrace\vec y_{\alpha}^{\text{S}}\rbrace_{\alpha = 1}^{N_{\text{S}}}$, $\lbrace\vec y_{\alpha}^{\text{T}}\rbrace_{\alpha = 1}^{N_{\text{T}}}$, $\lbrace\vec y_{\alpha}^{\text{R}}\rbrace_{\alpha = 1}^{N_{\text{R}}}$ and $\lbrace\vec y_{\alpha}^{\text{D}}\rbrace_{\alpha = 1}^{N_{\text{D}}}$, may vary with singularity type.

\begin{remark}[Stresslet strengths]\label{stress_rem} All the terms in the sums in \eqref{lin_comb} are of the form of a matrix-vector product, except for the terms involving the stresslet, $\mathbb T_{ijk}(\vec x-\vec y_{\alpha}^{\text{T}})C_{\alpha,{jk}}$. This can however be accomplished by restricting the unknown strength matrix $\vec C_{\alpha}$ to be of the form $C_{\alpha,{lm}}=h_{\alpha,l}q_{\alpha,m}$, for some known vector $\vec q_{\alpha}\in\mathbb R^3$ and unknown $\vec h_{\alpha}\in\mathbb R^3$. %With the representation in \eqref{lin_comb}, it then remains to determine weights, in terms of the vectors $\vec f^{\alpha}$, $\vec h^{\alpha}$, $\vec t^{\alpha}$ and $\vec d^{\alpha}$, and source locations $\vec y^{\alpha}$ for each singularity type inside the particle volume to satisfy the given boundary 
 %conditions $\vec u^{\text{bc}}$. 
 \end{remark} 

The numbers $N_{\text{S}}, N_{\text{T}}, N_{\text{R}}$ and $N_{\text{D}}$ may vary per particle, since particles have different numbers of neighbors, each at different distances.
To indicate this, we use $N^{(k)}$ to denote the total number of vector source unknowns for particle $k$, with 
\begin{equation}
\label{Nk}
 N^{(k)} = N^{(k)}_{\text{S}}+N^{(k)}_{\text{T}}+N^{(k)}_{\text{R}}+N^{(k)}_{\text{D}}.
 \end{equation}
As with our basic discretization, strengths for the different source types are determined in the least-squares sense via collocation at the particle surface (we postpone a precise description of collocation point parameters to Section~\ref{lub}).
This results in each matrix block $\vec G^{(kk')}$ in \eqref{Gmat} being replaced by a horizontal stack of four matrix blocks, one for each type of unknown in \eqref{lin_comb}, to give
\[
\vec G^{(kk')} = \bigl[ \vec G_\text{S}^{(kk')}, \vec G_\text{T}^{(kk')}, \vec G_\text{R}^{(kk')}, \vec G_\text{D}^{(kk')} \bigr].
\]
Explicitly, the elements in each block are $\vec G_{\text{S},ij}^{(kk')} := \mathbb S(\vec x_i^{(k)},\vec y_j^{(k')})$,
and analogously for types $\mathbb R$ and $\mathbb D$, while stresslets instead have elements
$\vec G_{\text{T},ij}^{(kk')} := \mathbb T(\vec x_i^{(k)},\vec y_j^{(k')}) \cdot \vec h^{(k')}_j$.
%\begin{equation} 
%\begin{aligned}
%   \vec G^{(kk')} = \left[
  % } 
%  \mathbb S(\vec x^{(k)},\vec y_1^{(j)}),\,\dots,\,\mathbb S(\vec x^{(i)},\vec y_{N_S}^{(j)}),\,  \mathbb T(\vec x^{(i)},\vec y_1^{(j)})\cdot \vec h^{(j)}_1,\, \dots, \,\mathbb T(\vec x^{(i)},\vec y_{N_T}^{(j)})\cdot\vec h^{(j)}_{N_T},\,\dots, \right. \\ \left.  \mathbb R(\vec x^{(i)},\vec y_1^{(j)}),\, \dots, \,\mathbb R(\vec x^{(i)},\vec y_{N_R}^{(j)}),\,\mathbb D(\vec x^{(i)},\vec y_1^{(j)}), \, \dots, \,\mathbb D(\vec x^{(i)},\vec y_{N_D}^{(j)})\right].
%      \end{aligned}
   %\end{bmatrix}
%\end{equation}
The full bare system matrix $\vec G$ is now of size 
$3\sum_{k=1}^P M^{(k)}$ by $3\sum_{k=1}^P N^{(k)}$, and the diagonal self-interaction blocks $\vec G^{(kk)}$ in the preconditioning strategy of Section \ref{solving} (previously all denoted by $\vec B$) are now in general different.

%Note particularly that the accuracy and stability of the representation in \eqref{lin_comb} rely on the source points' locations. The proper locations, in turn, depend on the singularities of the analytic continuation of the solution field into the particle interior \cite{Doicu2000,Barnett2007,Stein2022}.

 We now turn to the placement of these different types of fundamental solutions along the center-center line for each pair of close-to-touching particles. Our strategy is motivated by the reflection of points in spheres in the Laplace case, which we now review.
 %We here describe what is done in practice and introduce some necessary terminology. 
 The reflection of a point $\vec x$ exterior to a unit sphere with center $\vec c^{(k)}$ gives the \emph{image} point \cite{kelvin2spheres,Cheng2000}
	\begin{equation}\label{reflect}
		\vec x' = \vec c^{(k)}+ \frac{ \left(\vec x-\vec c^{(k)}\right)}{\|\vec x-\vec c^{(k)}\|^2}.
	\end{equation}
 The (unique) decaying exterior Laplace BVP solution for boundary data due to a charge at
 $\vec x$ is given by the potential from a (modified) charge at $\vec x'$.
 In other words, the continuation of this BVP solution into the sphere, as a Laplace solution,
 has a singularity only at $\vec x'$.
 This has two consequences for the MFS for exterior BVPs (e.g. \cite{Ka89}): i) such a singularity
 (if close to the surface) controls the convergence rate, and
ii) if all such singularities are not enclosed (``shielded'') by the proxy surface, the MFS coefficient
 norm will grow exponentially, causing a numerical round-off instability \cite{Doicu2000,Barnett2007}.

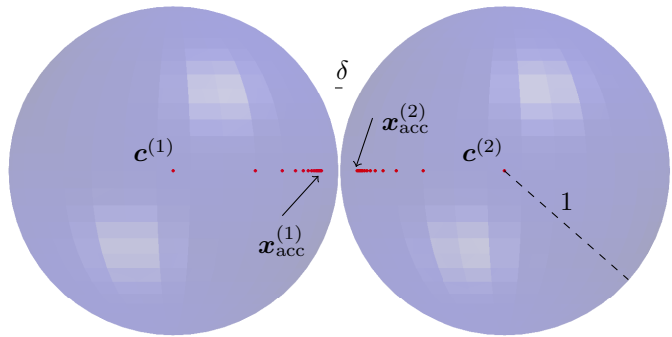
\begin{figure}[h!]
		\centering\begin{tikzpicture}
			
			% Draw a circle using TikZ
			%\draw (0,0) circle (1);
			
			%
			%% Label the start and target points
			%\node[below] at (start) {Start};
			%\node[above] at (target) {Target}
			
			% Insert an external graphic (e.g., a PNG file)
			\node at (3,0) {							\includegraphics[trim={1.9cm 2.4cm 1cm 2cm},clip,width=0.55\textwidth]{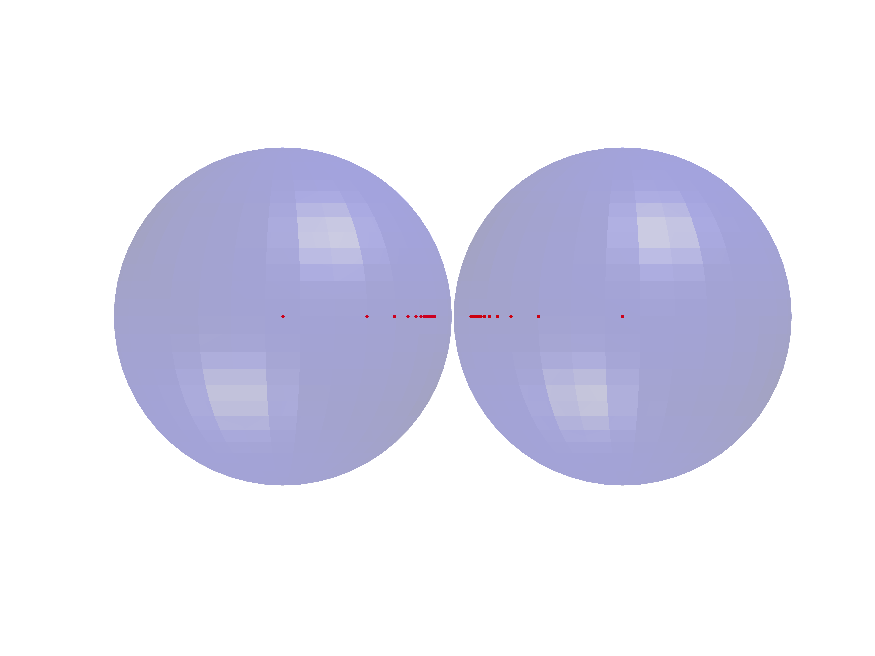}};
			
			\coordinate (start) at (3.3,0.7);
			\coordinate (target) at (3.08,0.05);
			%
			%% Draw the arrow
			\draw[->] (start) to (target);
			\node[right] at (start) {$\vec x^{(2)}_{\text{acc}}$};
			\coordinate (start) at (2.1,-0.6);
			\coordinate (target) at (2.58,-0.07);
			%
			%% Draw the arrow
			\draw[->] (start) to (target);
			\node[below] at (start) {$\vec x^{(1)}_{\text{acc}}$};
			\coordinate (start1) at (0.4,0);
			\coordinate (start2) at (4.77,0);
			\node[above] at (start1) {$\vec c^{(1)}$};
			\node[above] at (start2) {$\vec c^{(2)}$};

            \coordinate (start) at (2.9,1.1);
			\coordinate (target) at (2.8,1.1);
            \draw[-] (start) to (target);
            \node[above] at (start) {$\delta$};

            \coordinate (start3) at (5.05,0);
			\coordinate (target3) at (6.7,-1.45);
            \draw[-] (start3) to (target3) [dashed];
            \coordinate (mid) at (5.88,-0.65);
            \node[above] at (mid) {$1$};
   
		\end{tikzpicture}
		
		%	\caption{Source and target points for two spheres in 3D.}
		%	\label{imagespheres}
		%	\end{subfigure}
		\caption{Lines of image points obtained by repeated reflections of each sphere center in the other sphere via \eqref{reflect}. Accumulation points for the lines of images depend on the separation distance $\delta$ between the two spheres.}
		\label{imagespheres}			
	\end{figure}

 Consider now two disjoint unit spheres centered at $\vec c^{(1)},\vec c^{(2)}$ and separated by a small distance $\delta$, as illustrated in Figure \ref{imagespheres}. Starting with the center coordinates of each of the spheres reflected in the other sphere, and then computing reflections of reflections by applying \eqref{reflect} iteratively, an infinite sequence of image points is formed in each sphere, lying on a single line as indicated in the figure. The image points accumulate at
$\vec x_\text{acc}^{(1)}$ in sphere 1 and at 
$\vec x_\text{acc}^{(2)}$ in sphere 2.
For a given $\delta$, these points may be solved for as fixed points of
the iteration, giving
\begin{equation}\label{acc_points}
    \begin{aligned}
        \vec x_{\text{acc}}^{(1)} &=\frac{\vec c^{(1)}+\vec c^{(2)}}{2}-\sqrt{\delta+\delta^2/4}\cdot \frac{\vec c^{(2)}-\vec c^{(1)}}{\|\vec c^{(2)}-\vec c^{(1)}\|},\\
        \vec x_{\text{acc}}^{(2)} &=\frac{\vec c^{(1)}+\vec c^{(2)}}{2}+\sqrt{\delta+\delta^2/4}\cdot \frac{\vec c^{(2)}-\vec c^{(1)}}{\|\vec c^{(2)}-\vec c^{(1)}\|}.\\
    \end{aligned}
\end{equation}
As $\delta\to 0$, these accumulation points lie a distance ${\mathcal O}(\delta^{1/2})$ inside the surface,
which would ruin the convergence and stability of the MFS when using concentric proxy spheres.
Attempting to augment the MFS with the true multiple-reflection series (as in Cheng--Greengard \cite{Cheng1998,Cheng2000})
would result in too many images at small $\delta$,
and---as described in the introduction---the continuous line-densities generated in the Stokes
case render such analytic multiple reflections intractable.
Noticing that such a 2-sphere Stokes solution would involve
sources only on the lines from centers to the accumulation points motivates our proposal
to discretize these lines with a small number of quadrature nodes of unknown strengths.

 For unit spheres  discretized with concentric interior proxy-spheres of radius $R_p$, we first need a criterion for when additional sources along such an image line out to the accumulation point are needed at all.
A simple proxy-sphere of radius $R_p$ is stable for MFS (in the sense of having small coefficient norm) only when it contains the accumulation point (and thus all other images).
Let
\begin{equation}\label{Racc}
R_\text{acc}(\delta) = 1+ \delta/2 - \sqrt{\delta + \delta^2/4}
\end{equation}
denote the distance of the accumulation point to the center of the sphere.
Then, given $R_p$,
let $\delta^*$ denote the critical particle-particle separation where the image accumulation point lies on the sphere with radius $R_p$. 
Precisely, $\delta^*$ solves
\[
R_\text{acc}(\delta^*) = R_p, 
\]
and it is easy to verify that the unique solution to the above is
\begin{equation}
\delta^* = \frac{(1-R_p)^2}{R_p}.
\label{deltastar}
\end{equation}
If a pairwise distance $\delta$ is less than $\delta^*$,
then we will need to 
place additional image sources which discretize the part of the line from the center to $\vec x_\text{acc}$ lying outside of the proxy-sphere. %The interval of interest, $[R_p,R_{\text{acc}}]$ 
The line segment of interest in the interior of the first particle is hence parameterized by
%More specifically,  we approximate the distribution of images with half a Chebyshev grid in the interior of each particle, clustered towards the accumulation points. 
% In the interior of the first sphere, the additional sources are then placed on the line $\vec l^{(1)}(t)$ defined by 
	\begin{equation}\label{line_points}
		% \vec l^{1}(t) = \vec c_{1} + \left(R_p+t(\|\vec x^{1}_{\text{acc}}(\delta)-\vec c^{1}\|-R_p)\right)\frac{\vec c^{2}-\vec c^{1}}{\|\vec c^{2}-\vec c^{1}\|}, \quad t\in(0,1],
  		\vec l^{(1)}(t) = \vec c^{(1)} + t\frac{\left(\vec c^{(2)}-\vec c^{(1)}\right)}{\|\vec c^{(2)}-\vec c^{(1)}\|}, \quad t\in(R_p,R_{\text{acc}}(\delta)],
	\end{equation}
 and can be discretized with $n_{\text{im}}$ nodes $\vec l_j^{(1)} = \vec l^{(1)}(t_j)$, $j = 1,\dots,n_{\text{im}}$, with $t_j$ given by
 \begin{equation}\label{tj}
     t_j = R_p(1+\Delta)+\left(R_{\text{acc}}(\delta)-R_p(1+\Delta)\right)\cos\left(\frac{(j-1)}{2n_{\text{im}}}\pi\right).
 \end{equation}
The nodes are constructed from half
of a set of type-2 Chebyshev nodes,
%Chebyshev interval,
clustered towards (and including) the accumulation point; see Figure~\ref{colloc_images}.
This mimics a quadrature scheme for the interval.
% s well as the clustering of the true image series ...?
For conditioning reasons, the $n_{\text{im}}$th node is here given an offset from the proxy-surface via a small distance $\Delta$, which we fixed at $\Delta =0.05$. 

As an ansatz, a linear combination of $\mathbb S$, $\mathbb T$, $\mathbb R$ and $\mathbb D$ will be centered at each such image node, with an efficient combination chosen numerically in Section \ref{lub}.
 For all particle-particle distances larger than $\delta^*$, the Stokeslet sources on the proxy-surface alone will allow an MFS representation with small coefficient norm.
 Considering the target accuracy of $10^{-3}$, we will in Section \ref{param_sec} choose an $R_p$ such that $\delta^*<0.3$, meaning that image points are only required for separations less than this.

% The locations of source points in \eqref{line_points}-\eqref{tj} depend on the particle-particle distance for the pair and hence change depending on the particle configurations.
 
 % One-body preconditioning as introduced in Section \ref{solving} can still be applied, but with the difference that the total number of vector-valued sources per particle now can vary for each particle depending on resolution requirements set by the proximity to its neighbors. The self-interaction blocks of the target-from-source matrix will therefore be different for each particle and need to be recomputed. 
 In addition, the fundamental solutions $\mathbb T$, $\mathbb R$ and $\mathbb D$ are all singular at $\vec x=\vec y_{\alpha}$, similarly as the Stokeslet $\mathbb S$.	In a neighborhood of the additional clustered source points for every near contact of a sphere, collocation points  are therefore refined on the sphere surface, as in Figure \ref{colloc_points}. This has to be done locally, as any additional collocation points  will increase the number of unknowns to solve for in the preconditioned system. In our adaptive scheme presented in Section \ref{lub}, both source and collocation points will be set adaptively based on the particle-particle separation for every neighboring pair.
 
 \begin{remark}[Extracting net forces and torques]
 Having solved for the source strengths $\vec\lambda$ using preconditioned GMRES, the net force $\vec F^{(k)}$ and torque $\vec T^{(k)}$ on particle $k$ are %in a resistance problem determined in a post-processing step.
 extracted as follows.
Given point forces $\vec f_{\alpha}$ and torques $\vec t_{\alpha}$ corresponding to Stokeslets and rotlets in the particle interior,
	\begin{equation}
		\vec F^{(k)} = \sum_{\alpha=1}^{N_{\text{S}}}\vec f_{\alpha},\quad \vec T^{(k)}= \sum_{\alpha=1}^{N_{\text{R}}}\vec t_{\alpha}+\sum_{\alpha=1}^{N_{\text{S}}}\left(\vec y_{\alpha}^{\text{S}}-\vec c^{(k)}\right)\times \vec f_{\alpha}.
	\end{equation}
Conveniently, this avoids any need to integrate stresses over the sphere surfaces.
\end{remark}

	% \begin{remark}\label{rem1}
	% 	In any dynamic simulation, there is a risk for diminishing particle separations as an artifact of time-stepping. In this work, any \emph{physically relevant} particle separation is addressed, and by that we mean any separation larger than $10^{-3}R$. Practically, smaller separations would in a dynamic simulation be adjusted by a contact avoiding strategy. \question{Keep or move up?}
	% \end{remark}

	%The problem where net forces and torques are to be determined, given translational and angular velocities of the particles, is often referred to as the Stokes resistance problem. %Moved to intro 

	%in contrast to point-wise direct collocation \cite{Stein2022, Antunes2022}. 
	%\cite{Stein2022} and Antunes \cite{Antunes2022}. This follows a similar path as in the works by Stein \& Barnett \cite{Stein2022} and Antunes \cite{Antunes2022} and is different from a classic Nystr\"{o}m distcretisation of a boundary integral method. 

	%%%%%%%%%%%%%%%%%%%%%%%%%%%%%%%%%%%%%%%%%%%%%%%%%%%%%%%%%%%
	%	\clearpage
	\section{Parameters for the basic concentric proxy sphere}\label{param_sec}
	
	%%%%%%%%%%%%%%%%%%%%%%%%%%%%%%%%%%
%%	\subsection{Choice of proxy and collocation points}\label{proxy}
 This section serves to optimize the MFS proxy sphere parameters (the number of proxy-sources $N$ and the radius of the proxy-surface $R_p$) used in later experiments.  No image points are used, and the number of collocation points is fixed at $M \approx 1.2N$, as discussed in Section~\ref{MFS}.
 
As a quantification of the error, we will here and later use the pointwise relative residual, as defined by
 \begin{equation}\label{residual_x}
 \epsilon_{\text{res}}(\vec x)\coloneqq \|\vec u(\vec x)-\vec u_{\text{bc}}(\vec x)\|_{2}/\|\vec u_{\text{bc}}(\vec x)\|_{2}, 
 \end{equation}
 and 
  \begin{equation}\label{residual}
 \epsilon_{\text{res}}^{\max}\coloneqq \max\limits_{\vec x\in \partial\Omega} \left(\epsilon_{\text{res}}(\vec x)\right), 
 \end{equation} 
its maximum approximated using a dense set of points sampled on the particle surfaces. We also sometimes use the maximum force/torque error relative to a reference, as defined by 
\begin{equation}\label{force_err}
\epsilon_{\text{FT}}\coloneqq \|\vecmc F-\vecmc F_{\text{ref}}\|_{\infty}/\|\vecmc F_{\text{ref}}\|_{\infty},
\end{equation}
with $\vecmc F$ a stacked vector of forces and torques for all particles in the system. 
  %a Fibonacci grid \cite{Marques2021}, with spherical coordinates $\lbrace\theta,\phi\rbrace_{i = 1}^N$ of an $N$-point grid defined by 
	%		\begin{equation}
	%			\begin{aligned}
	%				\theta_j &= \arccos(1-2j/N),\\
	%				\phi_j &  = \frac{2\pi j}{\Phi}\mod 2\pi, \quad 0\leq j<N,
	%			\end{aligned}
	%		\end{equation}
	%		with $\Phi = (1+\sqrt{5})/2$ the golden ratio.  
%  \begin{figure}[h!]
%  \begin{minipage}[t]{0.5\textwidth}
% 		\centering
%   \hspace*{-3ex}
% 		\includegraphics[trim = {0cm 1cm 0cm 2cm},clip,width=1.2\textwidth]{figures/sphere.eps}
% 		\caption{Original source points (red) and collocation points (blue) for each sphere obtained from spherical design nodes.}
% 		\label{basic_disc}
% \end{minipage}~~
%  \begin{minipage}[t]{0.5\textwidth}
% 	   \centering
% 		\includegraphics[trim = {5cm 16.2cm 7cm 4.4cm},clip,width=0.95\textwidth]{figures/convergence_trans2.pdf}
% 		\caption{ $\delta$ apart, affected by two types of boundary conditions, as described in main text, and discretized with $N$ source points at inner proxy-surfaces of radius $R_p = 0.6$. Dashed curves show convergence  as predicted by $\mathcal O(R_{\text{acc}}^{\sqrt{N}})$.}
%   \label{convergence_test}
% \end{minipage}
% 	\end{figure}

\begin{figure}[h!]

%  \begin{minipage}[t]{0.2\textwidth}
% % \vspace*{-20.5ex}
 
% 		\centering
%   %\vspace*{-5.5ex}
%   \begin{tikzpicture}
%     \hspace*{-22ex}  
%  %
%   \includegraphics[trim = {4cm 0cm 3.5cm 2cm},clip,width=1.1\textwidth]{figures/sphere.eps}

%   \coordinate (start) at (-1.8,2.82);
% 			\coordinate (target) at (-0.95,3.34);
% 			%
% 			%% Draw the arrow
% 			\draw[|-|] (start) to (target);
% 			\node[right] at (start) {\small{$\medspace\medspace R_p$}};
% \tikz \fill [white] (-1.5,0.5) rectangle (-1.6,0.6);
%     \end{tikzpicture}

% 		% \includegraphics[trim = {4cm 1cm 3cm 2cm},clip,width=1.15\textwidth]{figures/sphere.eps}
%   \vspace*{5.9ex}
% 		\caption{Original source points (red) and collocation points (blue) for each sphere obtained from spherical design nodes.}
% 		\label{basic_disc}
% \end{minipage}~~~~
% \begin{minipage}[t]{0.8\textwidth}
 \begin{subfigure}[t]{0.49\textwidth}
	   \centering
		\includegraphics[trim = {5.3cm 16.2cm 8cm 4.4cm},clip,width=0.8\textwidth]{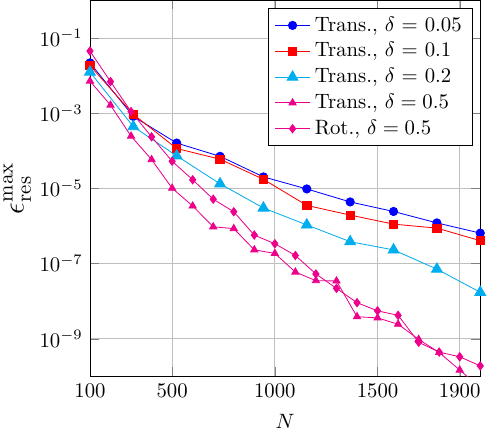}
		\caption{Translation or rotation BC:  $\epsilon_{\text{res}}^{\max}$.}
  \label{convergence_test}
  \end{subfigure}~~
  \begin{subfigure}[t]{0.49\textwidth}
	   \centering
		\includegraphics[trim = {5.3cm 16.2cm 8cm 4.4cm},clip,width=0.8\textwidth]{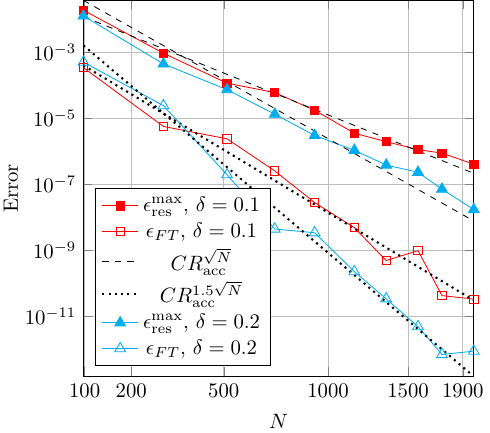}
		\caption{Translation BC: $\epsilon_{\text{res}}^{\max}$ and  $\epsilon_{FT}$, with axis linear in $\sqrt{N}$.}
  \label{convergence_testb}
  \end{subfigure}
  \caption{Convergence test for two unit spheres $\delta$ apart, discretized with $N$ source points at inner proxy-surfaces of radius $R_p = 0.6$, and affected by two types of boundary conditions, translation or rotation, as described in main text. Panel (a) displays rate dependence on $\delta$. In panel (b), dashed curves show convergence  as predicted by $\mathcal O(R_{\text{acc}}^{\sqrt{N}})$ for the max relative residual $\epsilon_{\text{res}}^{\max}$, while the max force/torque error $\epsilon_{FT}$, relative to a reference solution determined with $N=3529$, converges at a rate of  $\mathcal O(R_{\text{acc}}^{1.5\sqrt{N}})$.}
%\end{minipage}
	\end{figure}
 
Throughout this section, the geometrical setup is two unit spheres separated by a distance $\delta$. We start with a convergence test in $N$ with the two particles either rotating relative to each other, with velocities $\vec v^{(1)} = \vec 0$, $\vec v^{(2)} = \vec 0$, $\vec\omega^{(1)} =\vec e_{z}$ and $\vec\omega^{(2)} = -\vec e_{z}$, or translating together, with $\vec v^{(1)} = \vec v^{(2)} = \vec e_x$ and $\vec\omega^{(1)} = \vec\omega^{(2)} = \vec 0$. The resulting $\epsilon_{\text{res}}^{\max}$ as function of $N$ is visualized for a few different choices of $\delta$ in Figure \ref{convergence_test}.  The convergence appears to be root exponential, i.e.~spectral in $p\sim\sqrt{N}$. More specifically, we notice the rate is consistent with $\log(1/R_{\text{acc}})$, where $R_{\text{acc}}$ is given by \eqref{Racc}.
By analogy with rigorous 2D Laplace MFS theory \cite{Ka89},
the rate can be well predicted by the decay rate of spherical harmonic coefficients induced by a Stokes singularity at the radius $R_{\text{acc}}$ from the previous section.
One also must assume that the maximum harmonic degree captured by the MFS is
$p\approx\sqrt{N}$, corresponding to one proxy point per spherical harmonic term (as in Remark \ref{r:cond}).
However, the convergence of $\epsilon_{FT}$ is a bit faster:\footnote{We leave for future investigations to understand the exact convergence rates.}
Figure \ref{convergence_testb} fits a root-exponential rate of $ 1.5\log(1/R_{\text{acc}})$. In this test, $R_p = 0.6$. For larger $R_p$ relative to $R_{\text{acc}}$, $R_p$ instead starts to control the convergence rate (not shown), due to ``aliasing'' error in the MFS dominating \cite[Thm.~3 and Rmk.~4]{Barnett2007}.
%This is a consequence of aliasing, as smooth boundary data is approximated by high-frequency ``basis-functions'' located close to the surface (in terms of Stokeslets centered at the proxy-sources) \cite{Barnett2007}.
	
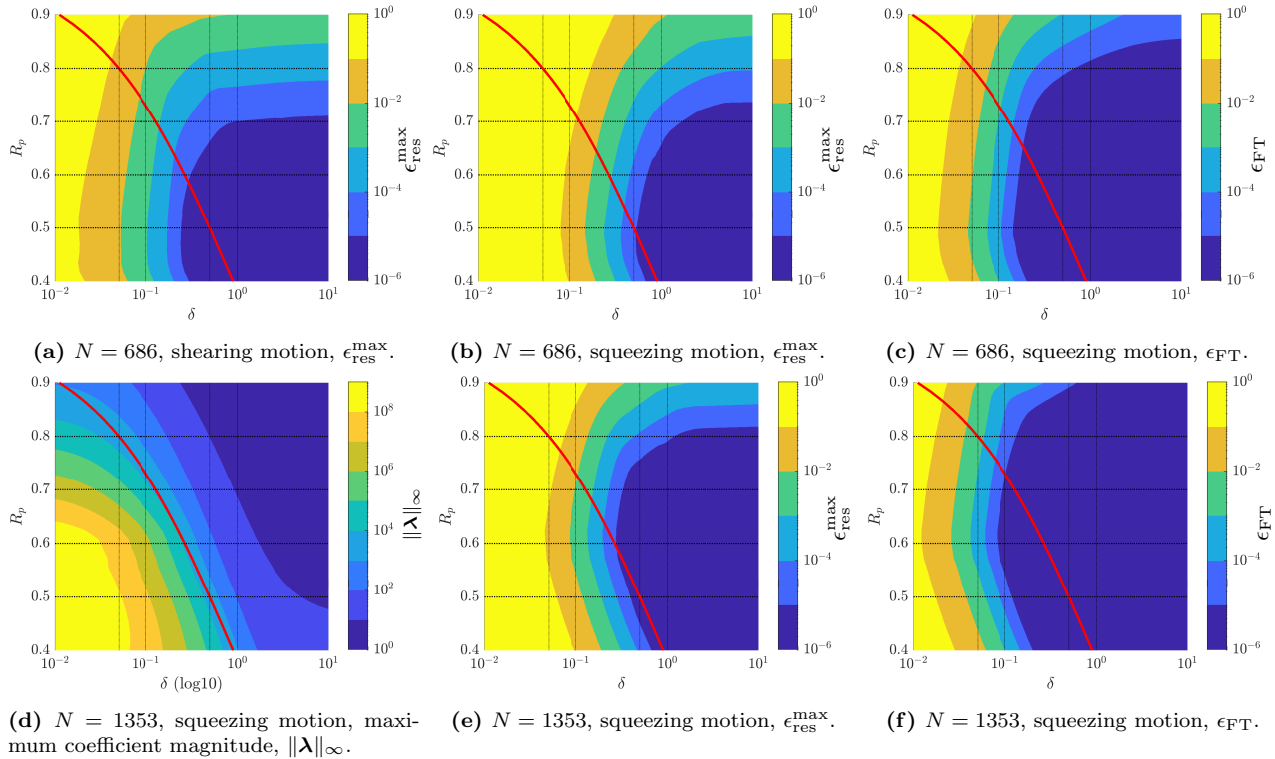
\begin{figure}[h!] %fffffffffffffffffffffffffffffffffffffffffffffffffffffffffffff
		\centering
				\begin{subfigure}[t]{0.33\textwidth}
    %\begin{minipage}{0.8\textwidth}
     \begin{overpic}[trim={0cm 0.2cm 0.2cm 0.2cm},clip,width=\textwidth]{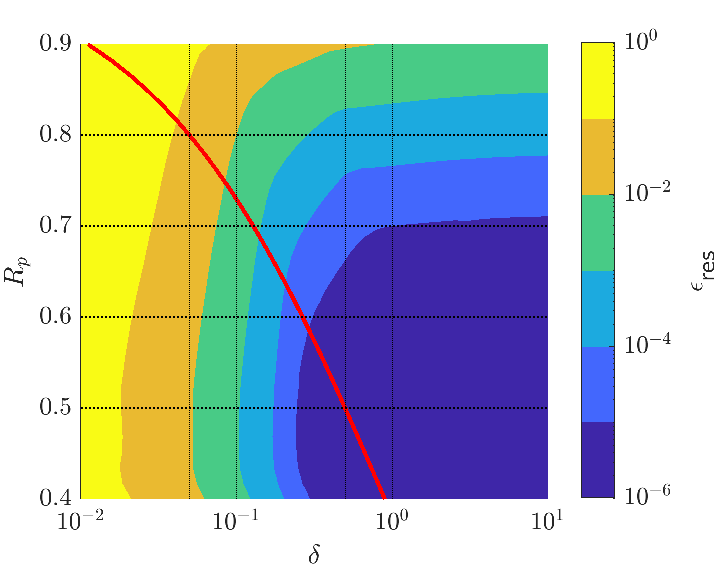}
\put(93,35){\rotatebox{90}{\colorbox{white}{\parbox{0.2\linewidth}{%
     \small$\epsilon_{\text{res}}^{\max}$}}}}
     \put(43,0){\colorbox{white}{\parbox{0.05\linewidth}{%
     \scriptsize$\delta$}}}
     \end{overpic}

			\caption{$N = 686$, shearing motion, $\epsilon_{\text{res}}^{\max}$.}
			\label{Rg_sweep_a}
  % \end{minipage}
   %\begin{minipage}{0.09\textwidth}
%        \begin{tikzpicture}
   
%     % Draw a rectangle with white fill as the background
%     \fill[white] (0,0) rectangle (2,1);
    
%     % Draw the rotated label
%     \node[rotate=90] at (1,0.5) {Rotated Label};
% \end{tikzpicture}
%    \end{minipage}
		\end{subfigure}~
			\begin{subfigure}[t]{0.33\textwidth}
   \begin{overpic}[trim={0cm 0.2cm 0.2cm 0.2cm},clip,width=\textwidth]{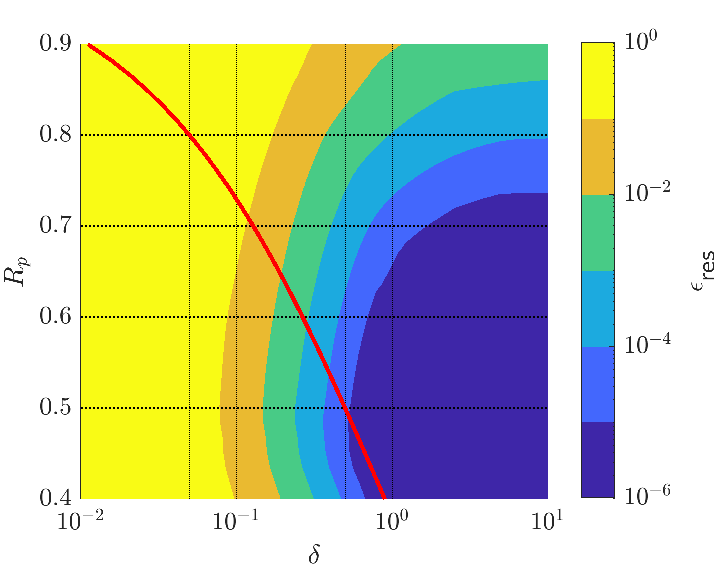}
\put(93,35){\rotatebox{90}{\colorbox{white}{\parbox{0.1\linewidth}{%
     \small$\epsilon_{\text{res}}^{\max}$}}}}
     \put(43,0){\colorbox{white}{\parbox{0.05\linewidth}{%
     \scriptsize$\delta$}}}
     \end{overpic}
			\caption{$N = 686$, squeezing motion, $\epsilon_{\text{res}}^{\max}$.}
			\label{Rg_sweep_b}
			\end{subfigure}~~
		\begin{subfigure}[t]{0.33\textwidth}
  \begin{overpic}[trim={0cm 0.2cm 0.2cm 0.2cm},clip,width=\textwidth]{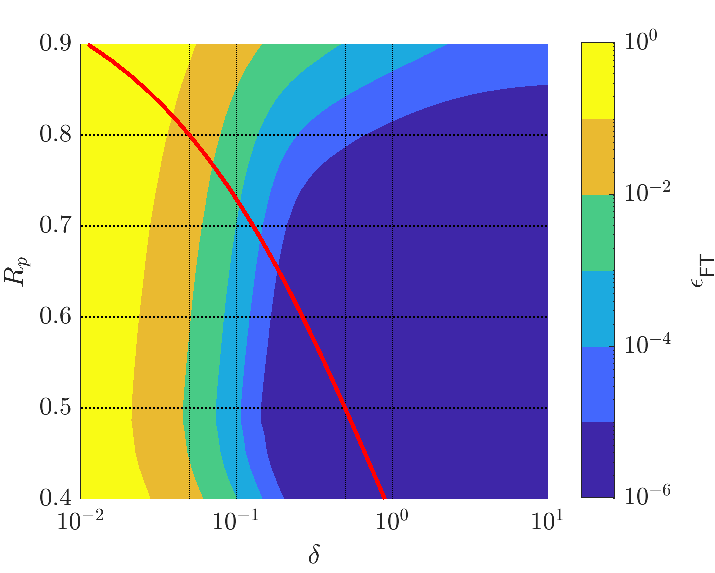}
\put(93,35){\rotatebox{90}{\colorbox{white}{\parbox{0.1\linewidth}{%
     \small$\epsilon_{\text{FT}}$}}}}
     \put(43,0){\colorbox{white}{\parbox{0.05\linewidth}{%
     \scriptsize$\delta$}}}
     \end{overpic}
			\caption{$N = 686$, squeezing motion, $\epsilon_{\text{FT}}$.}
			\label{Rg_sweep_c}
		\end{subfigure}
			\begin{subfigure}[t]{0.33\textwidth}
   \begin{overpic}[trim={0cm 0.2cm 0.2cm 0.2cm},clip,width=\textwidth]{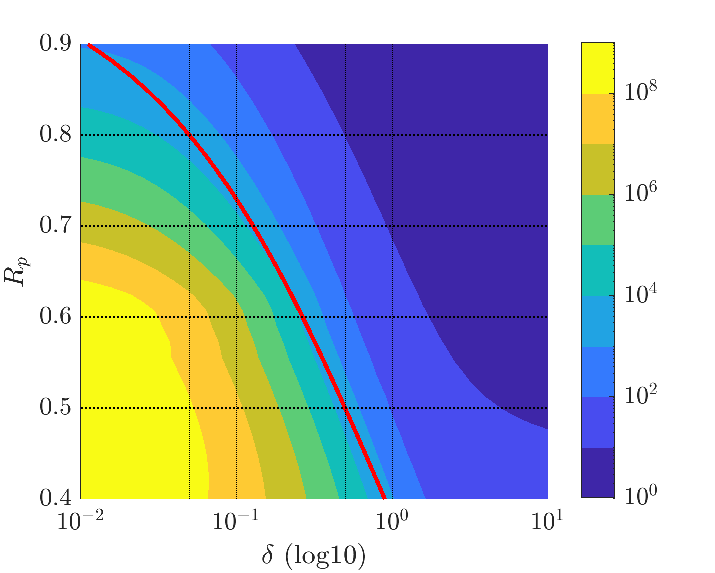}
\put(94,35){\rotatebox{90}{\colorbox{white}{\parbox{0.1\linewidth}{%
     \scriptsize{$\|\vec\lambda\|_{\infty}$}}}}}
     \put(37,0){\colorbox{white}{\parbox{0.12\linewidth}{%
     \scriptsize$\quad\delta$}}}
     \end{overpic}
		\caption{$N = 1353$, squeezing motion, maximum coefficient magnitude, $\|\vec\lambda\|_{\infty}$.}
		\label{Rg_sweep_d}
	\end{subfigure}~~
\begin{subfigure}[t]{0.33\textwidth}
\begin{overpic}[trim={0cm 0.2cm 0.2cm 0.2cm},clip,width=\textwidth]{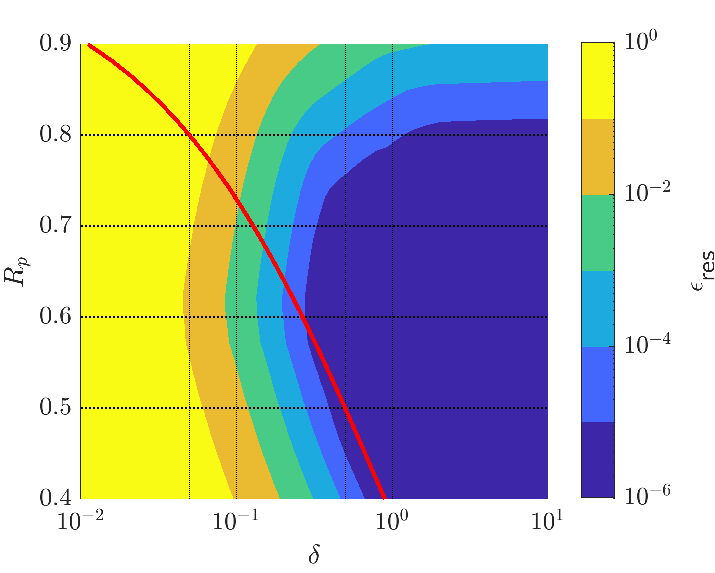}
\put(93,35){\rotatebox{90}{\colorbox{white}{\parbox{0.1\linewidth}{%
     \small$\epsilon_{\text{res}}^{\max}$}}}}
     \put(43,0){\colorbox{white}{\parbox{0.05\linewidth}{%
     \scriptsize$\delta$}}}
     \end{overpic}
	\caption{$N = 1353$, squeezing motion, $\epsilon_{\text{res}}^{\max}$.}
	\label{Rg_sweep_e}
\end{subfigure}~~
		\begin{subfigure}[t]{0.33\textwidth}
  \begin{overpic}[trim={0cm 0.2cm 0.2cm 0.2cm},clip,width=\textwidth]{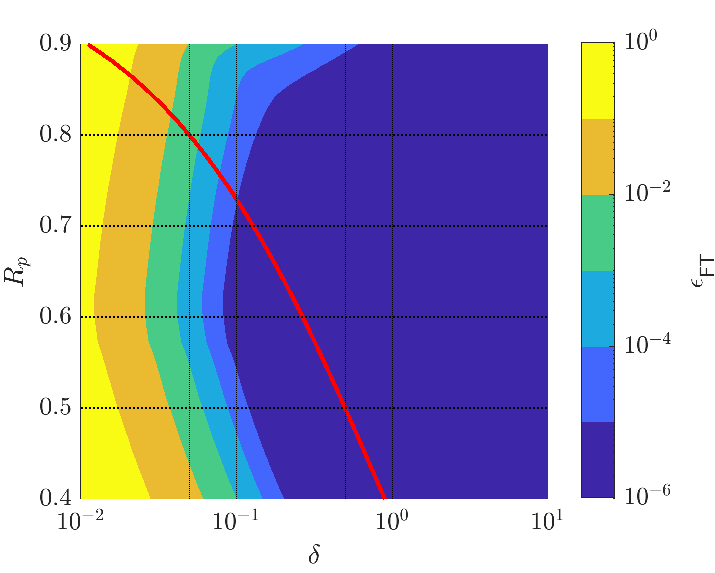}
\put(93,35){\rotatebox{90}{\colorbox{white}{\parbox{0.1\linewidth}{%
     \small$\epsilon_{\text{FT}}$}}}}
     \put(43,0){\colorbox{white}{\parbox{0.05\linewidth}{%
     \scriptsize$\delta$}}}
     \end{overpic}
	\caption{$N = 1353$, squeezing motion, $\epsilon_{\text{FT}}$.}
	\label{Rg_sweep_f}
\end{subfigure}
	%	forceerr_700parallel.eps
		%%OLD FIGS ARE HERE
%		\begin{subfigure}[t]{0.49\textwidth}
%			\includegraphics[trim={0cm 0.2cm 1cm 0.5cm},clip,width=\textwidth]{../figures/images/RpSweep700_parallel.eps}
%			\caption{$N = 700$, motion parallel to the line of centers}
%		\end{subfigure}
%		\begin{subfigure}[t]{0.49\textwidth}
%			\includegraphics[trim={0cm 0.2cm 1cm 0.5cm},clip,width=\textwidth]{../figures/images/RpSweep700_perp.eps}
%			\caption{$N = 700$, motion perpendicular to the line of centers}
%		\end{subfigure}
%		\begin{subfigure}[t]{0.49\textwidth}
%			\includegraphics[trim={0cm 0.2cm 1cm 0.5cm},clip,width=\textwidth]{../figures/images/RpSweep1400_parallel.eps}
%			\caption{$N = 1400$, motion parallel to the line of centers}
%		\end{subfigure}
%		\begin{subfigure}[t]{0.49\textwidth}
%			\includegraphics[trim={0cm 0.2cm 1cm 0.5cm},clip,width=\textwidth]{../figures/images/RpSweep1400_perp.eps}
%			\caption{$N = 1400$, motion perpendicular to the line of centers}
%		\end{subfigure}
		\caption{ A pair of unit spheres discretized with $N$ sources on proxy-surfaces of radius $R_p$ are separated by a small distance $\delta$. A sweep over $R_p$ and $\delta$ is made to investigate accuracy and MFS coefficient magnitudes. 
  The test is done for two settings: 
  %when the spheres travel with unit velocities towards each other 
  case 1 (squeezing motion), 
  %and in opposite directions perpendicular to the line of centers 
  and case 2 (shearing motion). Squeezing motion is harder to resolve, confirmed by comparing panels (a) and (b). The error in the computed force relative to the results of Brenner in \eqref{Brenner}, denoted by $\epsilon_{\text{FT}}$, is a more forgiving measure of the error than the maximum relative residual $\epsilon_{\text{res}}^{\max}$; compare panels (c) and (b) for $N=686$ and (f) and (e) for $N=1353$. For a given $R_p$, additional image points can only be added for $\delta<\delta^*$ (to the left of the red curve representing $\delta^*(R_p)$ in each panel). 
%			To pick a suitable choice of $R_p$, we consider a target tolerance of $10^{-3}$ for the maximum relative residual and choose $R_p$ as large as possible. The choices are marked with red horizontal lines in the four panels and the corresponding $\delta^*$ are marked with dashed red lines. The distance $\delta = 0.5$ is only plotted for reference.
			 For large $\delta$, $\epsilon_{\text{FT}}$ and $\epsilon_{\text{res}}^{\max}$ are capped at $10^{-6}$ which is the tolerance chosen for GMRES. Vertical black lines are drawn for reference at $\delta = 5\cdot 10^{-2}$, $\delta = 1\cdot 10^{-1}$, $\delta = 5\cdot 10^{-1}$ and $\delta = 1$.}
    %\in\lbrace5\cdot 10^{-2},1\cdot 10^{-1},5\cdot 10^{-1\rbrace}$.}		
		\label{Rg_sweep}
	\end{figure}

	Next, we describe how to optimize the parameter $R_p$ in a test where both $R_p$ and $\delta$ are varied.
 In the panels of Figure \ref{Rg_sweep}, the critical separation $\delta^*$, recalling \eqref{deltastar}, is drawn in red, as function of $R_p$. For a fixed $R_p$, $n_{\text{im}}$ image points will later be added when needed to resolve lubrication forces. %These will be placed at nodes $l_j$, recalling \eqref{tj}, along a line in the interior of each particle for any $\delta$ to the left of the red curve. 
 Looking for a small $\delta^*$ goes hand-in-hand with the wish for $R_p$ to be as large as possible, for the proxy-surface to better resolve challenging boundary data, and for the MFS coefficients $\vec\lambda$ to be reasonably small in magnitude. Two test cases are considered in this test to motivate that the discretization we choose will make a good choice for $\delta>\delta^*$ and general no-slip boundary data: 
	\begin{enumerate}
		\item The particles travel with equal and opposite velocities perpendicular to their line of centers (shearing motion). Lubrication forces are known to scale as $\mathcal O(\log(1/\delta))$ \cite{Kim1991}.
		\item The particles travel towards each other with equal and opposite velocities, $\vec v^{(1)} = -\vec v^{(2)}$, along their line of centers (squeezing motion). Lubrication forces are known to scale as $\mathcal O(1/\delta)$ \cite{Kim1991}, and the problem is hence more challenging to resolve than case 1. The resulting net force on each of the particles can be compared to an analytical expansion derived by Brenner, where \cite{Brenner1961}
		\begin{equation}\label{Brenner}
			\|\vec F^{(k)}\| = 8\pi \|\vec v^{(k)}\| \sinh \alpha \sum_{n=1}^\infty \frac{n(n+1)}{(2n-1)(2n+3)} \times
			\left[
			\frac{4 \cosh^2 (n+1/2)\alpha + (2n+1)^2\sinh^2\alpha}{2\sinh (2n+1)\alpha - (2n+1)\sinh 2\alpha} - 1
			\right],		
		\end{equation}
	and $\alpha$ is defined by $\cosh \alpha = 1 + \delta/2$ with $\alpha \sim \sqrt{\delta}$ as $\delta\to 0$.
	\end{enumerate} 
To proceed, we fix the number of points at the proxy-surface to $N=686$ or $N=1353$ (the largest number of precomputed spherical design nodes smaller than $N=700$ and $N=1400$ respectively). 
% The singular values of the self-interaction submatrices of the target-from-source matrix are all above $\epsilon_{\text{mach}}$ 
The ratio of the smallest to largest singular values of the self-interaction target-from-source matrix are all above $\epsilon_{\text{mach}}$ and no truncation is therefore in use in the one-body preconditioning scheme of Section \ref{solving}. The resulting $\epsilon_{\text{res}}^{\max}$, as given by \eqref{residual},  are visualized in Figures \ref{Rg_sweep_a}, \ref{Rg_sweep_b} and \ref{Rg_sweep_e} as functions of $R_p$ and $\delta$. Shearing motion in panel \ref{Rg_sweep_a} is easier to resolve than squeezing motion in \ref{Rg_sweep_b} and \ref{Rg_sweep_e}, and a difference in error levels can also be noted between the two grid resolutions, represented in panels \ref{Rg_sweep_b} and \ref{Rg_sweep_e} for the same test problem. This is as predicted by MFS theory: with a larger number of source points, aliasing errors, dominating for larger $R_p$, are decreased. Figure \ref{Rg_sweep_d} visualizes $\|\vec\lambda\|_{\infty}$  determined with the finer grid with $N=1353$. From panels \ref{Rg_sweep_a}, \ref{Rg_sweep_b}, \ref{Rg_sweep_d} and \ref{Rg_sweep_e}, note that a large $R_p$ corresponds to a small $\delta^*$, while a smaller $R_p$ generally gives a smaller relative residual $\epsilon_{\text{res}}^{\max}$ for a fixed particle separation $\delta$, however at the expense of a larger coefficient vector (as predicted by MFS theory).
%if singularities of the analytic continuation of the Stokes solution into the particle interior are not enclosed by the proxy-surface, accuracy can still be achieved, but with exponentially growing coefficient norm \cite{Barnett2007,Doicu2000}.
Note however that, due to round-off, such big coefficient norms typically cause loss of accuracy. 
The error in the computed force on one of the particles is shown in Figures \ref{Rg_sweep_c} and \ref{Rg_sweep_f}.
%In Figures \ref{Rg_sweep_c} and \ref{Rg_sweep_f}, $\epsilon_{\text{FT}}$ compares the computed force on one of the particles  to the result %in \eqref{Brenner} for the two different choices of $N$. 
Note that %$\epsilon_{\text{force}}<\epsilon_{\text{res}}$, %
the relative error in the force is generally smaller than the maximum relative residual, 
often by a few orders of magnitude. 

%
%For the coarser grid and given the target accuracy of $\text{TOL}=10^{-3}$, we pick $R_p = 0.59$. With $N=1400$, $R_p$ can be chosen larger for the same target accuracy, $R_p = 0.665$ and hence, $\delta^*$ is smaller. For $N=700$, $\delta^*\approx0.29$, while for $N=1400$, $\delta^*\approx0.17$.  
For the coarser grid and given a general target accuracy of $10^{-3}$, we pick $R_p = 0.63$, for which $\delta^*\approx 0.271 $. With $N=1353$, $R_p$ can be chosen larger for the same target accuracy. We set $R_p = 0.7$, where $\delta^* \approx 0.129 $. The coarse grid is used throughout the numerical experiments in this paper, while the dense grid will serve as a reference in some of the tests. 
By extrapolating from the data in Figure \ref{Rg_sweep}, it can be concluded that to resolve the interaction for particles with relative motion and a small separation distance $\delta \ll 0.1$, the number of sources on the proxy-surfaces would have to be intractably large. For such interactions, we will instead add  image points as given by \eqref{line_points}-\eqref{tj}. 
%The image line goes between the sphere center and the accumulation point  \eqref{acc_points}, and %image points can be placed between the proxy surface and the accumulation point, if the %accumulation point is outside of the proxy surface, i.e. whenever $\delta<\delta^*$.

% Note that it is possible to place such image points along 
% %Note that such image points \emph{can} always be placed on 
% the image line between the proxy-surface and the accumulation point, as defined in \eqref{acc_points}, if the accumulation point is outside of the proxy-surface, i.e.~whenever $\delta<\delta^*$. 

%This is however a choice we do only if needed to resolve given boundary conditions. 
%See Section~\ref{s:close} for further details of adding image points.

%We will for now assume that the boundary condition with the particles traveling towards each other  is representative for any other no-slip condition -- this question will also be tested numerically in Section \ref{lub}. 

	%\subsection{Image points, fundamental solutions and additional collocation points}\label{image_locations}
	
	%%%%%%%%%%%%%%%%%%%%%%%%%%%%%%%%%%%%%%%%%%%%%%%%%%%%%%%%%%%%%%%%%%%%%%%%%%%% 
	
	%´	\section{Numerical experiments }\label{StokesNum}	
	\section{Large scale examples with moderate lubrication forces}\label{large_ex}

 In this section we show the performance of the basic MFS scheme using uniform proxy spheres and no lubrication-adapted images. This can reach the target 3-digit accuracy
if the particles are well-separated (Example 1 below), or have no relative motion (Example 2), or both.
 %With $N$ only a couple of times larger, $\delta=0.1$ can be reached.
 
 We consider two example geometries. With no images, every particle is discretized equally, which means that only one SVD is needed to compute the particle pseudo-inverse used repeatedly to build the preconditioned linear system in \eqref{resistance}. To evaluate $\vec u_\text{All}$ in \eqref{uAll},
 and to evaluate the fluid flow in the exterior domain in a post-processing step, linear scaling is achieved using the MATLAB interface to the Stokes kernels in the multithreaded FMM3D library \cite{ChengH.1999,FMM3D}. In both examples in this section, the GMRES tolerance is set to $10^{-7}$ and the FMM tolerance is set to $10^{-8}$.
The implementations are in MATLAB R2022b, using an 8-core workstation with Intel\circledR~Xeon\circledR~CPU E5-2620 @ 2.00GHz and 128 GB of RAM.

	%Hence, the scaling for solving the least squares problem for the MFS coefficients $\vec \lambda$ and evaluating the fluid flow in the exterior domain for a given vector $\vec\lambda$ is linear in the number of particles.  
	\begin{example}\label{ex1}
	First, consider a dilute suspension of spheres placed randomly in a layer, with geometry as visualized in Figure \ref{layer}.  With the discretization choice $N = 686$ and $R_p=0.63$ from Section \ref{param_sec}, the well-separated condition means $\delta>\delta^*\approx0.271$, recalling \eqref{deltastar}.  To generate the suspension geometry, the position of each new sphere is drawn from a uniform distribution in a box of dimension $[L,L,2]$  until the distance to any neighbor is at least $\delta = \delta ^*$. 
The boundary conditions are generated by random angular and translational velocities for every particle. The computed flow magnitude in a plane cutting through an example suspension containing 2000 spheres is displayed in Figure \ref{flow_field}.   %$R_p = 0.63$, 
	A maximum relative residual  $\epsilon_{\text{res}}^{\max}$, recalling \eqref{residual}, %,as defined by \eqref{residual}, and 
 matching the target tolerance of $10^{-3}$ is obtained, as visualized in Figure \ref{Fig:random_error}. To determine $\epsilon_{\text{res}}^{\max}$, the residual is evaluated in 500 points per surface, different from the set of collocation points used to enforce boundary conditions. The force/torque error $\epsilon_{\text{FT}}$, determined relative to a reference computed with $N=1353$ via \eqref{force_err}, is nearly three digits more accurate than $\epsilon_{\text{res}}^{\max}$.  The cost of FMM is linear in the number of particles $P$, as illustrated in Figure \ref{Fig:time_random}, but the number of GMRES iterations grows weakly, which is the reason why the reported total solution time in Figure \ref{Fig:time_random} is slightly superlinear. The solution time is roughly 10 s per particle.
	For comparison, the time to compute the single particle pseudo-inverse with an SVD is 2.5 s. 
	\begin{figure}[h!]
		\centering
		%			\begin{subfigure}[t]{0.44\textwidth}
		%				\centering
		%				\hspace*{-2ex}
		%				\includegraphics[trim={0cm 0cm 0cm 0cm},clip,width=1.1\textwidth]{../figures/MFS/500_spheres_plane_flow.eps}
		%				%\vspace*{0.2ex}
		%				\caption{Suspension geometry.}
		%				\label{layer}
		%			\end{subfigure}\\
		% HERE IS THE SUSPENSION PLOT
		%					\begin{subfigure}[t]{0.65\textwidth}
		%			\centering
		%			\hspace*{-2ex}
		%			\includegraphics[trim={1cm 3cm 1.5cm 3cm},clip,width=1\textwidth]{../figures/MFS/config2000.eps}
		%			%\vspace*{0.2ex}
		%			\caption{Suspension geometry with 2000 spheres.}
		%			\label{layer}
		%		\end{subfigure}\\
		
		%			\begin{subfigure}[t]{0.44\textwidth}
		%				\centering
		%				%\hspace*{-1ex}
		%				\includegraphics[trim={1.5cm 3cm 3cm 3cm},clip,width=1.1\textwidth]{../figures/MFS/500_spheres_plane2.eps}
		%				%\vspace*{0.2ex}
		%				\caption{Flow magnitude in a plane.}
		%				\label{flow_field}
		%			\end{subfigure}	
		%			\begin{subfigure}[t]{0.09\textwidth}
		%				\vspace*{-30ex}
		%				\centering
		%				\hspace*{1ex}
		%				\includegraphics[trim={12cm 1cm 0cm 0cm},clip,width=1\textwidth]{../figures/MFS/500_spheres_plane2.eps}
		%				%	\caption{Flow magnitude in plane.}
		%				%	\label{imagelines}
		%			\end{subfigure}
		\hspace*{-4ex}
		\begin{subfigure}[t]{0.64\textwidth}
			\begin{minipage}[b]{\textwidth}
				\centering
				\hspace*{-2ex}
				\includegraphics[trim={1cm 3cm 1cm 3cm},clip,width=1\textwidth]{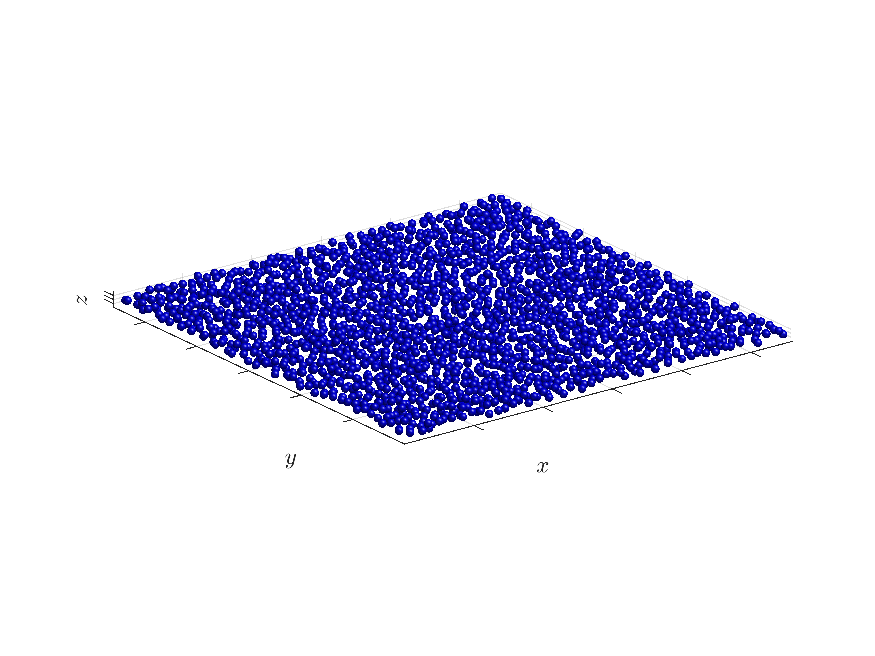}
				%\vspace*{0.2ex}
				\caption{Suspension geometry with 2000 spheres.}
				\label{layer}
			\end{minipage}
			\begin{minipage}[b]{\textwidth}
				\centering	
    \begin{overpic}[trim={2.7cm 0.7cm 0.5cm 0.5cm},clip,width=1.02\textwidth]{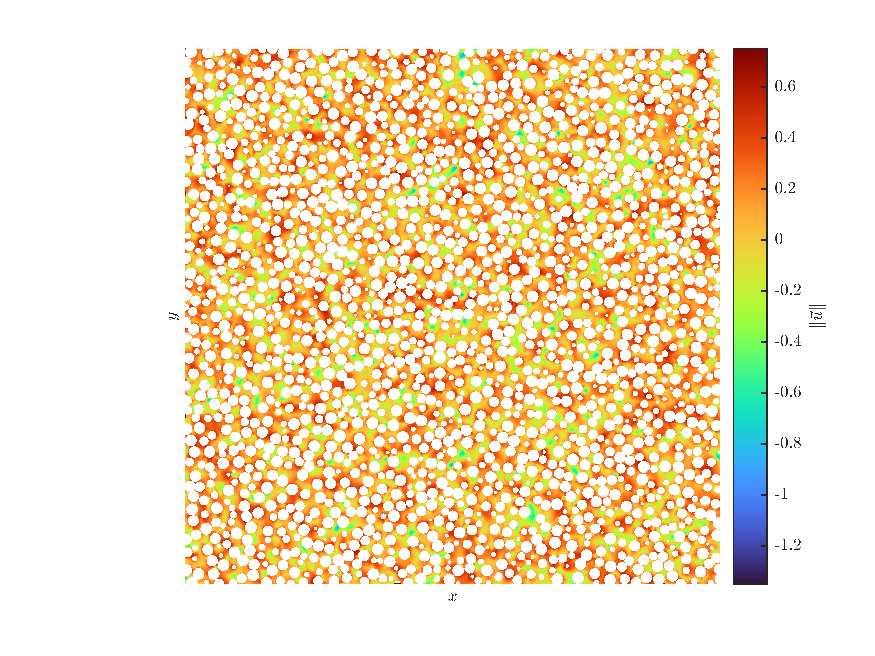}
\put(94,35){\rotatebox{90}{\colorbox{white}{\parbox{0.1\linewidth}{%
     \small$\log_{10}\left(\|\vec u\|_2\right)$}}}}
     \end{overpic}
				%\includegraphics[trim={2.7cm 0.7cm 0.5cm 0.5cm},clip,width=1.02\textwidth]{figures/velocity2000.eps}
				%\vspace*{0.2ex}
				\caption{Flow magnitude in a plane cutting through the configuration of spheres in (a). White space is occupied by particles.}
				\label{flow_field}
			\end{minipage}
		\end{subfigure}~~
		\begin{subfigure}[t]{0.35\textwidth}
			%\vspace*{-56ex}
			\vspace*{-32ex}
			\begin{minipage}[b]{\textwidth}
				\hspace*{-3ex}
				\includegraphics[trim = {5cm 15.8cm 7.7cm 4cm},clip,width=1.08\textwidth]{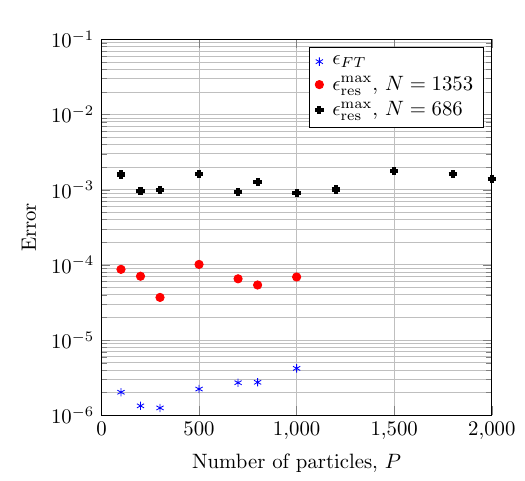}
				\caption{The error $\epsilon_{FT}$ in the force/torque vector determined with $N=686$ sources per particle is computed relative to a fine grid reference determined with $N=1353$. Displayed are also $\epsilon_{\text{res}}^{\max}$ for both choices of $N$.}
				\label{Fig:random_error}
			\end{minipage}
			\begin{minipage}[b]{\textwidth}
				\hspace*{-1ex}
				\includegraphics[trim = {5cm 14.5cm 7.7cm 4.0cm},clip,width=1.1\textwidth]{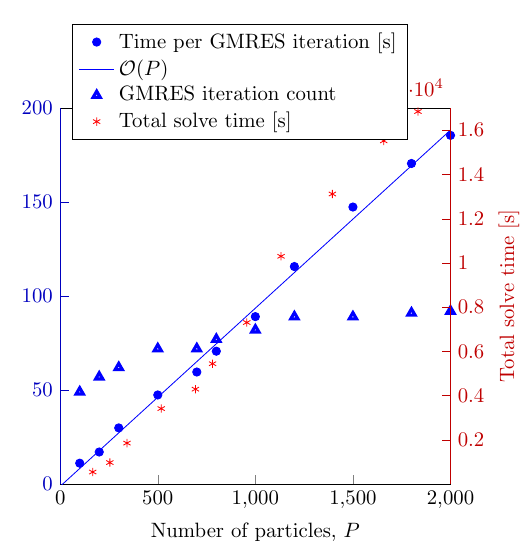}
				\caption{The computational cost per GMRES iteration scales linearly with the number of particles, but the iteration count grows weakly, giving a slightly worse than linear total solution time.}
				\label{Fig:time_random}
			\end{minipage}
		\end{subfigure}
		%\begin{subfigure}[t]{0.09\textwidth}
		%	\vspace*{-30ex}
		%	\centering
		%	\hspace*{1ex}
		%	\includegraphics[trim={12cm 1cm 0cm 0cm},clip,width=1\textwidth]{../figures/MFS/500_spheres_plane2.eps}
		%	%	\caption{Flow magnitude in plane.}
		%	%	\label{imagelines}
		%\end{subfigure}
		\caption{Example with a dilute suspension of spheres in a layer.}
	\end{figure}	
	%\subsubsection{A cluster of fixed spheres in a background flow}
%		The test is run on a 16-core 
	%		workstation with Intel(R) Xeon(R) CPU E5-2637 v3 @ 3.50GHz. %was on eight - updata CPU-timings for systems of the same density but different particle numbers are reported in Figure \ref{timings}. Note the linear scaling in the number of particles.
	%		\begin{figure}[h!]
	%			\centering
	%			\includegraphics[trim={0cm 0cm 0cm 0cm},clip,width=0.6\textwidth]{../figures/images/timings.eps}
	%			\caption{CPU timings}
	%			\label{timings}
	%			\end{figure}
	\end{example}
	\begin{example}
	Next, consider a cluster of fixed spheres in a background flow, as exemplified in Figure \ref{shear} for a shear flow and in Figure \ref{uniform} for a uniform flow. The cluster is grown so that each sphere is a distance of $\delta$ away from at least one other sphere, where $\delta = 10^{-2}$ or $\delta = 10^{-3}$. The spheres are fixed in space and do not move with the fluid, meaning that the \emph{total} velocity field at the sphere surfaces is zero. The boundary condition for the disturbance field due to the presence of the particles in \eqref{bc} in shear flow can be expressed as 
	\begin{equation}
		\vec u_{\text{bc}}(\vec x) = -\left(\vec E_{\infty}\vec x+\vec \omega_{\infty}\times \vec x\right),\,\quad\text{with }\quad \vec E_{\infty} = 0.5  \begin{bmatrix}
			0 & \dot{\gamma} & 0 \\ \dot{\gamma} & 0 & 0 \\ 0 & 0 & 0	
		\end{bmatrix},\,
		\quad \vec \omega_{\infty} = \begin{bmatrix}
			0 \\ 0 \\ -0.5\dot{\gamma}
		\end{bmatrix}, 
	\end{equation}
	and $\dot{\gamma} = 5$ the shear rate. For the uniform flow, we have
 \begin{equation}
    \vec u_{\text{bc}}(\vec x) = -\vec u_{\infty}, \quad \text{ with } \quad \vec u_{\infty} = \begin{bmatrix}
			5 & 0 & 0
		\end{bmatrix}^T. %\quad\xi = 5.
 \end{equation}
	An $\epsilon_{\text{res}}^{\max}$ smaller than the target accuracy of $10^{-3}$ can be obtained independently of the choice of $\delta$ or the number of particles for both flows, presented in Table \ref{cluster_table} for the shear case. 
 %\begin{figure}[h!]
 %\begin{minipage}[t]{0.6\textwidth}
		\begin{table}[h!]
		\centering
    \caption{Number of GMRES iterations to convergence, maximum relative residuals $\epsilon_{\text{res}}^{\max}$ and maximum coefficient magnitudes $\|\vec\lambda\|_{\infty}$ for simulations of $P$ particles in a shear flow (Figure \ref{shear}), using $N=686$. The GMRES tolerance is set to $10^{-6}$ and the FMM tolerance is $10^{-8}$.} 
		\label{cluster_table}
      %\vspace*{-37.5ex}   
\begin{tabular}{c|c|c|c}
			\textbf{Setting}  	& \multicolumn{3}{c}{\textbf{Quantity}} \\ \hline
			& \textbf{GMRES iters.} & $\vec \epsilon_{\text{res}}^{\max}$& $\|\vec\lambda\|_{\infty}$ \\ \hline \hline
			$\vec{P = 100}$ & \multicolumn{3}{c}{} \\ \hline 
   $\delta = 10^{-2}$, $R_p = 0.75$ & 92 & $1.1\cdot 10^{-4}$ & $9.83$\\ \hline 
    $\delta = 10^{-3}$, $R_p = 0.75$    & 122 & $9.6\cdot 10^{-5}$& $9.52$ \\ \hline
     $\delta = 10^{-3}$, $R_p = 0.65$  &  116 & $6.9\cdot 10^{-5}$& $316$ \\ \hline \hline
    $\vec{P = 2000}$ & \multicolumn{3}{c}{} \\ \hline $\delta = 10^{-2}$, $R_p = 0.75$    &   144 & $9.3 \times 10^{-5}$ &  $18.3$\\ \hline 
 $\delta = 10^{-3}$, $R_p = 0.75$    &     176 & $1.3 \times 10^{-4}$ & $18.9$ \\ \hline 
    $\delta = 10^{-3}$, $R_p = 0.65$    &     156 & $1.5 \times 10^{-4}$ &  $1013$\\ \hline 

		\end{tabular}
  %\vspace*{9.5ex}
  %\captionof{table}{

  \end{table}
 % \end{table}
 % \end{minipage}~~~
  % \begin{minipage}[t]{0.39\textwidth}
  %     \centering
		% 	\includegraphics[trim={5cm 15.8cm 8.5cm 4cm},clip,width=1\textwidth]{figures/largecluster_stats.pdf}
  %  \caption{Split of the CPU time for simulations of the type in Figure \ref{large_cluster}. 
  %  %The FMM is more dominant for larger systems. 
  %  The total CPU times for $P=2000$ are 12.8 h and 10.1 h for shear and uniform flow, while $P=100$ runs in 0.51 h and 0.45 h, respectively.}
  %  \label{CPUtime}
  % \end{minipage}
 % \vspace*{-2ex}
 % \end{figure}
  \begin{figure}[h!]
		\centering
		\begin{subfigure}[b]{0.49\textwidth}
  \hspace*{-3ex}
  \begin{minipage}{0.8\textwidth}
			\centering
			\includegraphics[trim={0cm 0cm 3cm 0cm},clip,width=1.1\textwidth]{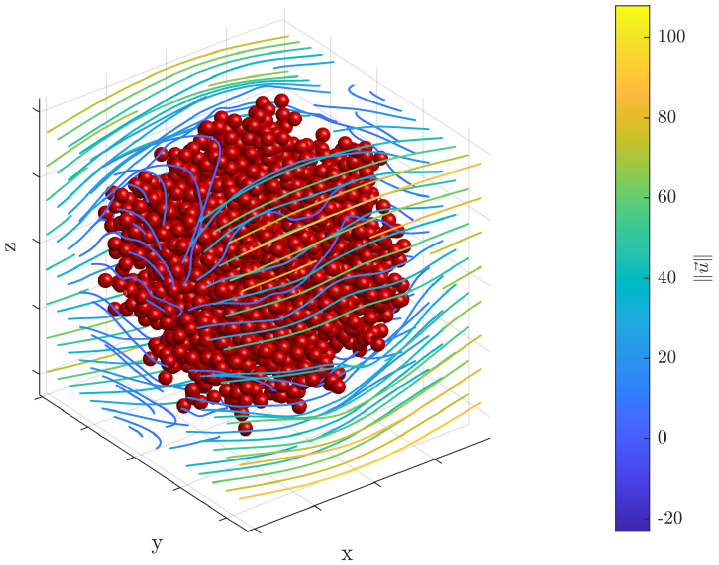}
			
    \end{minipage}
    \begin{minipage}{0.15\textwidth}
			\centering		%\includegraphics[trim={10cm 0cm 0cm 0cm},clip,width=1.25\textwidth]{figures/cluster4.pdf}	
   \begin{overpic}[trim={10cm 0cm 0.5cm 0cm},clip,width=1\textwidth]{figures/cluster4.pdf}	
\put(18,48){\rotatebox{90}{\colorbox{white}{\parbox{0.5\linewidth}{%
     \small$\|\vec u\|_2$ }}}}
     \end{overpic}
    \end{minipage}
    	\caption{Shear flow}
     \label{shear}
		\end{subfigure}
  %\vspace*{10ex}
  		\begin{subfigure}[b]{0.49\textwidth}
  \hspace*{-2.5ex}
    \begin{minipage}
    {0.8\textwidth}
      \vspace*{-15ex}
			\centering
			 \includegraphics[trim={1.5cm 0.8cm 4.2cm 1cm},clip,width=1.04\textwidth]{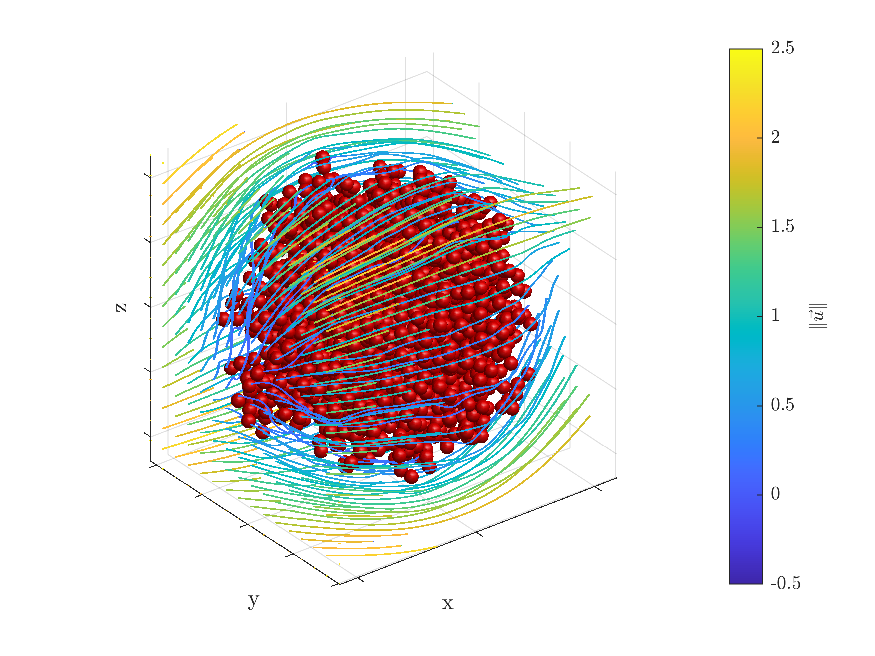}		
   \end{minipage}~~~
   \begin{minipage}
    {0.15\textwidth}
      \vspace*{-15ex}
			\centering
       \begin{overpic}[trim={11.7cm 0cm 1.2cm 0.2cm},clip,width=1.15\textwidth]{figures/cluster_uniform.eps}
\put(18,48){\rotatebox{90}{\colorbox{white}{\parbox{0.1\linewidth}{%
     \small$\|\vec u\|_2$}}}}
     \end{overpic}
			%\includegraphics[trim={11.4cm 0cm 0cm 0cm},clip,width=2\textwidth]{figures/cluster_uniform.eps}
			%\caption{Uniform flow}
   \end{minipage}
   \caption{Uniform flow}
   \label{uniform}
		\end{subfigure}
		%			\begin{subfigure}[b]{0.1\textwidth}
		%				%	\vspace*{-2ex}
		%				\centering
		%				\hspace*{-23ex}
		%				\includegraphics[trim={9.5cm 0cm 0cm 0cm},clip,width=0.8\textwidth]{../figures/images/cluster.eps}
		%				%	\caption{}
		%				%	\label{imagelines}
		%			\end{subfigure}
		\caption{Streamlines around clusters of 2000 fixed spheres in a background flow. Each sphere is separated by $\delta = 10^{-3}$ from at least one neighbor  and discretized with $N=686$ proxy-sources, resulting in a total number of $4.8$ million degrees of freedom in the preconditioned linear system. %4.806 
  No image points are needed as the particles do not move relative to each other.
  The solutions have about 4-digit accuracy (uniform boundary condition error), and take around 10 hours on one workstation.}
  \label{large_cluster}
	\end{figure}
 The results for the uniform flow are similar, but with slightly faster convergence. The FMM dominates the cost with about 90\% of the total solution time, and is more dominant for larger systems. For the largest suspensions with 2000 particles, the preconditioned systems solved via FMM have $4.81\times 10^6$ degrees of freedom. In comparison, applying the pseudo-inverse and computing the self-interaction (steps 1 and 3 in the algorithm described in Section \ref{solving}) each takes about 5\% of the solution time.  
  The total CPU times for $P=2000$ are 12.8 h and 10.1 h for shear and uniform flow, while $P=100$ runs in 0.51 h and 0.45 h, respectively.  
  Note in Table \ref{cluster_table} that the smallest particle-particle distance, $\delta = 10^{-3}$, is, surprisingly, no harder to resolve than the larger particle-particle distance $\delta = 10^{-2}$. This is a consequence of the fact that the force density---the coefficient vector $\vec\lambda$---is benign, despite the closeness of the particles; see the reported $\|\vec\lambda\|_{\infty}$ in the table.  This is in contrast to the prior test in Figure \ref{Rg_sweep}, where two particles travel relative to each other and $\|\vec\lambda\|_{\infty}$ gets very large for small $\delta$ and any relevant $R_p$. This contrast matches the observations in \cite{Lefebvre-Lepot2015}, who attribute the difficulty of resolving a peaked force density to problems with relative motion between particles. af Klinteberg \& Tornberg \cite{AfKlinteberg2016} in their study of 131 ellipsoids using QBX also implicitly make use of the smoothness of the force density for static particles.

  Despite this apparent lack of singularities at $R_{acc}$, we still find that decreasing $R_p$ results in a larger coefficient norm, as shown by comparing lines 2--3 or 5--6 of Table \ref{cluster_table}. The required number of GMRES iterations for clusters of two different sizes are also compared in this table. Similarly as in Example \ref{ex1}, and for different choices of $R_p$ and $\delta$, the number of iterations grows only weakly with $P$. 
	
Streamlines illustrated in Figure \ref{large_cluster} are computed by propagating a set of tracer particles in the fluid, using \textsc{Matlab}'s \texttt{ode45} to solve the ODE 
	\begin{equation}
		\frac{\mathrm d\vec x_t}{\mathrm d t} = \vec G(\vec x_t)\vec\lambda
	\end{equation}
 for the tracer particle locations $\vec x_t$ at time $t$.
	The MFS coefficient vector $\vec\lambda$ is computed once for the fixed particle geometry.  %\noteme{Residual in the uniform flow is $1.44\times 10^{-4}$.}
 %\begin{figure}
 
 % \end{figure}

\end{example}

%%%%%%%%%%%%%%%%%%%%%%%%%%%%%%%%%%%%%%%%
	\section{Closely interacting particles under relative motion}
 \label{s:close}

 We now turn to Stokes sphere problems with strong lubrication forces, arising for small separations with relative motions. We show that they are efficiently solved by augmenting the coarse proxy sphere with fundamental solutions placed along discretized image lines defined by \eqref{line_points}-\eqref{tj}. We start by studying which singularity types are needed along these lines, together with a strategy for choosing additional collocation points. The result is an adaptive strategy to set the number of image and collocation points for each pair-wise separation $\delta$. This is demonstrated on clusters of spheres in Section \ref{Num_res}. The coarse grid with $N = 686$ and $R_p = 0.63$ (from Section \ref{param_sec}) is used throughout for the proxy spheres.  Finally, in Section \ref{Gmres_count}, we discuss the GMRES iteration count for small $\delta$.% and test left-preconditioning (row scaling).
  
  One-body right preconditioning (as in Section \ref{solving}) can still be applied with the extra images, but with the difference that the number of vector-valued sources (and collocation nodes) now can vary for each particle, depending on the proximity to its neighbors. The self-interaction blocks of the target-from-source matrix will therefore be different for each particle, increasing the SVD costs. The clustered image and collocation points also cause more ill-conditioning of the self-interaction blocks, as their columns are nearly linearly dependent. %In addition, aliasing errors are larger with sources moved close to the particle surface \cite{Barnett2007}. Clustering of the additional collocation points reduces these problems to some extent.
  This means that truncation in the pseudo-inverses is needed: we set $\epsilon_{\text{trunc}} = 5\cdot 10^{-12}$ throughout, via a tradeoff between residual error vs coefficient norm (see Section \ref{sec_trunc} in the Appendix).  Moreover, since it slightly improves the conditioning of the target-from-source matrix, left-preconditioning (row-scaling) has also been used in the paper whenever image points are utilized, discussed in Remark \ref{left-precond}. The decay of the singular values with preconditioning is exemplified in Figure \ref{sing_decayb} in the Appendix, to be compared to the singular values without preconditioning shown in Figure \ref{sing_decay}.
%  To reduce the effect of ill-conditioning further, truncation of the self-interaction blocks is however necessary. % for the right-preconditioning strategy.  
  %High truncation generally results in larger residuals but smaller coefficient magnitudes.

	%	An example for the Stokes exterior BVP is outlined in Figure \ref{StokesnD}, where the relative error in the flow field is computed in the exterior of two circles in 2D and two spheres in 3D, as compared to a reference solution computed with a boundary integral method. 
	%	
	%	\begin{figure}[h!]
	%		\centering
	%		\begin{subfigure}[t]{0.49\textwidth}
	%			\centering
	%			\includegraphics[trim={3cm 2cm 0cm 1cm},clip,width=\textwidth]{../figures/images/error2discs.eps}
	%			\caption{2D}
	%		\end{subfigure}
	%		\begin{subfigure}[t]{0.49\textwidth}
	%			\centering
	%			%			\includegraphics[trim={2cm 19cm 10cm 1cm},clip,width=\textwidth]{../figures/MFS/single_sphere_Rg_stokes2b_theta001_a1.pdf}
	%			\caption{3D}
	%			%	\label{imagespheres_Stokes}
	%		\end{subfigure}
	%		\caption{Exterior Stokes BVPs are solved using the grids in Figure \ref{image_figs}, with no-slip boundary conditions enforced in each target point. Relative errors are displayed obtained with the method of fundamental solutions of this paper, as compared to reference solutions computed with boundary integral methods. The error is in the order of at most $10^{-5}$ everywhere in the domain, even with coarse representations of $\Gamma_p$.}
	%		\label{StokesnD}
	%		
	%	\end{figure}
 %%%%%%%%%%%%%%%%%%%%%%%%%%%%%
 	\subsection{Choice of fundamental solutions and extra collocation points}\label{lub}
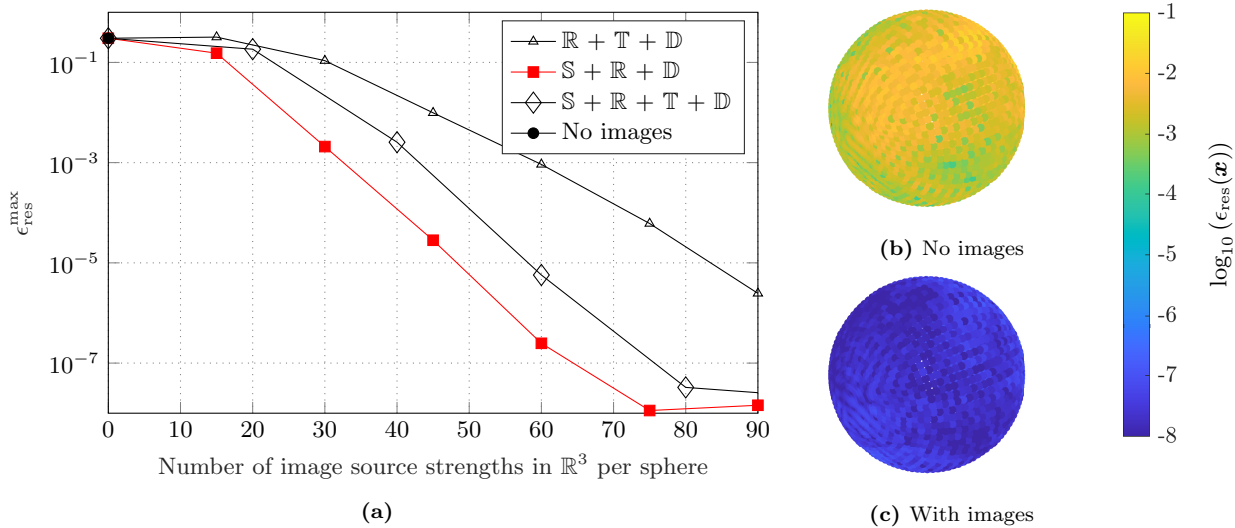
\begin{figure}[h!]
		\centering	
  	\begin{subfigure}[t]{0.6\textwidth}
   \begin{small}
			\begin{tikzpicture}		
				\begin{axis}[%
					width=3.4in,
					height=2.1in,
					at={(0.758in,0.481in)},
					scale only axis,
					xmin=0,
					xmax=90,
					xlabel style={font=\color{white!15!black}},
					xlabel={Number of image source strengths in $\mathbb R^3$ per sphere},
					ymode=log,
					ymin=1e-08,
					ymax=1,
					yminorticks=true,
					ylabel style={font=\color{white!15!black}},
					ylabel={$\epsilon_{\text{res}}^{\max}$},
					axis background/.style={fill=white},
					xmajorgrids,
					ymajorgrids,
					%yminorgrids,
					legend style={legend cell align=left, align=left, draw=white!15!black}
					]
					
					\addplot [color=black, mark=triangle, mark options={solid, black}, mark size = 2pt]
					table[row sep=crcr]{%
						%	0	0.217259017194829\\
						0	0.303475903538758\\
						15	0.31988613728762\\
						30	0.108455301218528\\
						45	0.00986611135581891\\
						60	0.000921045858516089\\
						75	6.09585440954213e-05\\
						90	2.46057150066078e-06\\
					};
					%0    15    30    45    60    75    90
					\addlegendentry{$\mathbb R$ + $\mathbb T$ +  $\mathbb D$}
					
					\addplot [color=red, mark=square*, mark options={solid, red}]
					table[row sep=crcr]{%
						%	0	0.228510641957193\\
						0	0.303475903538758\\
						15	0.152788890003357\\
						30	0.00209132143728974\\
						45	2.82773468702262e-05\\
						60	2.48267106013834e-07\\
						75	1.13107607029974e-08\\
						90	1.4410285842714e-08\\
					};
					\addlegendentry{$\mathbb S$ + $\mathbb R$ +  $\mathbb D$}
					
					\addplot [color=black, mark=diamond, mark options={solid, black}, mark size = 4pt]
					table[row sep=crcr]{%
						%0	0.303475903538758\\
						0	0.303475903538758\\
						20	0.183554058346321\\
						40	0.002557929860229\\
						60	5.72165787636253e-06\\
						80	3.2733351988227e-08\\
						100	2.01189420912737e-08\\
						120	2.03929033896039e-08\\
					};
					\addlegendentry{$\mathbb S$ + $\mathbb R$ + $\mathbb T$ +  $\mathbb D$}
					\addplot [color=black, mark=*, mark options={solid, black}, mark size = 2pt]
					table[row sep=crcr]{%
						0	0.303475903538758\\
					};
					\addlegendentry{No images}
				\end{axis}
			\end{tikzpicture}
   \end{small}
   \vspace*{-3ex}
   \caption{}
			%\caption{One combination of fundamental solutions reduces the residual at the particle surface more than any other combination, namely, a linear combination of Stokeslets, rotlets and potential dipoles.}
   \label{decay_im}
		\end{subfigure}~~~~
		\begin{subfigure}[t]{0.25\textwidth}
			\vspace*{-43ex}
			\begin{minipage}[b]{\textwidth}
				\includegraphics[trim = {2.5cm 2.3cm 2.8cm 2cm},clip,width=\textwidth]{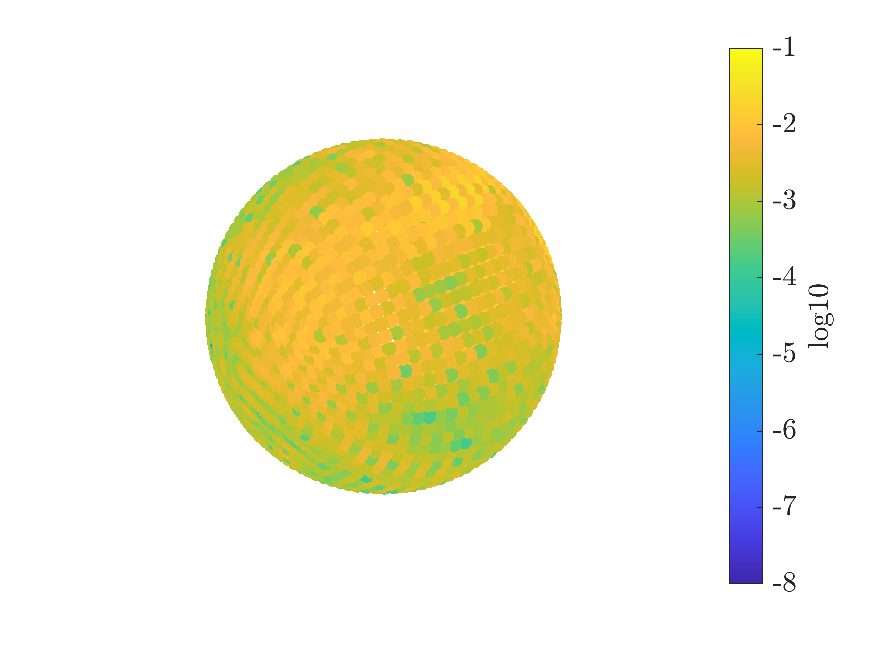}
				\caption{No images}
				\label{residuala}
			\end{minipage}
			\begin{minipage}[b]{\textwidth}
				\includegraphics[trim = {2.5cm 2.3cm 2.8cm 2cm},clip,width=\textwidth]{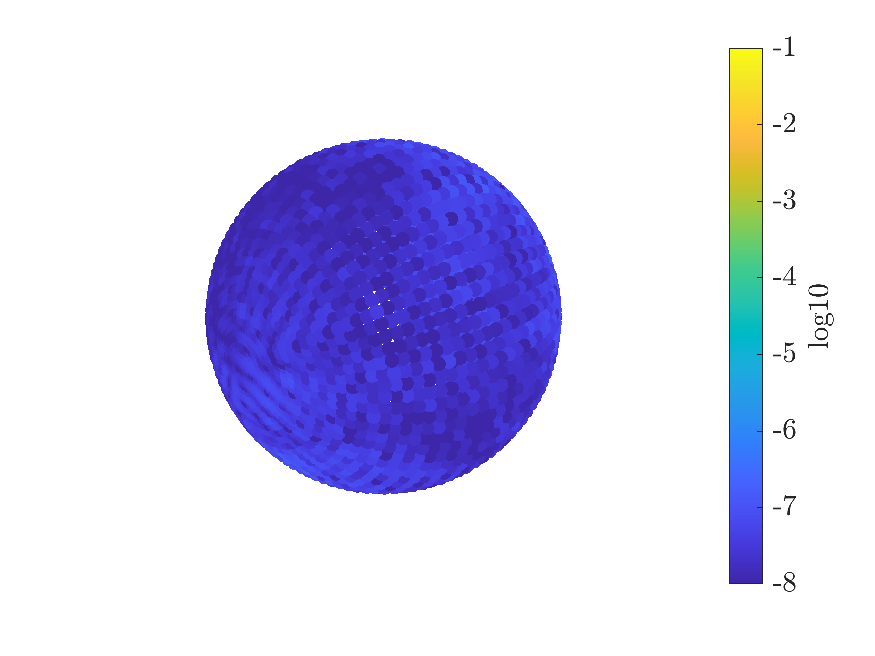}
				\caption{With images}
				\label{residualb}
			\end{minipage}
		\end{subfigure}
		\begin{subfigure}[t]{0.1\textwidth}
			\vspace*{-43ex}
   \begin{overpic}[trim = {12.2cm 0cm 1.2cm 0.5cm},clip,width=0.53\textwidth]{figures/with_images.eps}
% \put(17,16){\rotatebox{90}{\colorbox{white}{\parbox{\linewidth}{%
%      \small$\log_{10}\left(\|\vec u(\vec x)-\vec u_{\text{bc}}(\vec x)\|_{2}/\|\vec u_{\text{bc}}(\vec x)\|_{2}\right)$}}}}
     \put(17,40){\rotatebox{90}{\colorbox{white}{\parbox{\linewidth}{%
     \small$\log_{10}\left(\epsilon_{\text{res}}(\vec x)\right)$}}}}
     \end{overpic}
		\end{subfigure}

	%\vspace*{-3ex}
		\caption{Combinations of fundamental solutions along image lines reduce the maximum relative residual $\epsilon_{\text{res}}^{\max}$ on the particle surface. Panel (a) shows convergence for various combinations, for a pair of spheres separated by a distance of $\delta = 10^{-2}$. Panels (b) and (c) plot the relative residual, as defined in \eqref{residual_x}, over one of the particles, without and with the best combination of fundamental solutions at $n_{\text{im}}=20$ image points: the improvement is dramatic.  
		}
		\label{images_choice}
	\end{figure}

 We start by numerically testing different combinations of the four types of fundamental solutions at the approximate image points for a \emph{single} near-contact with separation $\delta$: Stokeslets $\mathbb S$, rotlets $\mathbb R$, stresslets $\mathbb T$ and potential dipoles $\mathbb D$.  Figure \ref{decay_im} compares the convergence with respect to the number of vector unknowns per image line, for a few such combinations, and shows that Stokeslets, rotlets and potential dipoles is the most efficient combination (shown in red).
 %The improvement of the maximum relative residual $\epsilon_{\text{res}}^{\max}$ defined in \eqref{residual} with an increasing total number of image source strengths (vectors in $\mathbb R^3$) is displayed.
 Here, a fixed large number of collocation points, $M=5000$, was used, sufficient to resolve lubrication forces.
 Figures \ref{residuala} and \ref{residualb} show that adding the converged image line to the coarse proxy sources reduces the residual by seven digits. 
 From now on, the number of image sources per particle per contact will be $3n_{\text{im}}$,
 accounting for three types of fundamental solutions at the $n_{\text{im}}$ image points $\vec l_j$ in \eqref{line_points}-\eqref{tj}. 
	
\begin{figure}[h!]
 \begin{subfigure}[t]{0.28\textwidth}
 %\begin{minipage}{0.82\textwidth}
 \hspace*{-3ex}
			\includegraphics[trim={1.5cm 0cm 4cm 0.5cm},clip,width=1.15\textwidth]{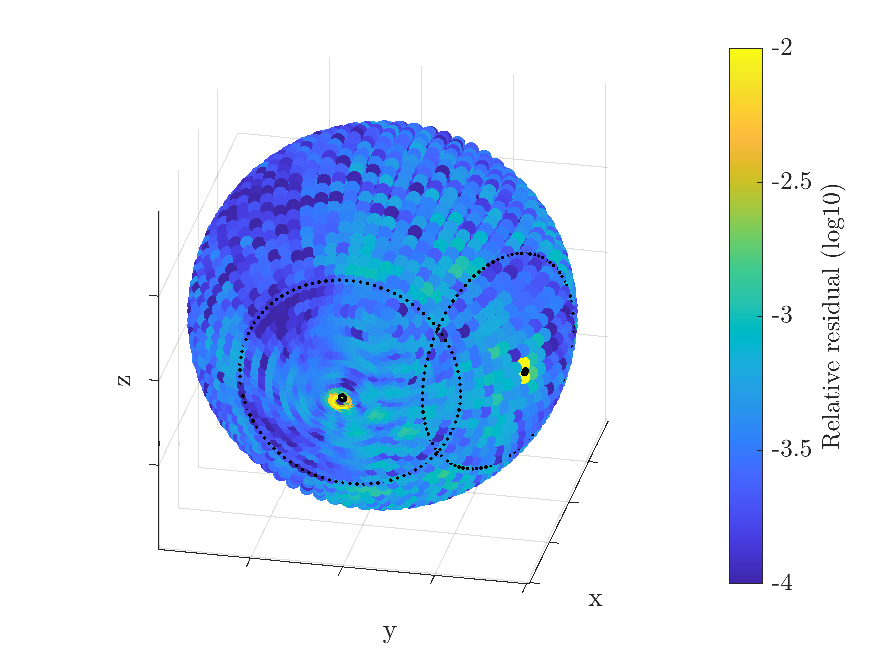}				
   %\end{minipage}
   \caption{$\beta_2 = \pi/200$}%: The angle of the second cap is chosen too narrow for the density to be sufficiently high in the region above the image lines. }
		\end{subfigure}~~
   \begin{subfigure}[t]{0.28\textwidth}
    \hspace*{-3ex}
			\includegraphics[trim={1.5cm 0cm 4cm 0.5cm},clip,width=1.2\textwidth]{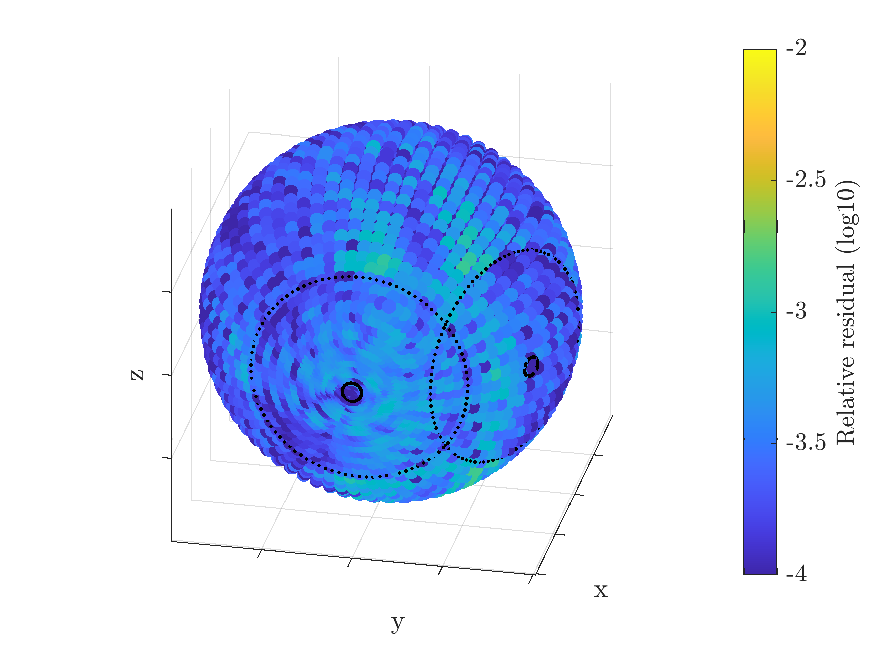}	
			\caption{$\beta_2 = \pi/60$}%: The maximum relative residual is below $10^{-3}$ everywhere on the particle surface.}
		\end{subfigure}~~
     \begin{subfigure}[t]{0.28\textwidth}
      \hspace*{-3ex}
			\includegraphics[trim={1.5cm 0cm 4cm 0.5cm},clip,width=1.2\textwidth]{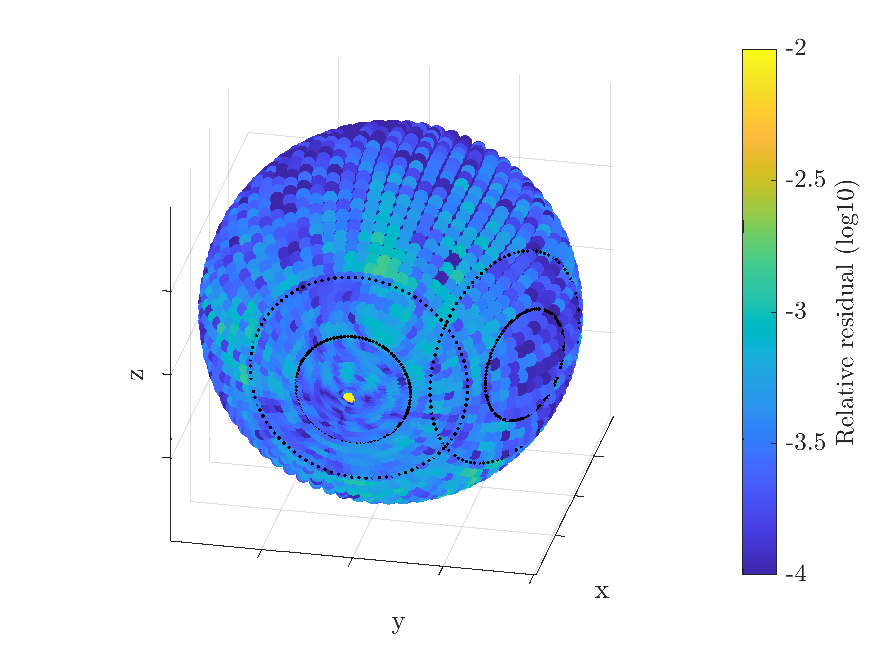}	
			\caption{$\beta_2 = \pi/10$}%:  The angle of the second cap is chosen too wide for the density to be sufficiently high in the region above the image lines. }
		\end{subfigure}~~~
  \vspace*{10ex}
      \begin{subfigure}{0.1\textwidth}
      \begin{overpic}[trim={11.2cm 1cm 1cm 0.5cm},clip,width=1.05\textwidth]{figures/capeffect_pi200.eps}	
% \put(17,16){\rotatebox{90}{\colorbox{white}{\parbox{\linewidth}{%
%      \small$\log_{10}\left(\|\vec u(\vec x)-\vec u_{\text{bc}}(\vec x)\|_{2}/\|\vec u_{\text{bc}}(\vec x)\|_{2}\right)$}}}}
     \put(26,35){\rotatebox{90}{\colorbox{white}{\parbox{\linewidth}{%
     \small$\log_{10}\left(\epsilon_{\text{res}}(\vec x)\right)$}}}}
     \end{overpic}
   \end{subfigure}
   \vspace*{-8ex}
		\caption{
 The dependence on the density of extra collocation points for the relative surface residual. Collocation points are added on two spherical caps per close contact (see 
  Figure \ref{colloc_points}) and three different choices of $\beta_2$ in \eqref{eq:thetacap} are investigated, with $\beta_1$ fixed at $\pi/5$. The sphere has two close neighbors ($n_c = 2$) at a distance $\delta = 10^{-3}$ in an equilateral triangle. The borders of the two additional spherical caps are shown with black lines. Here,  $M^{\text{cap}}_1=M^{\text{cap}}_2$ are fixed following \eqref{colloc_nbr}, with $\alpha_1 = \alpha_2 = 6$.
%Parameter study to set extra collocation points for closely interacting spheres, with two additional spherical caps of collocation points added for each near-contact, as in Figure \ref{colloc_points}.    The borders of these regions with angles $\beta_1$, $\beta_2$ are illustrated in black on one of the spheres in an equilateral triangle of particles in  each of the panels (a)-(c). Here, $\beta_1 = \pi/5$ and $M^{\text{cap}}_1=M^{\text{cap}}_2$ are fixed following \eqref{colloc_nbr}, with $\alpha_1 = \alpha_2 = 6$. 
  The density of collocation points has to be sufficiently high in the region above the image lines to resolve the boundary condition. This is the case in panel (b), where the relative residual $\epsilon_{\text{res}}(\vec x)$ defined in \eqref{residual_x} is below the target accuracy of $10^{-3}$, whereas in panel (a), $\beta_2$ is chosen too narrow and in panel (c), $\beta_2$ is chosen too wide. In the example, $N = 686$, $R_p = 0.6$, $n_{\text{im}}=30$ and $M=1.2\cdot N+2\cdot 6\cdot n_c\cdot 3n_{\text{im}}$.}
		\label{sweep_cone}
	\end{figure}

	%In addition to the proper combination of fundamental solutions, a strategy is needed to  distribute additional collocation points on the surfaces of the pair in which to enforce boundary conditions. As it is too costly to increase the density of collocation points globally over the surface of each sphere, the refinement has to be done locally to better resolve the peaked force densities. We choose to cluster the additional points near the points of closest approach between particles. More specifically,
 In order to minimize the number of collocation nodes $M$ while still resolving the solution, we add spherical caps of uniformly sampled collocation points ``above'' the line of images, as shown back in Figure \ref{colloc_points}(b-c).  For each cap, two parameters have to be set: the number of points and the angle of the cap, together determining an average node density. For robustness in configurations with more than two close-to-touching particles, we have found that the best results are obtained by superposing two such caps per particle per contact, one of which has a smaller opening angle $\beta_1$, and the other a wider angle $\beta_2$. The numbers of collocation points added by the two caps on one particle involved in a single near-contact are
 \begin{equation}\label{colloc_nbr}
     M^{\text{cap}}_1 = \alpha_1\cdot 3n_{\text{im}}\qquad \text{ and } \quad
     M^{\text{cap}}_2 = \alpha_2\cdot 3n_{\text{im}},
 \end{equation} 
 with $\alpha_1$ and $\alpha_2$ the respective upsampling parameters to be chosen. For cap $k=1$ or $2$, we choose a quasi-uniform Fibonacci grid \cite{Marques2021} in the spherical coordinates $(\theta_j,\phi_j)$ relative to the inter-particle axis, namely
 \begin{equation}
\begin{aligned}
    \theta_j &= \arccos\left(1-j\cos{\beta_k}/M_k^{\text{cap}}\right), \label{eq:thetacap}\\
    \phi_j&=2j\pi\Phi^{-1}\mod{2\pi},\quad j = 1,\dots,M_k^{\text{cap}},
\end{aligned}
 \end{equation}
 with $\Phi=(1+\sqrt{5})/2$ the golden ratio.
 %In Figure \ref{sweep_cone}, we make a sweep over the parameters ($\alpha_1$, $\theta_1$) with fixed $\theta_2 = \pi/5$ and a few different values of $\alpha_2$ 
 The node density has to be sufficiently high in a sufficiently large area above the accumulation point to resolve the boundary data.
 As before, we target a uniform surface residual of $10^{-3}$.
 The effect of different choices of $\beta_1$ and $\beta_2$ is visualized in Figure \ref{sweep_cone} for a configuration of three spheres forming an equilateral triangle. Here $\delta = 10^{-3}$, being the most challenging interaction that we would like to resolve, and we set $n_{\text{im}} = 30$. Based on experiments of this type, and to keep the parameter selection simple, we fix for later use $\alpha_1 = \alpha_2 = 6$, $\beta_1 = \pi/5$ and $\beta_2 =\pi/60 $, independently of $\delta$.

 \begin{remark}[Left-preconditioning]\label{left-precond}
       %With our image enhanced MFS technique, we apply preconditioning from the left. 
       To improve on the conditioning with clustered source and collocation points, the following preconditioning is applied from the left. Motivated by a desire to preserve the $L_2$-norm of the solution \cite{Bremer2013}, let $\vec W$ be a diagonal matrix formed by the local area elements corresponding to each collocation point, with the matrix blocks $$\left[\sqrt{w_1}\vec I_{3\times 3},\sqrt{w_2}\vec I_{3\times 3},\dots,\sqrt{w_i}\vec I_{3\times 3},\dots, \sqrt{w_{MP}}\vec I_{3\times 3}\right]$$ on the diagonal and $w_i$ the local area associated with collocation point $i$. For a particle involved in one near-contact, $w_i$ takes one of three values, depending on if the collocation point $\vec x_i$ belongs to the basic grid or lies on one of the spherical caps of additional uniformly distributed collocation nodes.  The preconditioned problem takes the form  
	\begin{equation}\label{lq2}
 \vec W\vec G\boldsymbol\lambda =\vec W\vec b,
		%\min\limits_{\boldsymbol\lambda\in\mathbb R^{3NP}}\|\vec W\vec G\boldsymbol\lambda-\vec W\vec u(\vec b)\|^2_2.
	\end{equation}
 using the notation from \eqref{lq}, to be solved for $\vec\lambda$ in the least-squares sense. The effect of the preconditioning is that matrix rows corresponding to the two spherical caps of collocation points will be weighted, with their number density taken into account. 
 \end{remark}

\begin{figure}[h!]
				\centering
  		\begin{subfigure}[t]{0.33\textwidth}
			\centering
   %\hspace*{-2ex}
			\includegraphics[trim={5.2cm 17.3cm 9.5cm 4.3cm},clip,width=\textwidth]{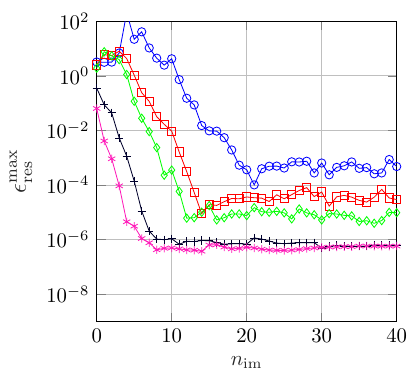}
			\caption{Squeezing motion.}		
      \label{squeeze_res}
		\end{subfigure}~
    		\begin{subfigure}[t]{0.33\textwidth}
			\centering
   %\hspace*{-1ex}
			\includegraphics[trim={5.2cm 17.3cm 9.5cm 4.3cm},clip,width=\textwidth]{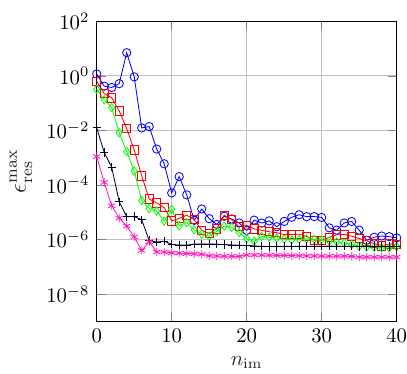}
			\caption{Shearing motion. }	
   \label{shear_res}
		\end{subfigure}~
  \begin{subfigure}[t]{0.33\textwidth}
			\centering
  % \hspace*{-2ex}
			\includegraphics[trim={0.79cm 17.3cm 13.5cm 4.3cm},clip,width=1.045\textwidth]{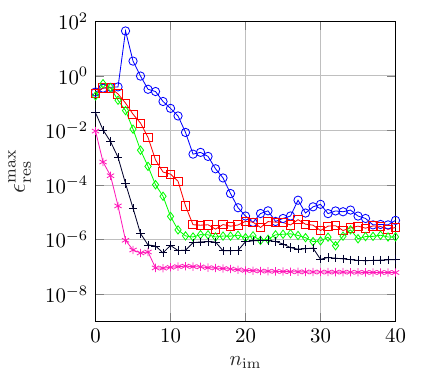}

			\caption{Equal and opposite rotation.}		
   \label{rot_res}
		\end{subfigure} \\
      		\begin{subfigure}[t]{0.33\textwidth}
			\centering
  % \hspace*{-2ex}
			\includegraphics[trim={5.2cm 17.3cm 9.5cm 4.3cm},clip,width=\textwidth]{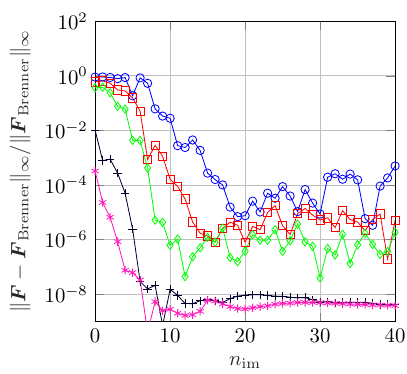}
			\caption{Squeezing motion: force error relative to the analytic result of Brenner in \eqref{Brenner}.}		
   \label{squeeze_force}
		\end{subfigure}~
  \begin{subfigure}[t]{0.8\textwidth}
			\centering
    \vspace*{-31ex}
			\includegraphics[trim={2cm 12.5cm 0cm 10.6cm},clip,width=1.1\textwidth]{figures/Two_rotation.pdf}
			%			\caption{Flow field in the plane of the triangle.}
			%			\label{flow_triangle}
			%		\end{subfigure}~~
			%		\begin{subfigure}[t]{0.33\textwidth}
			%		\centering
			%		\includegraphics[trim={2.9cm 1cm 0cm 0cm},clip,width=\textwidth]{../figures/images/triangle_zoom.eps}
			%		\caption{Zoomed in version of the flow-field in (a).}
			%		\label{flow_triangle2}
			%\caption{$\delta = 10^{-2}$, max relative residual}		
		\end{subfigure}
  % \hspace*{28ex}
  %   		\begin{subfigure}[b]{0.24\textwidth}
		% 	\centering
  %  \hspace*{-2ex}
		% 	\includegraphics[trim={5cm 16cm 8cm 4.2cm},clip,width=1.05\textwidth]{figures/Two_parallel_magn.pdf}
		% 	\caption{Squeezing motion: coefficient magnitude.}		
  %  \label{Fig:magn1}
		% \end{subfigure}~
  %   		\begin{subfigure}[b]{0.24\textwidth}
		% 	\centering
  %  \hspace*{-2ex}
		% 	\includegraphics[trim={5cm 16cm 8cm 4.2cm},clip,width=1.05\textwidth]{figures/Two_perp_magn.pdf}
		% 	\caption{Shearing motion: coefficient magnitude.}		
  %  \label{Fig:magn2}
		% \end{subfigure}~~
  % \begin{subfigure}[b]{0.24\textwidth}
		% 	\centering
  %  \hspace*{-2ex}
		% 	\includegraphics[trim={5cm 16cm 8cm 4.2cm},clip,width=1.05\textwidth]{figures/Two_rot_magn.pdf}
		% 	\caption{Equal and opposite rotation: coefficient magnitude.}		
  %  \label{Fig:magn3}
		% \end{subfigure}~
		\caption{Two spheres separated by $\delta$ require a different number of image points $n_{\text{im}}$ per particle depending on the boundary conditions; squeezing is the hardest to resolve. In all tests, $n_\text{im} = 20$ is enough to reach an error level of $10^{-3}$ or less, displayed for the max relative residual $\epsilon_{\text{res}}^{\max}$ in panels (a), (b) and (c), and the force error for squeezing in panel (d).}
	\end{figure}

\begin{figure}[h!]
 \centering
 \hspace*{-2ex}
 \includegraphics[trim = {5cm 16cm 6.2cm 4.5cm},clip,width=0.5\textwidth]{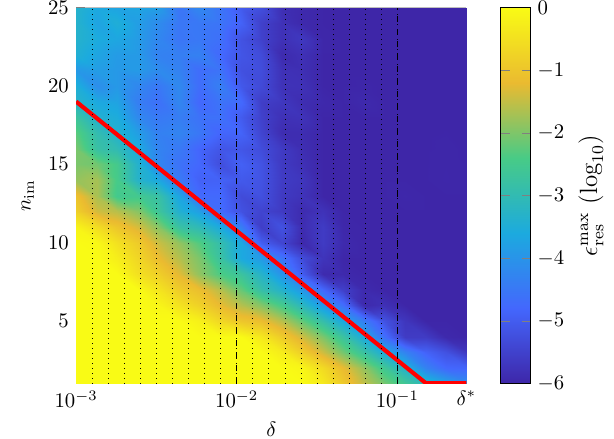}
 \caption{The number of image points per particle, $n_{\text{im}}$, varied for two spheres a varying distance $\delta$ apart undergoing squeezing motion, discretized with $N=686$ and $R_p = 0.63$. The red line is the fitted estimate \eqref{image_est} of the number of image points required per particle per close interaction, as a function of $\delta$, to reach the target accuracy $10^{-3}$. }
 \label{sweep_im}
 \end{figure}

 We now show that the above parameter choices achieve the target residual error for various boundary
 conditions.
 %The investigation will also lead to a strategy for choosing $n_{\text{im}}$, given $\delta$.
 % First, we determine the maximum relative residual from problems where the boundary conditions are generated from 40 different randomly sampled velocity vectors; see results visualized in Figure \ref{vary_bc}. The number of image points needed to reach a target accuracy of $10^{-3}$ is \note{Give interpretation.} In the same figure, the maximum relative residual corresponding to two static particles in a shear background flow is displayed. In contrast to with randomly sampled rigid body motions, no image points are then needed to resolve the interaction with sufficient accuracy. This result is in line with that presented in the previous section. Secondly, 
 The convergence with $n_\text{im}$ is measured for three distinct settings with lubrication:
	\begin{enumerate}
		\item Equal and opposite translational velocities, with the  spheres travelling towards each other (squeezing motion). Similarly to the test with squeezing motion without images in Figure \ref{Rg_sweep} in Section \ref{param_sec}, the computed forces on the particles are compared to the analytical expansion by Brenner in \eqref{Brenner}. The maximum relative residuals $\epsilon_{\text{res}}^{\max}$ over the particle surfaces are shown in Figure  \ref{squeeze_res}, together with the convergence of the relative force errors in Figure \ref{squeeze_force}.
		\item Equal and opposite translational velocities, with the spheres moving perpendicular to their line of centers (shearing motion), with residual convergence in Figure \ref{shear_res}.
		\item Equal and opposite angular velocities. The spheres rotate around the $z$-axis and are placed next to each other on the $x$-axis, with residual convergence in Figure \ref{rot_res}.
		\end{enumerate}
As convergence occurs in these tests, the coefficient norm drops from roughly $10^8$ down to $10^3$, possibly interpreted as an insufficient $n_{\text{im}}$ struggling to capture the true singularity line density.
The squeezing motion is the most difficult boundary condition to resolve, requiring the largest number of image points to reach the target accuracy of $10^{-3}$.
We investigate this further in Figure \ref{sweep_im}, where a sweep over $n_{\text{im}}$ is made for two spheres undergoing squeezing motion, a variable distance $\delta$ apart. Based on this, we may estimate the number of image points needed to reach a target accuracy of $10^{-3}$, given $N = 686$ and $R_p = 0.63$, as 
    \begin{equation}\label{image_est}
        n_{\text{im}}(\delta) = \left\{\begin{aligned} \ceil{-8.72\log_{10}(\delta) -6.15}, \quad &\delta < 0.15, \\ 
        0, \quad&\text{otherwise.} \end{aligned}\right.
    \end{equation}
   This is our adaptive rule for $n_{\text{im}}$, and it will be validated in Section \ref{Num_res}.  For all $\delta>10^{-3}$, note specifically that 20 image points per particle per contact are enough to resolve lubrication to the desired level of accuracy, corresponding to only 60 additional unknown vector source strengths.
   %Given the estimate in \eqref{image_est},

Finally, we may apply the above pair-wise recipe to the multi-particle case using the general singularity notation of Section \ref{image_enhance}.
Let the $k$th particle have $n_c \in \{0,1,\dots,\}$ near contacts, with separations $\delta_j$, $j=1,\dots,n_c$.
Its number of Stokeslets is then $N_\text{S}^{(k)} = N + \sum_{j}^{n_c}n_{\text{im}}(\delta_j)$, where the first term accounts for the coarse proxy points.
We also have
$N_\text{R}^{(k)} = N_\text{D}^{(k)} = \sum_{j}^{n_c}n_{\text{im}}(\delta_j)$,
while $N_\text{T}^{(k)} = 0$,
giving the total number of vector unknowns $N^{(k)}$ as in \eqref{Nk}.
The associated total number of collocation nodes for particle $k$ is
   \begin{equation}\label{Mest}
       M^{(k)} = 1.2\cdot N +2\cdot 6\cdot\sum_{j=1}^{n_c} 3n_{\text{im}}(\delta_j),
   \end{equation}
	the first term being the coarse nodes and the rest (recalling $\alpha_1=\alpha_2=6$) the two added node caps.
% the first term represents a uniform discretization of points over the surface of the sphere, as in our basic MFS setup, and the second term represents the sum of the two spherical caps of distance-dependent number of additional points per near-contact, surrounding the contact points.

 %We reiterate that each such spherical cap is discretized with a multiple of the number of sources, $3n_{\text{im}}$.
	% \begin{figure}[h!]
	% 	\centering
	% 	\begin{subfigure}[b]{0.49\textwidth}
	% 		%\includegraphics[trim = {12cm 0cm 0cm 0.5cm},clip,width=1\textwidth]{../figures/images/spherical_cone.eps}
	% 		\caption{The number of image points needed to resolve the interaction resulting from randomly sampled rigid body motions for a pair of spheres is reported. The result is compared to resolving the interaction for a static pair of spheres in a background shear flow. }
	% 		\label{vary_bc}
	% 	\end{subfigure}~~
	% 	\begin{subfigure}[b]{0.49\textwidth}
	% 		%\includegraphics[trim = {12cm 0cm 0cm 0.5cm},clip,width=1\textwidth]{../figures/images/spherical_cone.eps}
	% 		\caption{Three distinct boundary conditions are investigated for different $\delta$. The most difficult to resolve is \noteme{fill in description}.}
	% 		\label{force_err}
	% 	\end{subfigure}
	% 	\caption{Two spheres separated by $\delta$. The number of image points needed to resolve the interaction depends on the boundary condition. }
	% \end{figure}

	%%%%%%%%%%%%%%%%%%%%%%%%%%%%%%%%%%%%%%%%%%%%%%
	\subsection{Validating the pairwise-adaptive strategy in clusters of spheres}\label{Num_res}
 We now present four examples with clusters of spheres tightly packed together, where lubrication is resolved by choosing the number of image and collocation points adaptively, following the pairwise image line discretization \eqref{image_est}.
 We will find that simple unions of these pairwise image lines perform just as accurately for multi-sphere clusters, a remarkably useful result (whose explanation may be related to
 Cheng--Greengard's three-disk analysis \cite[Lem.~3.3]{Cheng1998}).
 
 Recall that the coarse proxy sphere always has $N=686$ and $R_p = 0.63$. In self-convergence tests, we compare against results determined with an image discretization tuned for a larger number of sources, where $N = 1353$, $R_p =0.7$, and 30 image points per particle per near-contact are used. In one of the examples, a BIE reference solution is used to validate the self-convergence error, but also to highlight the difficulties with BIE that we can overcome with the proposed MFS technique. This BIE solution is computed with a double layer formulation equipped with quadrature by expansion (QBX), developed by af Klinteberg \& Tornberg \cite{AfKlinteberg2016} and Bagge \& Tornberg \cite{Bagge2021}, where the quadrature tolerance is to $10^{-6}$.
 %%%%%%%%%%%%%%%%%%%%%%%%%%%%%%%%%%%%%%%%%%%
	\begin{example}[A self-convergence test for the exterior flow field] In the first example, a triangle of spheres a distance $\delta = 10^{-3}$ apart move with translational velocities of unit magnitude towards the center of the triangle. The flow field is evaluated in a large set of points in the plane cutting through the spheres, and in a self-convergence test, the maximum relative velocity difference in the exterior domain is $7.6\cdot 10^{-4}$, i.e., below the target accuracy of $10^{-3}$. The flow field is displayed in Figure \ref{flow_triangle}, while the relative difference in the flow field is displayed in Figure \ref{flow_triangle_diff}. 
	\begin{figure}[h!]
		\centering
		\begin{subfigure}[t]{0.58\textwidth}
			%			\centering
			%			\includegraphics[trim={2.9cm 1cm 4cm 0cm},clip,width=\textwidth]{../figures/images/triangle.eps}
			%			\caption{Flow field in the plane of the triangle.}
			%			\label{flow_triangle}
			%		\end{subfigure}~~
			%		\begin{subfigure}[t]{0.33\textwidth}
			%		\centering
			%		\includegraphics[trim={2.9cm 1cm 0cm 0cm},clip,width=\textwidth]{../figures/images/triangle_zoom.eps}
			%		\caption{Zoomed in version of the flow-field in (a).}
			%		\label{flow_triangle2}
			\begin{minipage}[b]{.42\textwidth}
				\centering
				\includegraphics[trim={3.15cm 1.7cm 4cm 1.7cm},clip,width=\textwidth]{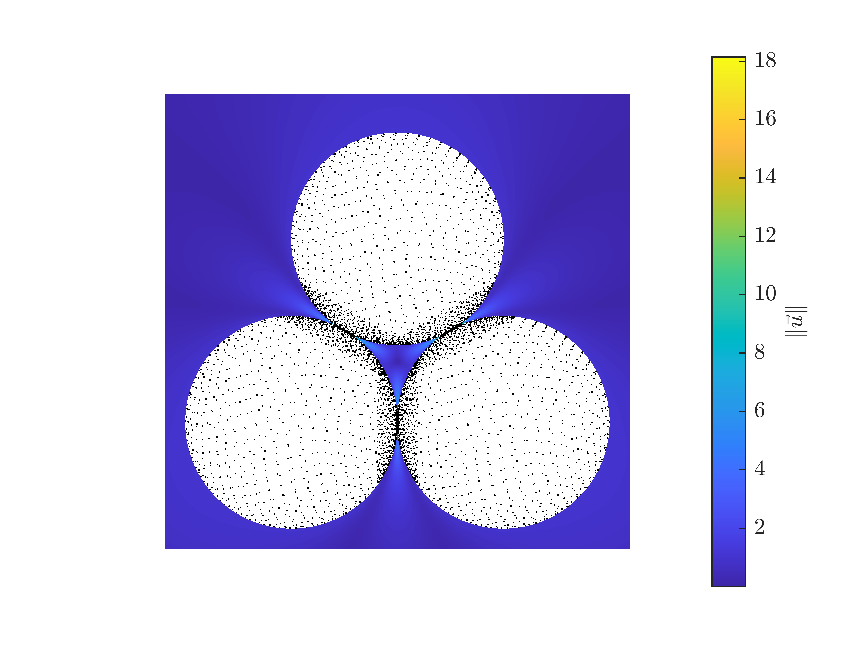}
				%	\caption{$dt=0.1$}
				%	\label{fig:prob1_6_2}
			\end{minipage}~
			\begin{minipage}[b]{0.6\textwidth}
				\centering
    \begin{overpic}[trim={2.9cm 1cm 0cm 0.8cm},clip,width=1.05\textwidth]{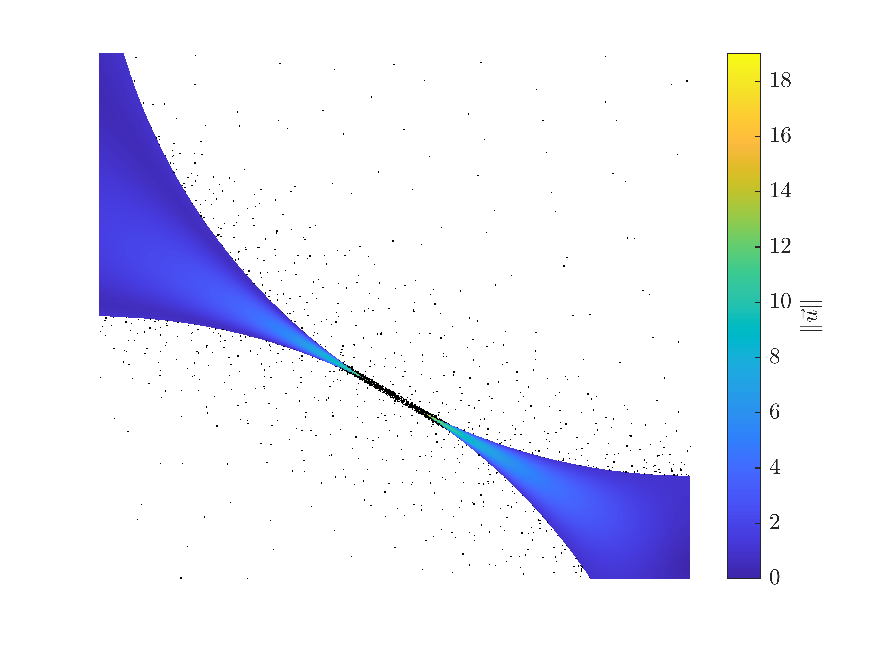}
\put(88,33){\rotatebox{90}{\colorbox{white}{\parbox{0.1\linewidth}{%
     \small$\|\vec u(\vec x)\|_2$}}}}
     \end{overpic}
				%	\includegraphics[width=0.3\linewidth, height=0.15\textheight]{prob1_6_1}
				%\includegraphics[trim={2.9cm 1cm 0cm 0.8cm},clip,width=1.1\textwidth]{figures/triangle_zoom.eps}
				%	\caption{$dt =$}
				%	\label{fig:prob1_6_1}
			\end{minipage}
			\caption{%Left: Flow field in the plane of the triangle. Right: zoomed in plot between two of the particles -- the flow field is the largest in the particle gap.
   }
			\label{flow_triangle}
		\end{subfigure}~~~~~
		\begin{subfigure}[t]{0.38\textwidth}
			\centering
			\vspace*{-31.5ex}
       \begin{overpic}[trim={2cm 1.2cm 0cm 0.6cm},clip,width=1.05\textwidth]{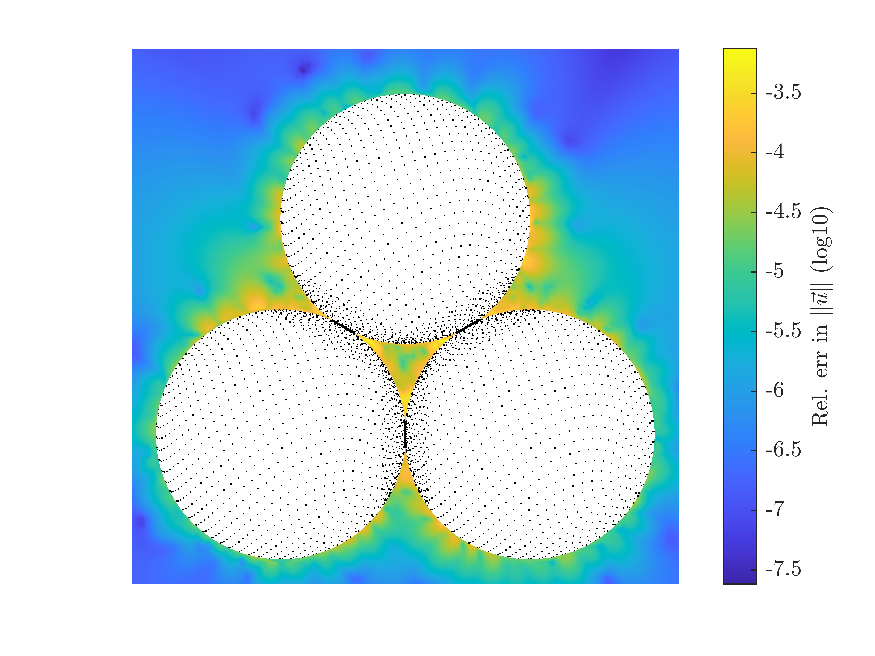}
\put(91,0){\rotatebox{90}{\colorbox{white}{\parbox{0.8\linewidth}{%
     \small$\log_{10}\left(\|\vec u-\vec u_{\text{ref}}\|_2/\max\limits_{\vec x}\|\vec u_{\text{ref}}\|_2\right)$}}}}
     \end{overpic}
			\caption{}
			\label{flow_triangle_diff}
		\end{subfigure}
		\caption{Self-convergence test for an example with three unit spheres forming a triangle of side length $2+10^{-3}$. No-slip boundary conditions are applied, with each sphere assigned a translational velocity pointing towards the center of the triangle, resulting in lubrication forces between the particles. In panel (a), the flow field is visualized in the plane of the triangle, with a zoomed-in plot between two of the particles to the right, where the flow velocity magnitude is the largest. In panel (b), the pointwise relative difference in the computed flow field using a coarser and a finer MFS grid is displayed. With the 20 image points per particle per near-contact used in the example, the maximum relative difference is $7.6\cdot 10^{-4}$.}
		%A 200 \times 200 grid is used to compute the flow field. 
		\label{image_figs}
		
	\end{figure}
	\end{example}
 %%%%%%%%%%%%%%%%%%%%%%%%%%%%%%%%%%%%%%
 \begin{example}[Accuracy gain with adaptive images]\label{ex_images}
 We display the difference in error level with the adaptively set image points compared to the basic MFS with proxy-sources only.  The geometry is a tetrahedron of unit spheres, as displayed in Figure \ref{tetra_conv}, with varying side length of the tetrahedron, given by $2+\delta$, and the  relative force/torque error $\epsilon_{\text{FT}}$, as defined by \eqref{force_err}, is quantified for each $\delta$ against the fine MFS grid in Figure \ref{noimage}. Errors $\epsilon_{\text{FT}}$ are less than the target tolerance of $10^{-3}$ with the use of images, while without image points, accuracy suffers when $\delta<\delta^*$, as expected.  For each $\delta$, 20 different problems are solved, with randomly sampled no-slip boundary conditions. %in the former with randomly sampled translational and angular velocities for each $\delta$, and in the latter via randomly sampled forces and torques. 
 \begin{figure}[h!]
		\centering
  		\begin{subfigure}[t]{0.44\textwidth}
			\centering
			%\vspace{-31ex}
			\includegraphics[trim={20cm 4cm 16cm 6.5cm},clip,width=0.9\textwidth]{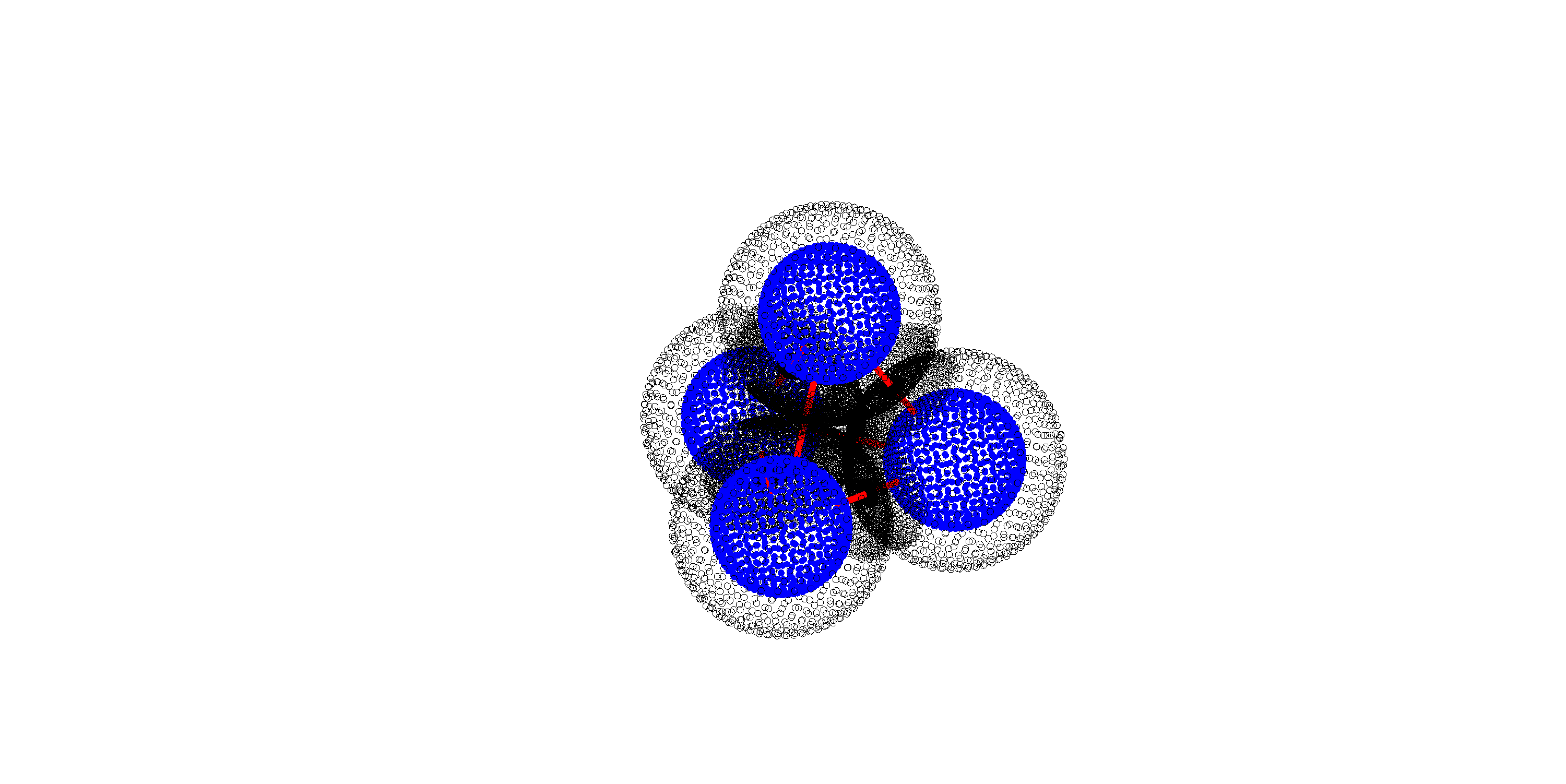}
			\caption{Test geometry in Examples \ref{ex_images} and \ref{ex_BIE}.}
			\label{tetra_conv}
		\end{subfigure}~~~
		\begin{subfigure}[t]{0.55\textwidth}
			\centering
	%\vspace{-31ex}
	\includegraphics[trim={0cm 0.1cm 1cm 0.8cm},clip,width=\textwidth]{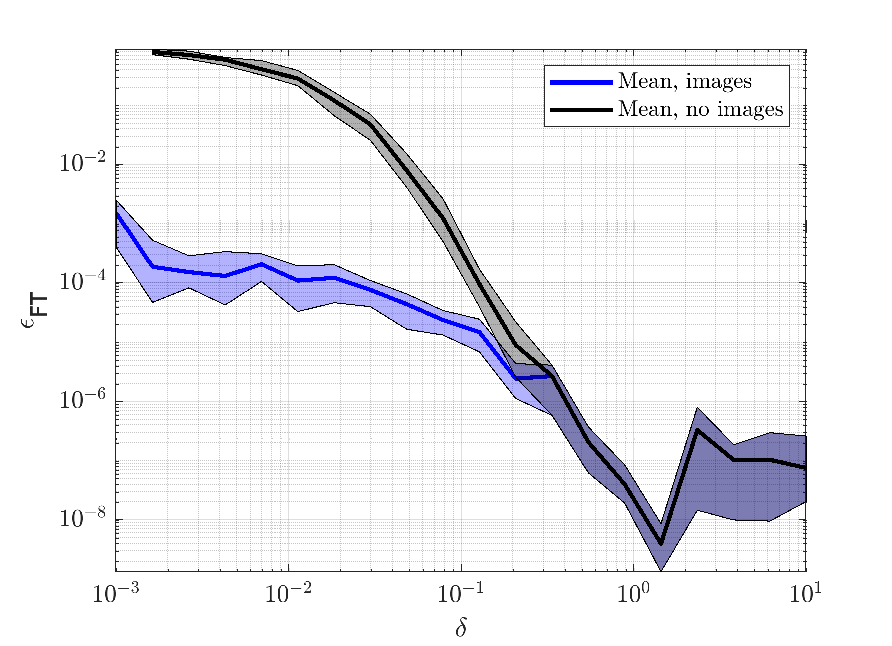}
	\caption{Self-convergence test for the geometry in (a).}
	\label{noimage}
		\end{subfigure}%~~
  \caption{In panel (a), four unit spheres form a tetrahedron of side length $2+\delta$. Source points on proxy-surfaces are visualized in blue, image points in red and collocation points in black. The picture corresponds to the coarse MFS grid, with $N = 686$ and $n_{\text{im}}=20$ image points per particle per near-contact. Extra collocation points are added close to the points of near-contact, following \eqref{Mest}. Panel (b) displays a self-convergence test with varying $\delta$. In blue, accuracy is controlled with image and collocation points set based on the estimate \eqref{image_est}. In black, without extra source and collocation points, the relative error in the computed force and torque, $\epsilon_{\text{FT}}$ defined by \eqref{force_err}, diverges with decreasing $\delta$.  For $\delta>\delta^*$, no image points are used in either case. Shaded regions display max and min for 20 different randomized boundary conditions.}
		%\caption{Convergence test for four unit spheres forming a tetrahedron of side length $2+\delta$. The reference computed with a boundary integral formulation is inaccurate for small $\delta$ due to the under-resolved peak in the layer density that causes increasing errors for decreasing particle distances. The effect is visible both for the coarse and fine QBX grid in (a). For larger $\delta$, the error level in the self-convergence test is similar in magnitude to the error in the comparison to a solution computed with QBX. %\note{Note that these results are with the adaptively set $n_{\text{im}}$.}
  %}
  \end{figure}
 \end{example}
 %%%%%%%%%%%%%%%%%%%%%%%%%%%%%%%%%%%%%%
 \begin{example}[Comparison of accuracy and efficiency relative to a BIE scheme]\label{ex_BIE}
With test geometry as in Example \ref{ex_images}, the following example serves two purposes:   to validate the self-convergence error against results determined with an established BIE solver equipped with QBX quadrature, and to highlight the ability of the MFS scheme to resolve the interaction for problems where a BIE struggles to reach accuracy \footnote{Since the mobility problem---where forces and torques on the particles are provided then rigid body motions of the particles are computed---generally converges faster than the resistance problem under study, we choose to generate reference solutions with QBX via the mobility formulation; forces and torques are assigned and the computed velocities are used as input to our MFS scheme.}. For moderately small $\delta$, where the peaks in the BIE layer density are not as severe, comparisons are possible, while for very small $\delta$, the force density peak of width $\mathcal O(\sqrt{\delta})$ \cite{Sangani1994} requires an impractically fine grid to be resolved in the BIE setting. In addition, the precomputations needed for employing QBX are more costly the finer the resolution. We test two grid-resolutions for use with QBX to demonstrate these claims: a coarser grid with $30\times 30$ discretization points per particle, and a finer grid with $60\times 90$ grid points
(5400 vector unknowns).
Due to CPU time constraints, the fine QBX-based reference is only determined for two sets of boundary conditions at each $\delta$: %\note{with an approximate cost of one hour per solve on a workstation, but don't know if I should write that out...}: 
a unit torque applied in the $z$-direction on every particle, or a unit force in the $x$-direction on every particle. In Figure \ref{convergence_tetra}, the error relative to the fine-grid BIE remains small for smaller values of $\delta$ compared to the coarser BIE-grid, yet both errors start diverging for sufficiently small $\delta$. This onset varies with boundary condition, but above this, comparisons against both the coarse and fine BIE-grids show similar errors as the self-convergence test (which displays the same data as in Example \ref{ex_images}), validating the MFS against the BIE solution. For large $\delta>0.3$, the BIE errors match the QBX tolerance $10^{-6}$. In the self-convergence test and for comparisons to the coarse BIE solution, 20 different problems are solved for each $\delta$, with randomly sampled no-slip boundary conditions. %in the former with randomly sampled translational and angular velocities for each $\delta$, and in the latter via randomly sampled forces and torques.
 %It can be concluded that, independently of $\delta$, $\epsilon_{\text{FT}}$ on the particles is less than the target tolerance of $10^{-3}$. 
	\begin{figure}[h!]
		\centering
		%\begin{subfigure}[t]{0.55\textwidth}
			\centering
			\includegraphics[trim={0.2cm 0cm 1cm 0.65cm},clip,width=0.55\textwidth]{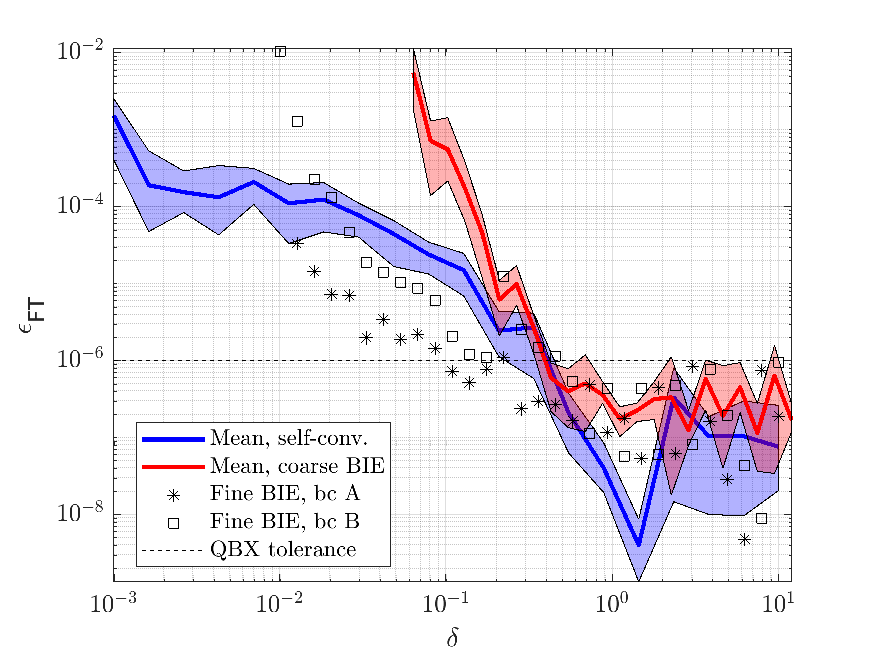}
			%			\caption{Flow field in the plane of the triangle.}
			%			\label{flow_triangle}
			%		\end{subfigure}~~
			%		\begin{subfigure}[t]{0.33\textwidth}
			%		\centering
			%		\includegraphics[trim={2.9cm 1cm 0cm 0cm},clip,width=\textwidth]{../figures/images/triangle_zoom.eps}
			%		\caption{Zoomed in version of the flow-field in (a).}
			%		\label{flow_triangle2}
			%\caption{}
			
		%\end{subfigure}~~
		% \begin{subfigure}[t]{0.4\textwidth}
		% 	\centering
		% 	%\vspace{-31ex}
		% 	\includegraphics[trim={20cm 4cm 16cm 6.5cm},clip,width=0.9\textwidth]{figures/tetra_geom.eps}
		% 	\caption{Source points on proxy-surfaces are visualized in blue, approximate image points in red and collocation points in black. The picture corresponds to the coarse grid evaluated in (a), with $N=686$ and 20 image points per near-contact. Extra collocation points are added close to the points of near-contact.}
		% 	\label{tetra_conv}
		% \end{subfigure}	
		\caption{Comparison to results from a BIE formulation for the tetrahedron test geometry of Figure \ref{tetra_conv};  forces and torques determined with the image enhanced grid  are compared to different sets of reference data. For small $\delta$, the BIE-reference is inaccurate  due to the under-resolved peak in the layer density. This  causes increasing errors for decreasing $\delta$ in comparisons of coarse MFS results both to solutions determined with coarse and fine BIE-grids (see details in main text), but with different onsets of the increase. For larger $\delta$, however, the error level  in the forces and torques from the self-convergence test is similar in magnitude to the error in the comparison to a solution computed with BIE. \revFour{For BIE, the QBX-tolerance is set to $10^{-6}$, which together with the GMRES tolerance restrict how small errors that can be detected for larger $\delta$.}    %\note{Note that these results are with the adaptively set $n_{\text{im}}$.}
  Red and blue represent the same type of comparisons, but using different reference solutions: for each separation distance $\delta$, the experiment is repeated 20 times with randomly generated boundary data, and the min, max and mean of $\epsilon_{\text{FT}}$ are visualized. Test problems generated
with a unit force or torque are labeled ’bc A’ and ’bc B’. The self-convergence data is the same as presented in Figure \ref{noimage}.}
		%A 200 \times 200 grid is used to compute the flow field. 
		\label{convergence_tetra}	
	\end{figure}
	
\end{example}
	%\clearpage
 %%%%%%%%%%%%%%%%%%%%%%%%%%%%%%%%%%%%%%%%%%%%%%%%

 %\subsection{A larger example}
 \begin{example}[Adaptive discretizations for randomized clusters]\label{adaptive}
	We study random clusters of the type in Figure \ref{image_example} consisting of 20 particles, with each sphere $\delta = 10^{-3}$ away from at least one neighbor.  In such a cluster, some spheres will have many close interactions, meaning that the adaptive number of image and collocation points will vary considerably between different particles.  The spheres have randomly sampled rigid body velocities, with the resulting surface velocity exemplified in Figure \ref{surf_vel} and relative residual over the sphere surfaces in Figure \ref{cluster_resim}. Statistics for the maximum relative residuals, $\epsilon_{\text{res}}^{\max}$, for 50 realizations of the experiment, are shown in Figure \ref{cluster_data}. The same figure reports statistics for the %number of GMRES iterations and the 
 number of close contacts for each realization. The adaptively determined number of image and collocation points for each sphere in each run are displayed in Figures \ref{image_adaptive} and \ref{collocation_adaptive}, where the short bars for the largest number of points correspond to the few particles that have the largest number of near contacts in the center of the cluster. Note, however, that the maximum number of collocation points never exceeds the 5400 vector unknowns in the fine-grid BIE in Example \ref{ex_BIE} (being non-adaptive), and the solution with MFS is also much faster to compute. A majority of the particles have only one close neighbor, resulting in moderate numbers of unknowns, which is reflected in the tall bars in Figures \ref{image_adaptive} and \ref{collocation_adaptive}.
 %A similar test case with 20\% packing fraction was studied in \cite{Yan2020}.
% \clearpage
	\begin{figure}[h!]
		\centering
		% \begin{subfigure}[t]{0.4\textwidth}
		% 	\centering
		% 	\includegraphics[trim={0.5cm 1cm 0cm 0cm},clip,width=1.1\textwidth]{figures/image_cluster.eps}
  %           \caption{Example configuration displaying the surface velocities. }		
		% 	\label{surf_vel}
		% \end{subfigure}~~
  % 		\begin{subfigure}[t]{0.4\textwidth}
		% 	\centering
		% 	\includegraphics[trim={0.5cm 1cm 0cm 0cm},clip,width=1.1\textwidth]{figures/residual_cluster.eps}
		% 	\caption{Example configuration displaying the relative residual at the particle surfaces. }
		% 	\label{cluster_resim}
		% \end{subfigure}\\
    		\begin{subfigure}[t]{0.37\textwidth}
			\centering
   \begin{minipage}[t]{0.46\textwidth}
   	\includegraphics[trim={5cm 16.5cm 12.7cm 4cm},clip,width=\textwidth]{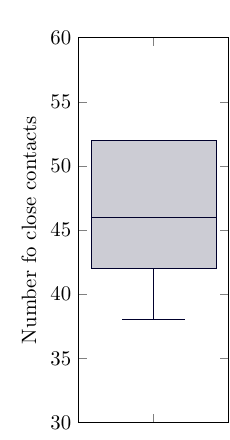}		
   \end{minipage}
      \begin{minipage}[t]{0.51\textwidth}
   	\includegraphics[trim={5cm 16.5cm 12.2cm 4cm},clip,width=1.02\textwidth]{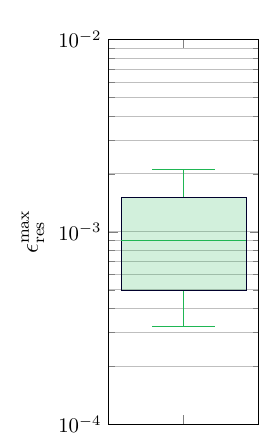}		
   \end{minipage}
   \caption{Statistics for the total number of close contacts (left) and the max relative residual (right).}
			\label{cluster_data}
		\end{subfigure}~~~
      		\begin{subfigure}[t]{0.3\textwidth}
			\centering
   \hspace*{-6ex}
			\includegraphics[trim={4.5cm 16cm 8.5cm 4cm},clip,width=1.1\textwidth]{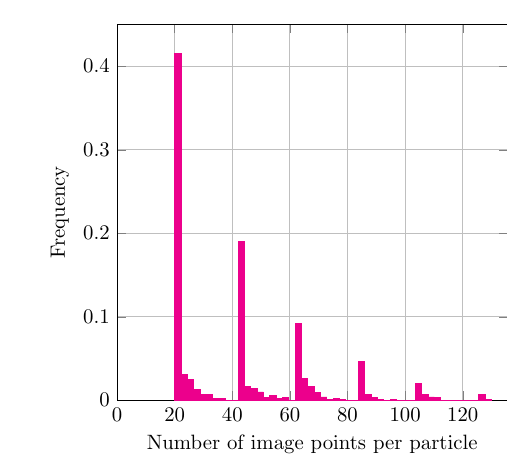}
			\caption{Statistics for the number of image points per sphere. }
\label{image_adaptive}
		\end{subfigure}~~
      		\begin{subfigure}[t]{0.3\textwidth}
			\centering
        \hspace*{-6ex}
			\includegraphics[trim={4.5cm 16cm 8.5cm 4cm},clip,width=1.1\textwidth]{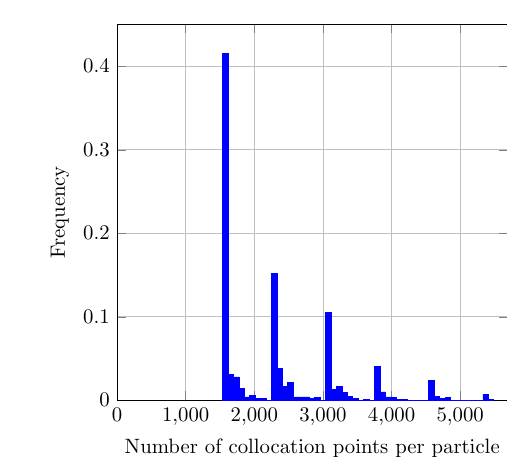}
			\caption{Statistics for the number of collocation points per sphere. }
			\label{collocation_adaptive}
		\end{subfigure}
  \caption{By using an adaptive number of image points per closely interacting neighbor, as in \eqref{image_est}, lubrication forces between $P=20$ spheres in a cluster with minimum separation $\delta = 10^{-3}$ are resolved.  Each sphere has a randomly sampled translational and angular velocity. The experiment is repeated 50 times. In panel (a), whiskers display min/max, box edges display the 10th and 90th percentile of results from the 50 runs, and the box center line displays the median.}
  \label{adaptive_test}
  \end{figure}
  \end{example}
	%	\clearpage	

 %%%%%%%%%%%%%%%%%%%%%%%%%%%%%%%%%%%%%%%%%%%%
 \subsection{GMRES iteration count}\label{Gmres_count}
We have demonstrated controllable accuracy with our scheme on challenging problems with strong lubrication forces. In terms of efficiency, the computational grid can be on average much coarser than what is needed in a uniformly-discretized BIE method to resolve the interactions. It however remains to discuss the cost of iteratively solving the preconditioned system for small $\delta$. Here, we show for two of the example geometries from the previous section that the
  GMRES iteration count grows with decreasing particle distances, and that this phenomenon is common to both MFS and BIE methods.
  %We will however demonstrate that this trend is not isolated to hold for our MFS technique, but holds also in a BIE-setting. 

 \begin{example}[continues=ex_BIE]%, Comparison relative to BIE scheme: GMRES iterations]
 For the tetrahedron of spheres pair-wise separated by $\delta$, the GMRES iteration count as a function of $\delta$ is reported in Figure \ref{iters}. The iteration count appears to grow like ${\mathcal O}(\delta^{-1/2})$ for all the tested boundary conditions (random velocity data, or, based on the reference solution, with either random forces and torques, or with a unit force/torque on each particle). The iteration counts required to provide the BIE reference with the fine grid are similar to those needed to compute the MFS solution with the same boundary condition. Note, however, that the numbers of BIE unknowns are larger -- $5400\cdot 3$ per particle vs.~$3M^{(i)}=3\cdot\left(801 + 12\cdot n_c\cdot3\cdot 20\right) = 3\cdot2960$ %2\times 3\times 6\times 20 = 2930$ 
	with MFS, following \eqref{Mest}. The displayed BIE GMRES iterations are for the inverse resistance problem (mobility), which we have found converges faster. %The GMRES tolerance is also stricter, $10^{-8}$ vs $10^{-6}$. 

	\begin{figure}[h!]
 \centering
  \begin{subfigure}[t]{0.7\textwidth}
		\centering
		%\vspace{-31ex}
		\includegraphics[trim={1.2cm 0.2cm 1.6cm 0.8cm},clip,width=\textwidth]{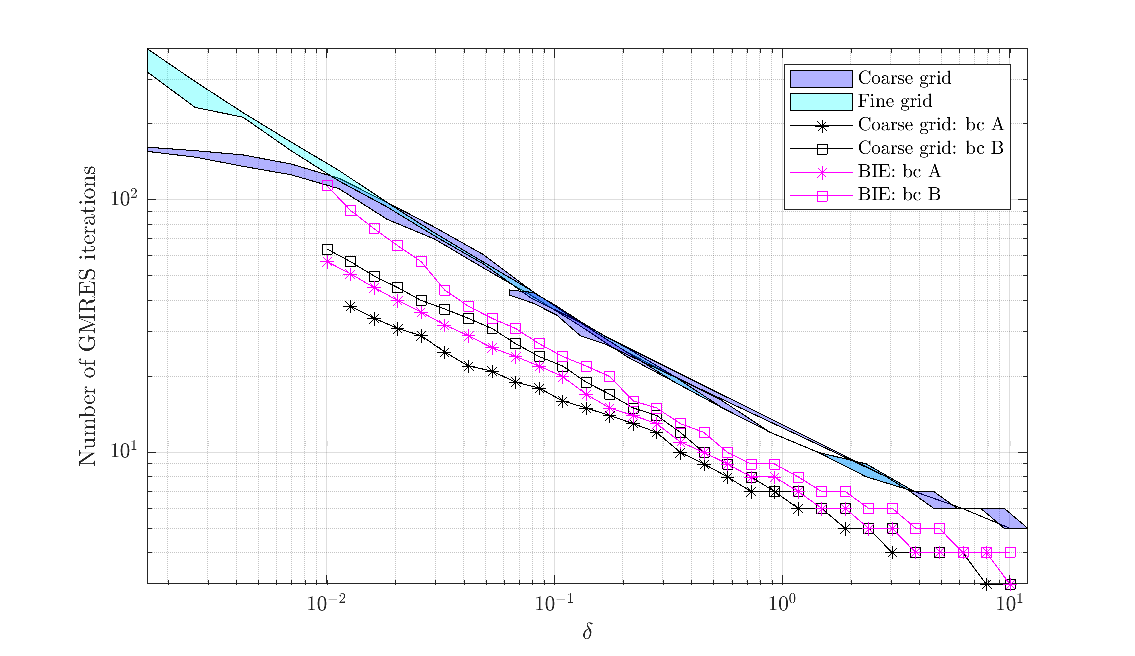}
		\caption{GMRES iteration count for the data in Figure \ref{convergence_tetra}. 
  }
		\label{iters}
   \end{subfigure}~~
   \begin{subfigure}[t]{0.25\textwidth}
   \vspace*{-44.2ex}
   	\includegraphics[trim={5cm 15.3cm 12.3cm 4cm},clip,width=0.8\textwidth]{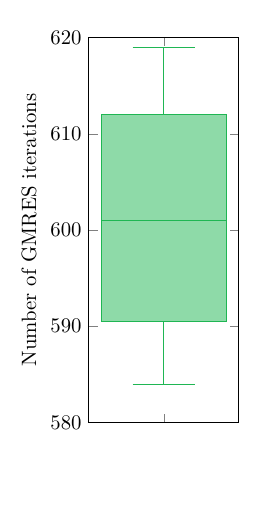}		
    \caption{GMRES iteration count for the data in Figure \ref{adaptive_test}.}
    \label{adaptive_GMRES}
   \end{subfigure}
   \caption{GMRES iteration count for the data in Examples \ref{ex_BIE} (panel (a)) and \ref{adaptive} (panel (b)), with the MFS  GMRES tolerance  set to $10^{-6}$. In panel (a) for the tetrahedron geometry, convergence is reported for different boundary data and compared to the convergence achieved with the BIE-formulation with QBX.  The same trend of an increasing iteration count for decreasing $\delta$ is noted for both methods. The shaded data  refers to max and min iteration counts corresponding to boundary conditions with  randomly sampled rigid body velocities, 20 different sets for %thfor
  the fine grid %'bc set 1' and velocities obtained from a mobility solve with BIE and 20 random forces and torques for 'bc set 2'. 
  and 40 sets for the coarse grid. Test problems generated with a unit force or torque are labeled 'bc A' and 'bc B' and converge slightly faster than problems with a boundary condition based on a randomly sampled velocity. % or force / torque vector.}
  Panel (b) displays iteration data for the clusters of 20 spheres with $\delta=10^{-3}$ and adaptively set image points, where whiskers display min/max, box edges the 10th and 90th percentile, and the box center lines display the median.}
	\end{figure}
 \end{example}

 \begin{example}[continues=adaptive] Statistics for the GMRES count required for the randomized clusters of 20 spheres with adaptively set image points is reported in Figure \ref{adaptive_GMRES}. 
 \end{example}
  To conclude, the GMRES count (for either MFS or BIE systems) becomes quite large for small $\delta$, and is consistent with ${\mathcal O}(\delta^{-1/2})$. Suitable preconditioning strategies are left to future work.

	\section{Conclusions and outlook}\label{conclusions}
	We consider exterior Stokes boundary value problems with no-slip boundary conditions applied to a collection of identical spherical particles immersed in a viscous fluid. The work presents a technique based on the method of fundamental solutions (MFS), with the flow in the basic setup represented as a linear combination of Stokeslets placed at proxy-surfaces in every particle interior.
  For problems with strong lubrication forces, this is augmented by a combination of Stokeslets, rotlets and potential dipoles that discretize each pair-wise ``image line'', providing controllable residuals for all separations down to $\delta = 10^{-3}$ for spheres of unit radii.
  We present a recipe, as a function of such a $\delta$, for the number of image points and accompanying ``caps'' of collocation nodes, to reach a target accuracy of $10^{-3}$ uniformly in the surface velocities.
  We show that, remarkably, {\em such line-image approximations developed for pairwise interactions also perform equally well for packed clusters with many close spheres}.
 
 Interactions for close-to-touching particles moving relative to each other are challenging to compute with established methods, such as a boundary integral (BIE) method with uniform Nystr\"om quadrature, because a sharply peaked force density has to be resolved. With the MFS technique, we demonstrate numerically that lubrication forces can be accurately resolved using much coarser particle representations (less unknowns) than a uniform BIE method at the same accuracy. It would be interesting to compare against spatially-adaptive BIE quadratures, which are less established. In the present work, by enforcing boundary conditions by collocation at the particle surfaces in the least-squares sense, the technical singular quadratures needed for BIE are avoided. In the special case without lubrication forces due to relative particle motions, with a backward-stable one-body preconditioning strategy, plain proxy sources are sufficient to reach the target accuracy. The sum of Stokeslets can then be rapidly evaluated using the fast multipole method, yielding close to linear scaling in the number of particles, which we demonstrate for systems with up to 4.8 million degrees of freedom.
	
	%\subsection{Outlook}
	A natural extension to the work reported here is to the mobility problem.
 %with the MFS, i.e.~the inverse to the resistance problem, where forces and torques are provided and rigid body velocities are to be computed.
 \ownchange{We report on a well-conditioned MFS formulation for the mobility problem separately in \cite{Broms2024c}.}\revThree{ In that work, we also extend to non-spherical smooth particles in the setting where proxy surfaces suffice to resolve interactions.
  %For simple particles with separations where proxy surfaces suffice to resolve interactions, such a generalization is straightforward.
 % with sources moved from the particle surface  in the normal direction into the particle interior. Proxy-surfaces of this kind were also discussed by Stein \& Barnett in \cite{Stein2022}. Similarly to the sphere case, self-interaction blocks could then be stored for each particle type for the one-body preconditioning strategy to be applied. 
  It is a topic of future research whether image representations could be useful for non-spherical particles.}
  %A naive approach of locally approximating a non-spherical particle by a sphere, image points might fall outside the particle, thus lie in the fluid domain. 
	In Section \ref{Gmres_count}, we note that the number of GMRES iterations grows with decreasing particle gaps. Future work includes investigating suitable preconditioning \revFour{or the use of a fast direct solver.}
 %Here, the challenge is to devise a preconditioning strategy of limited computational cost, that can substantially decrease the number of GMRES iterations even for spheres having different numbers of close interactions at varying separation distances.   
 As for the size of the systems considered in this work, all tests fit on a standard workstation. Larger systems could be studied via a distributed memory implementation and FMM, as in \cite{Yan2020}.
	
	\section*{Acknowledgements}
	Broms and Tornberg acknowledge the support from the Swedish Research Council: grant no.~2019-05206 and the research environment grant INTERFACE (biomaterials), no.~2016-06119. The Flatiron Institute is
a division of the Simons Foundation. Broms appreciates the research visit at CCM Flatiron, which served as a catalyst for the initial ideas in this work. \ownchange{We would also like to thank the referees for their thorough reading of the manuscript and for their valuable
comments which have helped us to improve the text. }

 \appendix
 \section{Singular values}\label{sec_trunc}
 As for the numerical examples in Section \ref{lub}, the coarse MFS grid with $N=686$ and $R_p=0.63$ is here used for numerical experiments. Further, for a sphere close-to-touching with one other sphere such that image points are used, $M = 1.2\cdot N+2\cdot 6\cdot 3 n_{\text{im}}$.
 
 \subsection{Truncation of singular values with image enhancement}
 It is difficult to get the approximation power of the image points as large as possible for a general choice of no-slip boundary conditions, while keeping the ill-conditioning under control. Truncation of singular values of the self-interaction blocks in the target-from-source matrix is hence necessary and the effect of the truncation level, $\epsilon_{\text{trunc}}$, relative to the largest singular value, is here investigated.

 As discussed in Section \ref{base_disc}, clustering of source points worsens the conditioning of the self-interaction blocks in the target-from-source matrix. Two clustering properties of the image points defined by \eqref{line_points}-\eqref{tj} impact this. As $\delta$ decreases, image points approach the particle surface ($R_{\text{acc}}(\delta)\to 1$), making the fundamental solutions near-singular when evaluated near the closest approach points, while for larger $\delta$, image points are further from the surface, reducing near-singularity issues. Additionally, with Chebyshev sampling and equal $n_{\text{im}}$, image points are more clustered for larger $\delta$, increasing ill-conditioning.

\begin{figure}[h!]
		\centering
    		\begin{subfigure}[b]{0.8\textwidth}
			\centering
  % \hspace*{-7.5ex}
   \hspace*{-10.5ex}
			\includegraphics[trim={0cm 16.3cm 0cm 10.7cm},clip,width=1.23\textwidth]{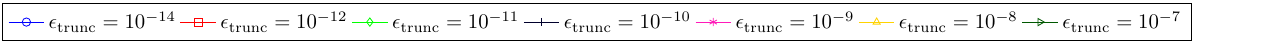}	
		\end{subfigure}\\
  		\begin{subfigure}[b]{0.24\textwidth}
			\centering
   \hspace*{-5ex}
			\includegraphics[trim={5cm 17.3cm 9.5cm 4.5cm},clip,width=1.15\textwidth]{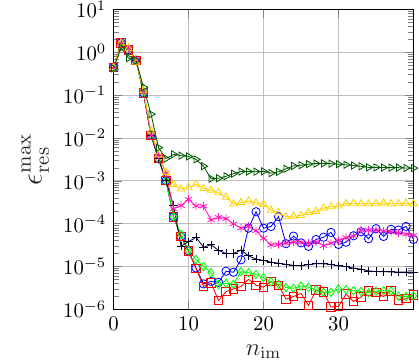}
			\caption{$\delta = 10^{-2}$, max relative residual $\epsilon_{\text{res}}^{\max}$.}			
		\end{subfigure}~~
    		\begin{subfigure}[b]{0.24\textwidth}
			\centering
   \hspace*{-5ex}
			\includegraphics[trim={0.5cm 17.3cm 14cm 4.5cm},clip,width=1.15\textwidth]{figures/Truncate_two_spheres_gap1e-2_magn.pdf}
			\caption{$\delta = 10^{-2}$, max coefficient magnitude $\|\vec\lambda\|_{\infty}$.}			
		\end{subfigure}~~
  \begin{subfigure}[b]{0.24\textwidth}
			\centering
   \hspace*{-3ex}
			\includegraphics[trim={5cm 17.3cm 9.5cm 4.5cm},clip,width=1.15\textwidth]{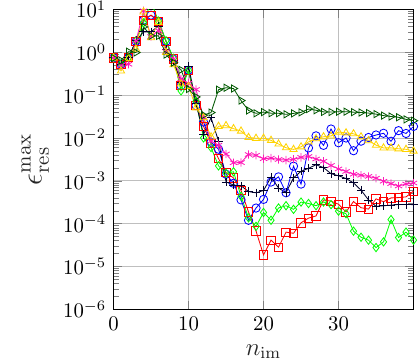}
			\caption{$\delta = 10^{-3}$, max relative residual $\epsilon_{\text{res}}^{\max}$.}			
		\end{subfigure}~~
    		\begin{subfigure}[b]{0.24\textwidth}
			\centering
   \hspace*{-1ex}
			\includegraphics[trim={5cm 17.3cm 9.5cm 4.5cm},clip,width=1.15\textwidth]{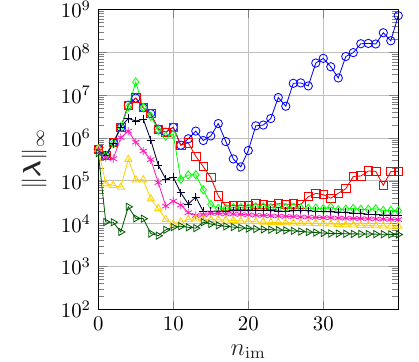}
			\caption{$\delta = 10^{-3}$, max coefficient magnitude $\|\vec\lambda\|_{\infty}$.}			
		\end{subfigure}
  
		% \begin{subfigure}[b]{0.42\textwidth}
		% 	\centering
		% 	\includegraphics[trim={3.0cm 18cm 9cm 2cm},clip,width=\textwidth]{figures/Truncate_two_spheres_gap1e-2.pdf}
		% 	%			\caption{Flow field in the plane of the triangle.}
		% 	%			\label{flow_triangle}
		% 	%		\end{subfigure}~~
		% 	%		\begin{subfigure}[t]{0.33\textwidth}
		% 	%		\centering
		% 	%		\includegraphics[trim={2.9cm 1cm 0cm 0cm},clip,width=\textwidth]{../figures/images/triangle_zoom.eps}
		% 	%		\caption{Zoomed in version of the flow-field in (a).}
		% 	%		\label{flow_triangle2}
		% 	\caption{$\delta = 10^{-2}$}			
		% \end{subfigure}
		% \begin{subfigure}[b]{0.57\textwidth}
		% 	\centering
		% 	\includegraphics[trim={3.0cm 18cm 5.3cm 2cm},clip,width=1\textwidth]{figures/Truncate_two_spheres_gap1e-3.pdf}
		% 	%			\caption{Flow field in the plane of the triangle.}
		% 	%			\label{flow_triangle}
		% 	%		\end{subfigure}~~
		% 	%		\begin{subfigure}[t]{0.33\textwidth}
		% 	%		\centering
		% 	%		\includegraphics[trim={2.9cm 1cm 0cm 0cm},clip,width=\textwidth]{../figures/images/triangle_zoom.eps}
		% 	%		\caption{Zoomed in version of the flow-field in (a).}
		% 	%		\label{flow_triangle2}
		% 	\caption{$\delta = 10^{-3}$}
		% \end{subfigure}						
		\caption{Effect of truncation of singular values in the self-interaction blocks of the target-from-source matrix, investigated for two spheres separated by $\delta$. Notably, $\epsilon_{\text{res}}^{\max}$ decreases  and $\|\vec\lambda\|_{\infty}$ stays bounded with the number of image points per particle $n_{\text{im}}$, only if the smallest singular values are truncated. However, high truncation of singular values come at a cost, as information then is lost. A reasonable level of truncation is chosen to be $\epsilon_{\text{trunc}} = 5\cdot 10^{-12}$. }
		\label{Tol_dep}
	\end{figure}

  For a fixed no-slip boundary condition on two spheres $\delta$ apart, the maximum relative residual, $\epsilon_{\text{res}}^{\max}$, as defined by \eqref{residual}, is computed as a function of the number of image points, $n_{\text{im}}$, for a few different choices of $\epsilon_{\text{trunc}}$. The result is visualized in Figure \ref{Tol_dep} for $\delta = 10^{-2}$ and $\delta = 10^{-3}$; high truncation generally results in larger residuals but smaller coefficient magnitudes. Given experiments of this type,  we pick $\epsilon_{\text{trunc}}= 5\cdot 10^{-12}$ for the numerical examples presented in the results section of the  paper, Section \ref{s:close}.  

 %%%%%%%%%%%%%%%%%%%%%%%%%
 \subsection{Decay of the singular value spectrum}\label{decay}
\begin{figure}[h!]
		\centering
		\hspace*{-8ex}
			\begin{subfigure}[t]{0.33\textwidth}
			\includegraphics[trim={3.1cm 15.6cm 8.5cm 2.5cm},clip,width=1.17\textwidth]{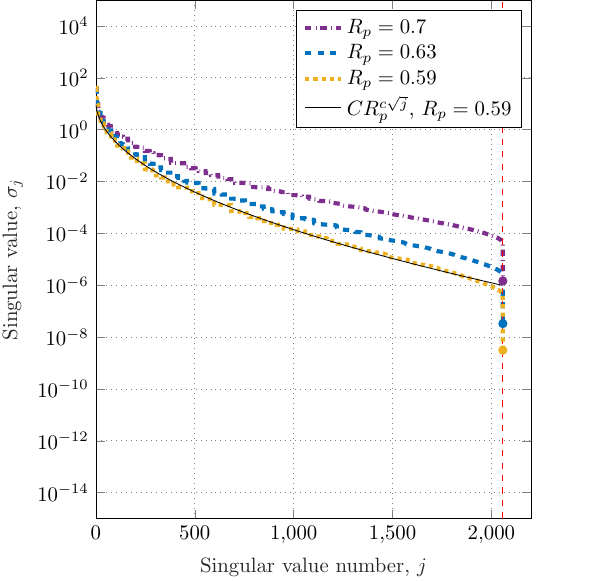}	
			\caption{Decay of singular values for a single particle without image points.}
   \label{sing_without}
		\end{subfigure}~~
  \hspace*{0.5ex}
		\begin{subfigure}[t]{0.33\textwidth}
		\includegraphics[trim={3.1cm 15.6cm 8.5cm 2.5cm},clip,width=1.17\textwidth]{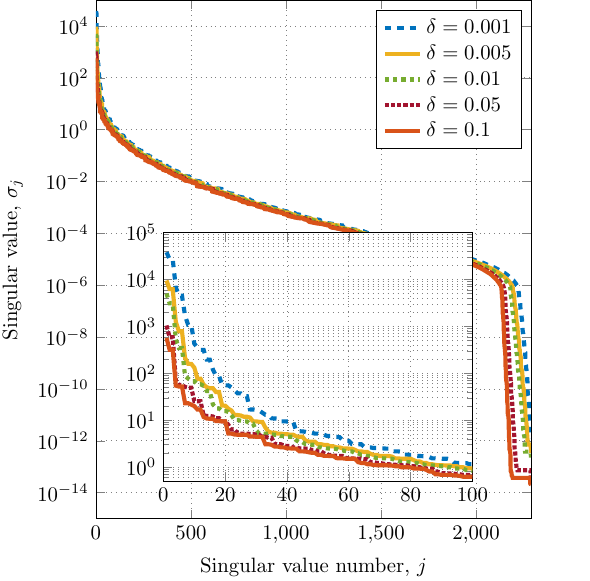}		
\caption{With image points: Example without preconditioning from the left.}
\label{sing_decay}	
		\end{subfigure}~~
  \hspace*{1ex}
				\begin{subfigure}[t]{0.33\textwidth}
		\includegraphics[trim={3.1cm 15.6cm 8.5cm 2.5cm},clip,width=1.17\textwidth]{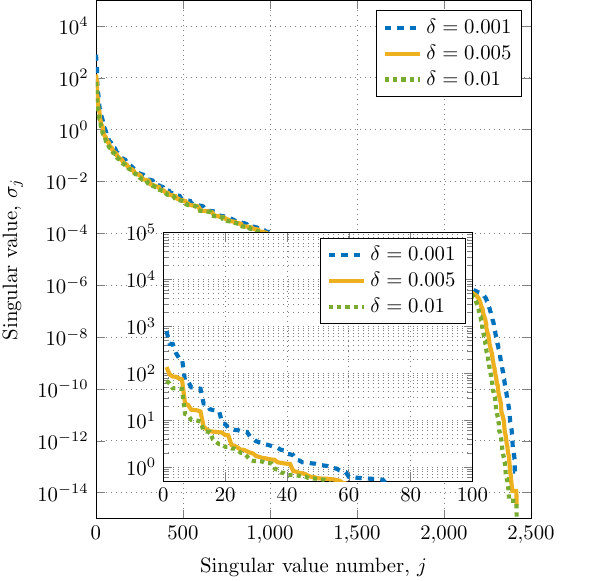}		
		\caption{With image points: Example with preconditioning from the left.}
		%A 200 \times 200 grid is used to compute the flow field. 
  \label{sing_decayb}	
		\end{subfigure}	
		\caption{Decay of singular values for one particle with a single close neighbor at a distance $\delta$. In panel (a), without images, the radius of the proxy-surface, $R_p$, has an impact on the magnitude of the singular values; see also the discussion in Section \ref{base_disc}. %and the small last singular value reflects the null space of the sum of Stokeslets for a sphere (see Remark \ref{nullspace}). 
  With image points in panels (b) and (c), the magnitude of the largest singular values are increased, the smaller %interparticle distance 
  $\delta$ is.  These singular values are in panel (c) decreased by preconditioning from the left.  }
		
	\end{figure}

 Fixing $N=686$ proxy-sphere points and no added image points, Figure \ref{sing_without} plots the
  singular value spectrum for the target-from-source (self-interaction) matrix denoted $\vec B$ in Section \ref{solving}. The $j$th singular value is well-predicted as $\mathcal O(R_p^{c\sqrt{j}})$; see Remark \ref{r:cond}. Hence, the ill-conditioning is more severe, the smaller the radius of the proxy-surface, $R_p$. Then adding images, singular values for one particle in a pair separated by $\delta$ are displayed in Figures \ref{sing_decayb} and \ref{sing_decay} without and with left-preconditioning, as described in Remark \ref{left-precond} of Section \ref{s:close}. The more closely interacting the particles are, the more ill-conditioned is the self-interaction matrix. Left-preconditioning does not remove the numerical nullspace, but does decrease the largest few singular values by 1--2 orders of magnitude. 
	
% bib comes after Appendix...	
%	\bibliographystyle{myIEEEtran} %including doi
    
   % \bibliographystyle{elsarticle-num}    % elsarticle recommends.
 \bibliographystyle{IEEEtranDOIfix} %including doi
	\bibliography{close_refs}
	
\end{document}